\documentclass[a4paper,10pt]{amsart}
\usepackage{amsmath}
\usepackage{amsfonts}
\usepackage{amssymb}
\usepackage{amsthm}
\usepackage[all]{xy}
\usepackage[CJKbookmarks]{hyperref}
\usepackage{fancyhdr}
\usepackage{tikz-cd}
\usepackage{longtable}
\usepackage{array}
\usepackage{diagbox}

\makeatletter
\newcommand{\Rmnum}[1]{\expandafter\@slowromancap\romannumeral #1@}
\makeatother
\newtheorem{theorem}{Theorem}[section]
\newtheorem{lemma}[theorem]{Lemma}
\newtheorem{proposition}[theorem]{Proposition}
\newtheorem{corollary}[theorem]{Corollary}

\theoremstyle{definition}
\newtheorem{definition}[theorem]{Definition}

\newtheorem{hypothesis}[theorem]{Hypothesis}

\newtheorem{remark}[theorem]{Remark}

\numberwithin{equation}{section}
\AddToHook{cmd/theequation/before}{E}

\DeclareMathOperator{\Aut}{Aut}

\DeclareMathOperator{\cd}{cd}

\DeclareMathOperator{\depth}{depth}

\DeclareMathOperator{\gEnd}{\underline{End}}
\DeclareMathOperator{\Ext}{Ext}
\DeclareMathOperator{\gExt}{\underline{Ext}}
\DeclareMathOperator{\fd}{fd}

\DeclareMathOperator{\gr}{gr}

\DeclareMathOperator{\gldim}{gldim}

\DeclareMathOperator{\Gr}{Gr}

\DeclareMathOperator{\Hom}{Hom}
\DeclareMathOperator{\gHom}{\underline{Hom}}

\DeclareMathOperator{\injdim}{injdim}

\DeclareMathOperator{\im}{Im}

\DeclareMathOperator{\pdim}{pdim}

\DeclareMathOperator{\RHom}{RHom}

\DeclareMathOperator{\soc}{soc}

\DeclareMathOperator{\Tor}{Tor}

\DeclareMathOperator{\CMreg}{CMreg}
\DeclareMathOperator{\cmreg}{cmreg}
\DeclareMathOperator{\Exreg}{Ex-reg}
\DeclareMathOperator{\exreg}{ex-reg}
\DeclareMathOperator{\extreg}{extreg}
\DeclareMathOperator{\torreg}{torreg}
\DeclareMathOperator{\Torreg}{Torreg}
\DeclareMathOperator{\Extreg}{Extreg}
\DeclareMathOperator{\reg}{reg}
\DeclareMathOperator{\sdeg}{sup.deg}
\DeclareMathOperator{\ideg}{inf.deg}
\DeclareMathOperator{\ASreg}{ASreg}
\DeclareMathOperator{\asreg}{asreg}

\DeclareMathOperator{\D}{D}
\DeclareMathOperator{\Sq}{Sq}
\newcommand\uxm[2]{u^{#1}{({#2})}}
\newcommand\lxm[2]{l^{#1}{({#2})}}
\newcommand{\abs}[1]{\vert #1\vert}
\newcommand{\qqq}{\mathfrak{q}}

\usepackage{makecell}

\begin{document}
\title{Numerical Homological Regularities over positively graded algebras}

\author{Quanshui Wu}
\address{School of Mathematical Sciences \\
Fudan University \\
Shanghai, 200433 \\
 China}
\email{qswu@fudan.edu.cn}

\author{Bojuan Yi}
\address{School of Mathematical Sciences, Fudan University, Shanghai 200433, China}
\email{21110180026@m.fudan.edu.cn}

\thanks{This research has been supported by the NSFC (Grant No. 12471032) and the National Key Research and Development Program of China (Grant No. 2020YFA0713200).}



\dedicatory{}

\begin{abstract}
We study numerical regularities for complexes over noncommutative noetherian locally finite 
$\mathbb{N}$-graded algebras $A$ such as CM-regularity, Tor-regularity (Ext-regularity) and Ex-regularity, which are the supremum or infimum degrees of some associated canonical complexes. We
introduce their companions---lowercase named regularities, which are defined by taking the infimum or supremum degrees of the respective canonical associated complexes.   
We show that for any right bounded complex $X$ with finitely generated cohomologies, the supremum degree of $R\gHom_A(X, A_0)$ coincides with the opposite of the infimum degree of $X$ if $A_0$ is semisimple. If $A$ has a balanced dualizing complex and $A_0$ is semisimple, we prove that the CM-regularity of $X$ coincides with the supremum degree of $R\gHom_A(A_0,X)$ for any left bounded complex $X$ with finitely generated cohomologies. 

Several inequalities concerning the numerical regularities and the supremum or infimum degrees of derived Hom or derived tensor complexes are given for noncommutative noetherian locally finite $\mathbb{N}$-graded algebras. Some of these are generalizations of J\o rgensen's results on the inequalities between the CM-regularity and Tor-regularity, some are new even in the connected graded case. Conditions are given under which the inequalities become  equalities by establishing two technical lemmas. 

Following Kirkman, Won and Zhang, we also use the numerical AS-regularity (resp. little AS-regularity) to study Artin-Schelter regular property (finite-dimensional property) for noetherian $\mathbb{N}$-graded algebras.
We prove that the numerical AS-regularity of $A$ is zero if and only if that $A$ is an $\mathbb{N}$-graded AS-regular algebra under some mild conditions,
which generalizes a result of Dong-Wu and a result of Kirkman-Won-Zhang. 
If $A$ has a balanced dualizing complex and $A_0$ is semisimple, we prove that the little AS-regularity of $A$ is zero if and only if $A$ is finite-dimensional. 
\end{abstract}

\maketitle

\section{Introduction}
 In 1966, Mumford \cite{M} established a vanishing theorem concerning a coherent sheaf $\mathcal{M}$ on $\mathbb{P}^n$, which says that if $\mathcal{M}$ is $m$-regular in the sense of Castelnuovo that $H^i(\mathcal{M}(m-i))=0$ for all $i \geqslant 1$, then $\mathcal{M}$ is $m'$-regular for all $m' > m$. Mumford's vanishing theorem inspired the development of a concept now known as Castelnuovo-Mumford regularity (abbreviated as CM-regularity). In particular, a notion of (CM-)regularity was introduced by Ooishi \cite{Ooi}, and Eisenbud-Goto \cite{EG} for graded modules over commutative graded algebras (especially, polynomial algebras). This notion of regularity is
closely related to the regularity of sheaves, also to the existence of linear free resolutions of sufficiently high truncations
of graded modules, and to the degrees of generators of the syzygies. If $A =k[x_1, x_2, \dots, x_{n+1}]$ is a polynomial algebra, Proj $A = \mathbb{P}^n$, and $M$ is an $A$-module of the form $
\bigoplus_{i \geqslant 0} H^0(\mathbb{P}^n, \mathcal{O}_X (i))$ for some scheme
$X\subset \mathbb{P}^n$, then the regularity of $M$ is the regularity of $X\subset \mathbb{P}^n$ in the sense of Castelnuovo. 
 Tor-regularity is another numerical regularity introduced by Eisenbud-Goto \cite{EG} and Avramov-Eisenbud \cite{AE} which coincides with the CM-regularity for graded modules over polynomial algebras. 
 Eisenbud and Goto \cite{EG} proved that if $A$ is a polynomial algebra with standard grading then $\CMreg(M)=\Torreg(M)$ for any $M \in \gr A$.
 By using the Tor-regularity it was proved in \cite{AE} that sufficiently high truncations of any finitely generated module over a commutative Koszul algebra have linear free resolutions.

Inspired by the preceding studies in a commutative setting, J\o rgensen first defined CM-regularity for graded modules over noncommutative noetherian connected graded algebras \cite{Jo3}. Furthermore, J\o rgensen  defined Castelnuovo-Mumford regularity for sheaves on a non-commutative
projective scheme (see \cite{AZ}) and proved a noncommutative version of Mumford's vanishing theorem.
Eisenbud-Goto's theorems on linear resolutions and syzygies were also generalized to quantum polynomial algebras. 
 
Subsequently, J\o rgensen \cite{Jo4} established two inequalities relating CM-regularity to Tor-regularity, that is, for any finitely generated graded module $M$, 
$$\CMreg(M)\leqslant \Torreg(M) + \CMreg(A)\textrm{ and } \Torreg(M) \leqslant \CMreg(M) + \Torreg(k),$$ 
which spurred numerous intriguing research efforts. The result that  sufficiently high truncations of finitely generated graded modules have linear free resolutions, was proved in \cite{EG} for commutative standard graded Koszul algebras, in \cite{Jo3} for  noncommutative Koszul connected graded AS-regular algebras, and in \cite{Jo4} for noncommutative Koszul algebras with a balanced dualizing complex. 
 
R\"{o}mer \cite{R} studied CM-regularity and Tor-regularity for commutative standard graded algebras later. One result of particular interest is the characterization (see \cite[Theorem 4.1]{R}) that a commutative standard graded algebra is a polynomial algebra if and only if $\CMreg(M)=\Torreg(M)$ holds for any non-zero finitely generated graded $A$-module $M$ (which is a converse of \cite{EG}); if and only if either of the two inequalities relating CM-regularity to Tor-regularity in \cite{Jo4} is always an equality for any finitely generated graded module.
 
Dong and the first named author \cite{DW} of this paper extended R\"{o}mer's result to noncommutative case, and gave a criterion that a noetherian connected graded algebra with a balanced dualizing complex is Koszul AS-regular if and only if  $\CMreg(M)=\Torreg(M)$ holds for all non-zero finitely generated graded A-module $M$; if and only if $\Torreg(k)=0$ and $\CMreg(A)=0$.  Hence the CM-regularity $\CMreg(A)$ can be considered as an invariant that measures how far away $A$ is from being AS-regular for any Koszul noetherian connected graded $k$-algebra $A$ with a balanced dualizing complex. In fact, Koszulity of a connected graded $k$-algebra $A$ can be characterized by the Tor-regularity $\Torreg(k)=0$ of the trivial module $k$. Based on these facts, Kirkman, Won and Zhang \cite{KWZ1} introduced a numerical invariant $ \ASreg(A):=\CMreg(A)+\Torreg(k)$, called the AS-regularity of $A$, for any noetherian connected graded algebra $A$.
They proved that $A$ is AS-regular if and only if $ \ASreg(A) = 0$, which provides a far-reaching generalization of \cite[Theorem 5.4]{DW}. Weighted versions of these regularities were studied in \cite{KWZ2}. On one hand, Kirkman, Won and Zhang \cite{KWZ2} proved the weighted version of \cite[Theorem 2.5, Theorem 2.6]{Jo4} and \cite[Thoerem 2.8]{KWZ1} which give more relations about these regularities. On the other hand, Kirkman, Won and Zhang \cite{KWZ2} proved the weighted version of \cite[Theorem 3.2 and Theorem 0.8]{KWZ1}, that is, if $A$ is a noetherian connected graded $k$-algebra with a balanced dualizing complex, then $A$ is AS-regular if and only if $\ASreg_{\xi}(A)=0$ for some $\xi \leqslant 1$. 
Meanwhile, Kirkman, Won and Zhang \cite{KWZ3} studied semisimple Hopf algebra actions on AS-regular algebras and showed several upper bounds on the degrees of the minimal generators of the invariant subring, and on the degrees of syzygies of modules over the invariant subring.

In this paper we study numerical regularities over $\mathbb{N}$-graded algebras. Some numerical regularities are used to characterize $\mathbb{N}$-graded AS-regular algebras. All $\mathbb{N}$-graded algebras considered in this paper are assumed to be locally finite  over a field $k$. 

The various numerical regularities appeared in literature, such as CM-regularity, Tor-regularity, Ext-regularity and Ex-regularity (see Definitions \ref{CM-reg-def}, \ref{torreg-def}, \ref{Ext-reg-def}  and \ref{exreg-def}) are defined by using the supremum degrees or infimum degrees of some canonically associated complexes. They are introduced as invariants to study the complexity of complexes.   We also introduce the corresponding lowercase character named regularities of aforementioned regularities.

The {\it supremum degree} and {\it infimum degree} of $X \in \mathbf{D}(\Gr{A})$ are defined respectively as 
$$\sdeg(X)=\sup\{ i + j \mid H^i(X)_j\neq 0,i,j\in \mathbb{Z} \},$$
$$ \ideg(X)=\inf\{ i+ j \mid H^i(X)_j\neq 0, i,j\in \mathbb{Z} \}.$$
\begin{definition} (See Definitions \ref{CM-reg-def}, \ref{torreg-def}, \ref{Ext-reg-def}, \ref{exreg-def})
Let $A$ be an $\mathbb{N}$-graded algebra.
\begin{itemize}
    \item [(1)] For any $X\in \mathbf{D}^+(\Gr A)$, the CM-regularity, cm-regularity, Ex-regularity and ex-regularity of $X$ are defined respectively as
$$\CMreg(X):=\sdeg (R\Gamma_A(X)), \, \cmreg(X):=-\ideg (R\Gamma_A(X)),$$
$$\Exreg(X):=\sdeg (R\gHom_A(S, X)), \, \exreg(X):=-\ideg (R\gHom_A(S, X)).$$

   \item [(2)] For any $X\in \mathbf{D}^-(\Gr A)$, the Tor-regularity, tor-regularity, Ext-regularity and ext-regularity of $X$ are defined respectively as
$$\Torreg(X):=\sdeg (S \,\, {}^L\!\otimes_A X), \, \torreg(X):=-\ideg (S \,\, {}^L\!\otimes_A X),$$
$$\Extreg(X):=-\ideg (R\gHom_A(X,S)), \, \extreg(X):=\sdeg (R\gHom_A(X,S)).$$
\end{itemize}
\end{definition}

Since $S:=A/J$ is a finite-dimensional symmetric algebra, where $J$ is the graded Jacobson radical, 
it follows directly from the definition and the Hom-Tensor adjunction that $\Torreg(X) = \Extreg(X)$ and $\torreg(X) = \extreg(X)$ (see Lemma \ref{Extreg-Torreg}).

The following result is one of the main results, which plays a key role in this paper. In fact, it is a general form of \cite[Lemma 5.2]{DW}.
\begin{theorem} [Lemma \ref{CM and Ext inequality} and Proposition \ref{CM and Ext}] \label{main-result-1}
Let $A$ be a noetherian $\mathbb{N}$-graded algebra. Then, for any $0\neq X \in \mathbf{D}^-(\gr A),$
$$-\ideg(X)\leqslant \extreg(X).$$
If, further, $A_0$ is semisimple, then $-\ideg(X) = \extreg(X).$
\end{theorem}

In the case that $A$ is a commutative standard graded algebra, it was proved in \cite[Theorem 3.9]{Ngu} that the CM-regularity
is the same as the Ex-regularity for any $0\neq X \in \mathbf{D}^b(\gr A)$. By using  Theorem \ref{main-result-1}, we prove that the same conclusion holds if $A$ is a noetherian locally finite $\mathbb{N}$-graded algebra with $A_0$ semisimple and $A$ has a balanced dualizing complex. 

\begin{theorem}[Theorem \ref{ex}]\label{ex 1}
  Suppose $A$ is a noetherian $\mathbb{N}$-graded algebra with $A_0$ semisimple and $A$ has a balanced dualizing complex. Then, for any $0 \neq X \in \mathbf{D}^+_{fg}(\Gr A)$, 
  $\CMreg(X)=\Exreg(X).$
\end{theorem}

Theorem \ref{extreg=-inf} in the following is a general form of \cite[Proposition 1.1]{Tr}, where the algebras $A$ are assumed to be polynomial algebras. The proof here relies on Proposition \ref{CM and Ext}.

\begin{theorem}[Lemma \ref{Exreg-is-greater} and Theorem \ref{Extreg}]\label{extreg=-inf}
Suppose $A$ is a noetherian $\mathbb{N}$-graded algebra. For any  $0\neq X\in \mathbf{D}^b(\gr A)$,  
$\Extreg(X) \geqslant -\ideg(R\gHom_A(X,A)).$

If, further, $A_0$ is semisimple and $\pdim X < \infty$, then $$\Extreg(X)=-\ideg(R\gHom_A(X,A)).$$
\end{theorem}

Table $1$ in \S \ref{relationship between the regularities} provides a synopsis of several fundamental relationships between these regularities when $A_0$ is semisimple.

We prove several inequalities about the homological regularities by using the cohomological spectral sequences induced by $R\gHom_A(X,Y)$ or $Y\, {}^L\!\otimes_A X$ in Propositions  \ref{sdeg-ideg inequality 1}, \ref{sdeg-ideg inequality 1-1}, \ref{sdeg-ideg inequality 1-4}, \ref{sdeg-ideg inequality 2} for noetherian $\mathbb{N}$-graded algebras. For example, 

\begin{theorem} \label{inequalities}Let $A$ be a noetherian $\mathbb{N}$-graded algebra.
 \begin{itemize}  
 \item [(1)] For any $0\neq X\in \mathbf{D}^-(\gr A)$, $0\neq Y\in \mathbf{D}^+(\Gr A)$, 
   $$ -\ideg(R\gHom_A(X,Y))\leqslant \Extreg(X)-\ideg(Y).$$
 \item [(2)] For any $0\neq X\in \mathbf{D}^-(\Gr A)$ with  torsion cohomologies,
        $0\neq Y\in \mathbf{D}^+(\Gr A)$,  
$$\sdeg(R\gHom_A(X,Y))\leqslant -\ideg(X)+\Exreg(Y).$$
\item [(3)] For any $0\neq X\in \mathbf{D}^-(\Gr A)$ 
        and $0\neq Y\in \mathbf{D}^-(\gr A^o)$,
$$-\ideg(R\gHom_A(X,\D(Y)))\leqslant \sdeg(X)+\Extreg(Y).$$
        \item [(4)] For any $0\neq X\in \mathbf{D}^-(\gr A)$, $0\neq Y\in \mathbf{D}^-(\Gr A^{o})$,
   $$    -\ideg(Y\, {}^L\!\otimes_A X) \leqslant \extreg(X)-\ideg(Y).$$
\end{itemize}     
\end{theorem}    

We explore the conditions under which the inequalities in Theorem \ref{inequalities} are equalities by developing two technical results in Lemma \ref{pre for equality 1} and Lemma \ref{pre for equality 2}. 
Suppose $0\neq X\in \mathbf{D}^-(\gr A)$, $p=-\ideg(X)$ and $P^{\bullet}$ is the minimal graded projective resolution of $X$. Lemma \ref{pre for equality 1} says that at least one generator in a minimal generating subset of some $P^{-\alpha}$ is a $(-\alpha)$-cocycle of degree $\alpha - p$.  Lemma \ref{pre for equality 2} says that if $p=-\ideg(X) =-\ideg(e X)$ for some primitive idempotent $e$ in $A_0$, then there is some  $(-\alpha)$-cocycle $y$ of the minimal degree in $P^{-\alpha}$ such that $e y \neq 0$. 
Lemmas \ref{pre for equality 1} and \ref{pre for equality 2} play a key role in proving Propositions \ref{sdeg-ideg equality 1} and \ref{sdeg-ideg equality 2} concerning the identities of the supremum and infimum degrees among some associated canonical complexes, or among numerical homological regularities. 
For example, 
concerning the inequalities in Theorem \ref{inequalities} (1) and (4), we have the following propositions.

\begin{theorem}[Proposition \ref{sdeg-ideg equality 1}]
   Suppose $A$ is a noetherian $\mathbb{N}$-graded algebra with $A_0$ semisimple. For any $0\neq Y \in \mathbf{D}^b(\gr {A^o})$, the following are equivalent.
 \begin{itemize}
     \item [(1)] $\ideg(Y\, {}^L\!\otimes_A X) =  \ideg(X)+\ideg(Y)$ for all $0\neq X\in \mathbf{D}^b(\gr A)$.
     \item [(2)] $\ideg(Ye)=\ideg(Y)$ for any primitive idempotent $e \in A_0$.
 \end{itemize} 
\end{theorem}

\begin{proposition}[Proposition \ref{sdeg-ideg equality 3}]
 Suppose $A$ is a noetherian $\mathbb{N}$-graded algebra with $A_0$ semisimple. For any $0\neq Y \in \mathbf{D}^b(\gr {A})$, the following are equivalent.
 \begin{itemize}
    \item [(1)] $-\ideg(R\gHom_A(X,Y))=\Extreg(X)-\ideg(Y)$ for all $0\neq X\in \mathbf{D}^b(\gr A)$ with $\pdim(X)<\infty$.
    \item [(2)] $\ideg(eY)=\ideg(Y)$ for any primitive idempotent $e \in A_0$.
    \end{itemize}
\end{proposition} 

If $A$ is connected graded, then the second conditions in Propositions \ref{sdeg-ideg equality 1} and  \ref{sdeg-ideg equality 3} are trivially true and so the equalities in (1) hold. We get several identities in  Propositions \ref{sdeg-ideg equality 1}-\ref{sdeg-ideg equality 6} and Corollaries \ref{sdeg-ideg equality 9} and \ref{sdeg-ideg equality 10}, which are new even in the connected graded case.

Following the idea in \cite{KWZ1}, we study the numerical AS-regularities, and introduce CM-regularity homogeneous property and ex-regularity homogeneous property for $\mathbb{N}$-graded algebras in the Section 5. 
The  numerical AS-regularity $\ASreg(A)$  and little numerical AS-regularity $ \asreg(A)$ of a noetherian $\mathbb{N}$-graded algebra $A$ are defined respectively as
$$\ASreg(A)=\CMreg(A)+\Torreg(S), \, \asreg(A)=\cmreg(A)+\torreg(S).$$

We give some inequalities among CM-regularity, Tor-regularity and ex-regularity by using Proposition \ref{sdeg-ideg inequality 1}. The following theorem is a generalization of \cite[Theorems 2.5 and 2.6]{Jo4}.

\begin{theorem}[Theorem \ref{relation between CM and Tor} and Proposition \ref{relation between CM and Tor-p}]\label{CM-Tor-inequalities}
 Suppose $A$ is a noetherian $\mathbb{N}$-graded algebra with a balanced dualizing complex. Then
   \begin{itemize}
       \item [(1)] $\CMreg(X)\leqslant \Torreg(X)+\CMreg(A)$ for all $0\neq X \in \mathbf{D}^-(\gr A)$.
       \item [(2)] $\Torreg(X)\leqslant \CMreg(X)+\Torreg(S)$ for all $0\neq X \in \mathbf{D}^-(\gr A)$.
       \item [(3)] $\ASreg(A)\geqslant 0$.
       \item [(4)] $\CMreg(A)+\exreg({}_AA)\geqslant 0$, and $\CMreg(A)+\exreg(A_A)\geqslant 0$.
       \item [(5)] If $A_0$ is semisimple, then for any $0\neq X \in \mathbf{D}^b(\gr A)$ with $\pdim X< \infty$,
       $$\Torreg(X) \leqslant \CMreg(X)+\exreg({}_AA).$$
\end{itemize}
\end{theorem}
Note that in Theorem \ref{CM-Tor-inequalities} the inequality (4) is stronger than (3). 
By using Theorem \ref{CM-Tor-inequalities}, we prove that sufficiently high truncations of $M \in \gr A$ have linear projective resolutions if $A_0$ is semisimple and $\Torreg(S)<\infty$. 

The following Proposition \ref{linear resolution co} was proved in \cite{EG} for polynomial algebras $A$, in \cite{AE} for  Koszul commutative graded algebras $A$, in \cite[Theorem 3.1]{Jo4} for Koszul connected graded algebras $A$ and in \cite[Theorem 3.13]{KWZ2} for connected graded algebras $A$. 
\begin{proposition}[Proposition \ref{linear resolution}]\label{linear resolution co}
Suppose that $A$ is a noetherian $\mathbb{N}$-graded algebra with a balanced dualizing complex. If $\Torreg(S)\leqslant p$, then for any $0\neq M \in \gr A$ with $\CMreg(M) \leqslant r$,
$\Torreg(M_{\geqslant r}(r+p))\leqslant 0$.

Moreover, if $A_0$ is semisimple, then $\Torreg(M_{\geqslant r}(r+p)) = 0$, and $M_{\geqslant r}(r+p)$ has a linear projective resolution.
\end{proposition}

The fact that the inequality (1) or (2) in Theorem \ref{CM-Tor-inequalities} is always an equality is equivalent to that $\ASreg(A)=0$, and (5) becomes to an equality if and only if $\CMreg(A)+\exreg({}_AA)=0$, 
see Corollaries \ref{ASreg=0} and \ref{CMreg+Torreg=0}.  
The dual versions of Theorem \ref{CM-Tor-inequalities}, as well as Corollaries \ref{ASreg=0} and \ref{CMreg+Torreg=0} are given respectively in Proposition \ref{Extreg+CMreg}, Corollary \ref{ASreg=0 2}-\ref{CMreg+Torreg=0 4} when $A$ has a balanced dualizing complex. 

Similarly, we prove some inequalities between cm-regularity and tor-regularity. 
\begin{theorem}[Proposition \ref{extreg+sdeg}] \label{Prop5.4} 
    Suppose $A$ is a noetherian $\mathbb{N}$-graded algebra with a balanced dualizing complex. Then 
    \begin{itemize}
        \item [(1)] $\extreg(X)\leqslant \cmreg(X)+\extreg(S)$ for all $0\neq X \in \mathbf{D}^-(\gr A)$.
        \item [(2)] $\cmreg(X) \leqslant \extreg(X)+\cmreg(A)$ for all $0\neq X \in \mathbf{D}^-(\gr A)$.
        \item [(3)] $\asreg(A) \geqslant 0$.
    \end{itemize}
    \end{theorem}
    
Then we use $\asreg(A)=0$ to characterize when $A$ is a finite-dimensional algebra.
\begin{proposition}[Corollary \ref{cmreg-finite}]
 Suppose $A$ is a noetherian $\mathbb{N}$-graded algebra with a balanced dualizing complex. If $A_0$ is semisimple, then the following are equivalent. 
\begin{itemize}
        \item [(1)] $\cmreg(A)<\infty$.
        \item [(2)] $\cmreg(X)< \infty$ for all $X \in \mathbf{D}^b(\gr A)$.
        \item [(3)] $A$ is finite-dimensional.
        \item [(4)] $\cmreg(A)=0$.
        \item [(5)] $\asreg(A)=0$.
        \item [(6)] $\asreg(A)< \infty$.
\end{itemize}     
\end{proposition}

We also use little AS-regularity $\asreg(A)=0$ to characterize when the inequalities in Theorem \ref{Prop5.4}  become equalities, see Corollary \ref{extreg+sdeg=0}. 
The dual versions of Proposition \ref{extreg+sdeg} and Corollary \ref{extreg+sdeg=0} are given in  Proposition \ref{exreg+sdeg} and Corollary \ref{exreg+sdeg=0} respectively when $A$ has a balanced dualizing complex. 

For $\mathbb{N}$-graded algebras, we introduce CM-regularity homogeneous property and ex-regularity homogeneous property to study the relations among the previously defined regularities. 
\begin{definition} Suppose $A$ is an $\mathbb{N}$-graded algebra.
 \begin{itemize}
      \item [(1)] $A$ is called {\it  left (resp.\@ right) CM-regularity homogeneous} if 
     $\CMreg(A)=\CMreg(Ae)$  (resp. $\CMreg(A)=\CMreg(eA)$) for any primitive idempotent $e \in A_0$.
     \item [(2)] $A$ is called {\it left (resp.\@ right) ex-regularity homogeneous}  if 
 $\exreg({}_{A}A)=\exreg(Ae)$ (resp. 
 $\exreg(A_A)=\exreg(eA)$) for any primitive idempotent $e \in A_0$.
    \end{itemize}
    \end{definition}

Obviously, any connected graded algebra $A$ is CM-regularity homogeneous and ex-regularity homogeneous. 
If $A$ is an $\mathbb{N}$-graded AS-Gorenstein algebra such that all the Gorenstein parameters of $A$ are equal, then $A$ is left (right) CM-regularity homogeneous and left (right) ex-regularity homogeneous, see Remark \ref{left-right homogeneous}.
    
CM-regularity homogeneous property is used in the following  Theorem \ref{CM-homo-char} to characterize when the inequality (2) in Theorem \ref{CM-Tor-inequalities} becomes an equality, which generalizes \cite[Theorem 0.7]{KWZ1} and \cite[Proposition 5.6]{DW}. 
 \begin{theorem}[Theorem \ref{relation between CM and Tor 2}]\label{CM-homo-char} 
Suppose $A$ is a noetherian $\mathbb{N}$-graded algebra with $A_0$ semisimple and $A$ has a balanced dualizing complex. Then the following  are equivalent.
   \begin{itemize}
       \item [(1)] $A$ is  left CM-regularity homogeneous.
       \item [(2)] $\CMreg(X)=\Torreg(X)+\CMreg(A)$ for all $0\neq X\in \mathbf{D}^b(\gr A)$ with finite projective dimension.
   \end{itemize}
\end{theorem}
Proposition \ref{right CM-reg homogenous} is the dual version of Theorem \ref{CM-homo-char}. 

Similarly, we characterize when $\Torreg(X) \leqslant \CMreg(X)+\exreg({}_AA)$ in Proposition \ref{relation between CM and Tor-p}
and $-\ideg(R\gHom_A(X,A))\leqslant \exreg({}_AA)+\CMreg(X)$ in Proposition \ref{hom-ineq} (2) become equalities.

 \begin{proposition}[Proposition \ref{relation between CM and Tor 3}]\label{111}
    Suppose $A$ is a noetherian $\mathbb{N}$-graded algebra with $A_0$ semisimple and $A$ has a balanced dualizing complex. If $\injdim {}_{A}A=\injdim A_A <\infty$, then the following statements are equivalent.
    \begin{itemize}
        \item [(1)] $A$ is right ex-regularity homogeneous. 
        \item [(2)] $\Extreg(X)=\CMreg(X)+\exreg(A)$ for all $0\neq X\in \mathbf{D}^b(\gr A)$ with $\pdim X < \infty$.
       \item [(3)] $-\ideg(R\gHom_A(X,A))=\CMreg(X)+\exreg(A)$ for all $0\neq X\in \mathbf{D}^b(\gr A)$.
    \end{itemize}
    \end{proposition}
    If further, $\gldim A < \infty$, then
     \begin{corollary}[Corollary \ref {right ex-homo}]\label{222}
       The following are equivalent.
       \begin{itemize}
        \item [(1)] $\Torreg(X)=\CMreg(X)+\Torreg(A_0)$ for all $0\neq X\in \mathbf{D}^b(\gr A)$.
       \item [(2)] $A$ is right ex-regularity homogeneous. 
    \end{itemize}
    \end{corollary}
 The dual versions of Proposition \ref{111} and Corollary \ref {222} are given in Corollaries \ref{the right Tor-regularity homogeneous} and \ref{dual of ex-hom} when $A$ has a balanced dualizing complex.

Following \cite{KWZ1}, we characterize $\mathbb{N}$-graded AS-regular algebras $A$ by  $\ASreg(A)=0$ for noetherian $\mathbb{N}$-graded algebras. We add the condition that $A$ satisfies the Auslander-Buchsbaum formula. The Auslander-Buchsbaum formula holds for connected graded algebras satisfying $\chi$ \cite[Theorem 3.2]{Jo2}. \cite{LW} explored the conditions when a noetherian $\mathbb{N}$-graded algebra satisfying the Auslander-Buchsbaum formula. 

We use the technique in \cite[Theorem 4.7]{IKU} to adjust Gorenstein parameters of $\mathbb{N}$-graded AS-Gorenstein algebras in Theorem \ref{AGP}. The average Gorenstein parameters were introduced also in \cite{IKU}. By using Theorem \ref{AGP}, we prove

 \begin{theorem}[Theorems \ref{AS regular 11} and \ref{AS regular 22}] \label{AS reg}
Suppose $A$ is a ring-indecomposable, basic noetherian $\mathbb{N}$-graded algebra with a balanced dualizing complex $R$, $A_0$ is semisimple and $A$ satisfies the left and right Auslander-Buchsbaum Formula. Then the following are equivalent.
\begin{itemize}
    \item [(1)] $A$ is $\mathbb{N}$-graded AS-regular of dimension $d$ with the average Gorenstein parameter $\ell_{av}^A \in \mathbb{Z}$.
    \item [(2)] For some integers $p_1, p_2,\dots, p_n$,  $B: =\gEnd_A(\mathop{\bigoplus}\limits_{i=1}^n Ae_{i}(p_i))$ is an $\mathbb{N}$-graded CM-algebra of dimension $d$ such that $B_0$ is semisimple and $\ASreg(B)=0$.
    \item[(3)] $B=\gEnd_A(\mathop{\bigoplus}\limits_{i=1}^n Ae_{i}(p_i))$ is an $\mathbb{N}$-graded algebra for some integers $ p_1, p_2,..., p_n$, such that $B_0$ is semisimple and $\ASreg(B)=0$.
\end{itemize}
\end{theorem}


The following theorem is the generalization of \cite[Theorems 3.2 and 0.8]{KWZ1}.

\begin{theorem}[Theorem \ref{AS regular 33}] \label{AS-reg 2}
   Suppose $A$ is a noetherian $\mathbb{N}$-graded algebra with a balanced dualizing complex $R$, $A_0$ is semisimple and $A$ satisfies the left and right Auslander-Buchsbaum Formula. Then the following are equivalent.
   \begin{itemize}
       \item [(1)] $A$ is $\mathbb{N}$-graded AS-regular of dimension $d$
       such that the Gorenstein parameters of $A$ are equal. 
       \item [(2)] $A$ is a graded CM-algebra of dimension $d$ such that $\ASreg(A)=0$.
       \item [(3)] $\ASreg(A)=0$.
   \end{itemize}
    If, furthermore, $A$ is basic, then $\ASreg(A)=0$ if and only if $A$ is AS-regular over $A_0$ in the sense of \cite{MM}.
\end{theorem}

        
        
        
        
\section{Preliminaries}
\subsection{Notations} Let $k$ be a field. An $\mathbb{N}$-graded $k$-algebra $A=\mathop{\bigoplus}\limits_{i=0}^{\infty}A_{i}$ is called {\it locally finite} if $\dim_{k}A_{i}<\infty$ for all $i\geqslant 0$.  An $\mathbb{N}$-graded $k$-algebra $A$ is called {\it connected graded} if $A_0 = k$. 

Throughout the paper, when we say an $\mathbb{N}$-graded algebra it means a locally finite $\mathbb{N}$-graded algebra. For the basic theory about graded algebras we refer \cite{Oys}.
When we say a graded $A$-module, it means a graded left $A$-module. Usually, a graded right $A$-module is viewed as a graded left $A^o$-module, where  $A^o$ is the opposite ring of the $\mathbb{N}$-graded algebra $A$.
A graded $A$-$B$-bimodule is viewed as a graded left $A\otimes B^o$-module for graded algebras $A$ and $B$. Let $A^e = A\otimes A^o$. So, graded $A$-$A$-bimodules are identified as $A^e$-modules.

The category of all graded $A$-modules and degree-zero morphisms is denoted by $\Gr A$. The full subcategory of $\Gr A$ consisting of all finitely generated graded $A$-modules is denoted by $\gr A$. The derived category of $\Gr A$ is denoted by $\mathbf{D}(\Gr{A})$. The right bounded, left bounded and bounded derived category of $\Gr A$ are denoted by $\mathbf{D}^-(\Gr{A})$, $\mathbf{D}^+(\Gr{A})$ and $\mathbf{D}^b(\Gr{A})$ respectively. 
For any $X \in \mathbf{D}(\Gr A)$, $H^i(X)$ denotes the $i$th cohomology of $X$.
The full subcategory of $\mathbf{D}(\Gr{A})$ consisting of complexes with finitely generated (resp. locally finite) cohomologies is denoted by $\mathbf{D}_{fg}(\Gr{A})$ (resp. $\mathbf{D}_{lf}(\Gr{A})$). If $A$ is noetherian, then the bounded derived category of $\gr A$ (resp. category of finite-dimensional modules) is denoted by $\mathbf{D}^b(\gr{A})$ (resp. $\mathbf{D}^b(\fd A)$).
If $A$ is noetherian, then $\mathbf{D}_{fg}^-(\Gr{A}) = \mathbf{D}^-(\gr A)$.

Let $\ell$ be an integer. For any graded $A$-module $M$, the shifted $A$-module $M(\ell)$ is defined by $M(\ell)_m=M_{m+\ell}$,
for all $m\in \mathbb{Z}$. For any cochain complex $X$, we define two kinds of shifting: $X(\ell)$ and $X[\ell]$, where $X(\ell)$ shifts the degrees of each graded module by $(X(\ell))^i=X^i(\ell)$ for all $i\in \mathbb{Z}$, while $X[\ell]$ shifts the complex by $(X[\ell])^i=X^{i+\ell}$ for all $i\in \mathbb{Z}$.
For graded modules $M$ and $N$, we write
\begin{align*}
\gHom_A(M,N)=\mathop{\bigoplus}\limits_j \Hom_{\Gr A}(M,N(j)), \quad
\gExt_A^i(M,N)=\mathop{\bigoplus}\limits_j \Ext_{\Gr A}^i(M,N(j)).
\end{align*}
The derived functors of graded $\gHom$ and $\otimes$ are denoted by $R\gHom$ and ${}^L\!\otimes$ respectively. 
For $X\in \mathbf{D}^-(\Gr A), Y\in \mathbf{D}^+(\Gr A)$, the $i$th cohomology of $R\gHom_A(X,Y)$ is denoted by $\gExt^i_A(X,Y)$, and $$\gExt^i_A(X,Y)\cong H^i(\gHom_A(P^{\bullet},Y))\cong H^i(\gHom_A(X,I^{\bullet}))$$ where $P^{\bullet}$ is a graded projective resolution of $X$ and $I^{\bullet}$ is a graded injective resolution of $Y$.

For $X\in \mathbf{D}^-(\Gr A), Y\in \mathbf{D}^-(\Gr A^o)$, the $(-i)$th cohomology of $Y \, {}^L\!\otimes_A X$ is denoted by $\Tor_i^A(Y,X)$, and $$\Tor_i^A(Y,X)\cong H^{-i}(Y \otimes_A P^{\bullet})\cong H^{-i}(Q^{\bullet} \otimes_A X)$$ where $P^{\bullet}$ and $Q^{\bullet}$ are graded projective resolutions of $X$ and $Y$ respectively.

Let $A$ be an $\mathbb{N}$-graded algebra and $S=A/J$ where $J$ is  the graded Jacobson radical of $A$. A graded projective complex $(P^{\bullet},d^{\bullet})$ is minimal if $\im d^m \subseteq JP^{m+1}$ for all $m \in \mathbb{Z}$. A minimal graded projective resolution of a graded module $M\in \Gr A$ is a graded projective resolution $P^{\bullet}$ of $M$ such that $P^{\bullet}$ is minimal.  If $M \in \Gr A$ is bounded-below, then $M$ has a minimal graded projective resolution \cite{MM}. Furthermore, if $X \in \mathbf{D}^-(\Gr A)$ with $X^i$ is bounded-below, for all $i\in \mathbb{Z}$, then $X$ also has a minimal graded projective resolution.  

If $(P^{\bullet},d^{\bullet})$ is a minimal complex, then the differentials of the complexes $S\otimes_A P^{\bullet}$ and $\gHom_A(P^{\bullet},S)$ are zero.

Let $M$ be a graded $A$-module, $N$ be a graded submodule of $M$. $N$ is called an essential submodule of $M$ if for any graded submodule  $0\neq M'\leqslant M$, $M' \cap N \neq 0$.

A minimal graded injective resolution of $X\in \mathbf{D}^+(\Gr A)$ is a graded injective resolution $(I^{\bullet},d^{\bullet})$ such that $\ker d^m$ is an essential submodule of $I^m$ for all $m\in \mathbb{Z}$. Note that the minimal graded injective resolution of $X\in \mathbf{D}^+(\Gr A)$ always exists.

If $(I^{\bullet},d^{\bullet})$ is a minimal graded injective resolution of $X\in \mathbf{D}^+(\Gr A)$, then the differentials of $\gHom_A(S,I^{\bullet})$ are zero.

The {\it supremum degree} of a graded $A$-module (or more generally, a graded vector space) $M$ is the supremum of the degrees of nonzero homogeneous elements in $M$,
namely, 
$$\sdeg(M)=\sup \{ i \mid M_i \neq 0 \} \in \mathbb{Z} \cup \{ \pm \infty \}.$$ 
Similarly, the {\it infimum degree} of $M$ is the minimum of the degrees of nonzero homogeneous elements in $M$, 
$$ \ideg(M)=\inf \{ i \mid M_i \neq 0 \} \in \mathbb{Z} \cup \{ \pm \infty \}.$$
By convention, $\sup(\emptyset)=-\infty$ and $\inf(\emptyset)=+\infty$.

The {\it supremum degree} of a graded vector space complex $X$ or $X \in \mathbf{D}(\Gr{A})$ is defined as 
$$\sdeg(X)=\sup\{ i + j \mid H^i(X)_j\neq 0,i,j\in \mathbb{Z} \}
=\sup\{i +\sdeg(H^i(X)) \mid  i\in \mathbb{Z} \}.$$

 Similarly, the {\it infimum degree} of a graded vector space complex $X$ or $X \in \mathbf{D}(\Gr{A})$ is defined as
  $$\ideg(X)=\inf\{ i+ j \mid H^i(X)_j\neq 0, i,j\in \mathbb{Z} \}
 =\inf\{i + \ideg(H^i(X)) \mid i\in \mathbb{Z} \}.$$

Let $\D(-):=\gHom_k(-,k)$ be the Matlis dual functor and $(-)^*:=\gHom_A(-,A)$.
The  Matlis dual functor induces an duality between $\mathbf{D}_{lf}(\Gr A)$ to $\mathbf{D}_{lf}(\Gr A^o)$ (see for example \cite[Proposition 3.1]{VdB}).

Obviously, for any $X \in \mathbf{D}(\Gr{A})$,
\begin{align}
\sdeg(X)=-\ideg(\D(X)) \textrm{ and } \ideg(X)=-\sdeg(\D(X)).
\end{align}
Sometimes, we use the notation
\begin{align}
   \sup(X)=\mathop{\sup}\{ i\in \mathbb{Z} \mid  H^i(X)\neq 0 \} \textrm{ and } 
   \inf(X)=\mathop{\inf}\{  i\in \mathbb{Z}\mid  H^i(X)\neq 0 \}.
\end{align}

\subsection{Local Duality and Balanced Dualizing Complexes}
\begin{definition}
    The \emph{socle} of a graded $A$-module $M$ is the largest graded semisimple submodule of $M$, denoted by $\soc M$.
\end{definition}
Obviously, $\soc M \cong \gHom_A(S,M)$ as graded $A$-modules.

For any graded $A$-module $M$, $\Gamma_A(M)= \{ m\in M \mid A_{\geqslant i}m=0, i\gg 0 \}$ is called the {\it torsion submodule} of $M$. A graded module $M$ is called {\it torsion} (resp. {\it torsion-free)} if $\Gamma_A(M)=M$ (resp. $\Gamma_A(M)=0$). Any $x\in \Gamma_A(M)$ is called a torsion element of $M$. As $A$ is locally finite, a finitely generated graded $A$-module is torsion if and only if it is a finite-dimensional module. 

Clearly, $\Gamma_A(M)\cong \mathop{\lim}\limits_{n\to \infty} \gHom_A(A/A_{\geqslant n},M)$, and $\Gamma_A : \Gr A \to \Gr A$ is a left exact functor. Let $R\Gamma_A$ be the right derived functor of $\Gamma_A$.  $R^i\Gamma_A(M):=H^i(R\Gamma_A(M))$ is called the $i$th {\it local cohomology} of $M$. In fact, 
$R^i\Gamma_A(M)\cong \mathop{\lim}\limits_{n\to \infty} \gExt^i_A(A/A_{\geqslant n},M)$ for any $i \geqslant 0$. Similarly, we define $\Gamma_{A^o}$ and $R\Gamma_{A^o}$.

 We say $A$ or $\Gamma_A$  has {\it cohomological dimension} $d$ if for all $M\in \Gr A$ and $i> d$, $R^i\Gamma_A(M)=0$ and there exists some $N \in \Gr A$ such that $R^d\Gamma_A(N)\neq 0$. 

 \begin{lemma}\cite[Lemma 2.9]{LW}
    Let $A$ be a locally finite $\mathbb{N}$-graded $k$-algebra and $M \in \Gr A$. Then $\soc M\neq 0$ if and only if $M$ contains a non-zero torsion element.
\end{lemma}

\begin{lemma}\cite[Lemma 4.1]{LW} \label{inj-decom}
    Suppose $A$ is a noetherian $\mathbb{N}$-graded algebra. Let $I$ be a graded injective module. Then $\Gamma_A(I)$ is the injective hull of $\soc I$ and $I=\Gamma_A(I)\oplus I'$ where $I'$ is a torsion-free graded injective module. 
\end{lemma}
 \begin{definition}
The projective dimension of $X \in \mathbf{D}^-(\Gr A)$ is defined by $$\pdim(X)=\mathop{\sup} \{i\in \mathbb{Z} \mid \gExt_A^i(X,Y) \neq 0 \textrm{ for some }\, Y\in \Gr A\}.$$
\end{definition}

If $P^{\bullet}$ is a minimal graded projective resolution of $X$, then $\pdim(X)=-\mathop{\inf} \{ i\in \mathbb{Z} \mid P^i \neq 0 \}$ (see \cite[Definition 2.1]{DW} in the connected graded case).

\begin{definition}
    The injective dimension of $X \in \mathbf{D}^+(\Gr A)$ is defined by $$\injdim(X)=\mathop{\sup} \{i\in \mathbb{Z} \mid \gExt_A^i(Y,X) \neq 0 \textrm{ for some }\, Y\in \Gr A\}.$$
\end{definition}

If $I^{\bullet}$ is a minimal graded injective resolution of $X$, then $\injdim(X)=\mathop{\sup} \{ i\in \mathbb{Z} \mid I^i \neq 0 \}$ (see \cite[Definition 2.2]{DW} in the connected graded case).

\begin{definition}
Suppose $A$ is an $\mathbb{N}$-graded $k$-algebra and $S=A/J$. The depth of $X \in \mathbf{D}^+(\Gr A)$ is defined to be
$$\depth X =\mathop{\inf} \{ i\in \mathbb{Z}\mid \gExt^i_A (S,X)\neq 0 \}.$$
\end{definition}
The depth of a complex is closely related to the local cohomology as in the connected graded case. 

\begin{lemma}\cite[Lemma 2.12]{LW}
    Suppose $A$ is an $\mathbb{N}$-graded algebra and $X \in \mathbf{D}^+(\Gr A)$. Then $$\depth X=\mathop{\inf} \{ i\in \mathbb{Z} \mid R^i\Gamma_A(X)\neq 0 \}.$$
\end{lemma}
The $\chi$-condition is defined in \cite[Definition 3.2]{AZ}.
\begin{definition}
   Let $A$ be an $\mathbb{N}$-graded algebra. $A$ is called satisfying the $\chi$-condition (resp. the $\chi^{o}$-condition) if for any $M\in \gr A$ (resp. $M\in \gr A^{o}$) and $i\geqslant 0$, $\gExt^i_A(A/A_{\geqslant 1},M)$ (resp. $\gExt^i_{A^{o}}(A/A_{\geqslant 1},M)$) is bounded-above.
\end{definition}

\begin{lemma}\cite{LW}
   Let $A$ be an $\mathbb{N}$-graded algebra. Then the following are equivalent.
   \begin{itemize}
       \item [(1)] $A$ satisfies the $\chi$-condition.
       \item [(2)] $\gExt^i_A(S,M)$ is bounded-above for any $M\in \gr A$ and $i\geqslant 0$.
       \item [(3)] For any graded simple $A$-module $X$, any $M\in \gr A$ and $i\geqslant 0$, $\gExt^i_A(X,M)$ is bounded-above.
   \end{itemize}
\end{lemma}
The following facts due to Ven den Bergh for connected graded algebras are also true for locally finite graded algebras.
\begin{lemma}\cite[Corollary 4.8]{VdB} \label{local-cohom}
  Let $A$ be a noetherian $\mathbb{N}$-graded algebra satisfying the $\chi$-condition and $\chi^{o}$-condition. Then $R\Gamma_A(A)\cong R\Gamma_{A^{\circ}}(A)$ in $\mathbf{D}(\Gr A^e)$. Furthermore, $\depth({}_{A}A)=\depth(A_A)$.
\end{lemma}
\begin{theorem}{(Local Duality)}\label{LD} \cite[Theorem 5.1]{VdB} Suppose $A$ is a noetherian $\mathbb{N}$-graded algebra and $\Gamma_A$ has finite cohomological dimension. Then
\begin{itemize}
    \item [(1)] $\D(R\Gamma_A(A))$ has finite injective dimension as an object in $\mathbf{D}(\Gr A)$.
    \item [(2)] For any graded algebra $B$, and $M\in \mathbf{D}(\Gr A\otimes B^o)$, $$\D(R\Gamma_A(M))\cong R\gHom_A(M,\D(R\Gamma_A(A)))$$ in $\mathbf{D}(\Gr B\otimes A^o)$.
\end{itemize}
\end{theorem}
Dualizing complexes are first defined in the noncommutative case in \cite{Ye}. Here is the definition.
\begin{definition} \label{defi-bdc}(\cite[Definition 3.3]{Ye})
    Let $A$ be a noetherian $\mathbb{N}$-graded algebra. A complex $R \in \mathbf{D}^b(\Gr A^e)$ is called a dualizing complex of $A$, if it satisfies the following conditions:
    \begin{itemize}
        \item [(1)] $R$ has finite injective dimension over $A$ and $A^o$ respectively;
        \item [(2)] The cohomologies of $R$ are finitely generated as $A$-module and $A^o$-module;
        \item [(3)] The natural morphisms $\Phi$: $A \to R\gHom_A(R,R)$ and $\Phi^o$: $A \to R\gHom_{A^o}(R,R)$ are isomorphisms in $\mathbf{D}(\Gr A^e)$.
    \end{itemize}
    Furthermore, if $R\Gamma_A(R) \cong \D(A)$ and $R\Gamma_{A^o}(R) \cong \D(A)$ in  $\mathbf{D}(\Gr A^e)$, then $R$ is called a balanced dualizing complex of $A$.
\end{definition}
If $R$ is a dualizing complex of $A$, then $R\gHom_A(-,R)$ and $R\gHom_{A^o}(-,R)$ define a duality between $\mathbf{D}_{fg}(\Gr A)$ and $\mathbf{D}_{fg}(\Gr A^o)$ (see \cite[Proposition 1.3]{YZ}).

The following existence theorem is also due to Van den Bergh.
\begin{theorem}\cite[Theorem 6.3]{VdB} Let $A$ be a noetherian $\mathbb{N}$-graded algebra. Then $A$ admits a balanced dualizing complex if and only if $A$ satisfies the following two conditions:
\begin{itemize}
    \item [(1)] $A$ satisfies both the $\chi$-condition and $\chi^{o}$-condition;
    \item[(2)] Both $\Gamma_A$ and $\Gamma_{A^o}$ have finite cohomological dimension.
    \end{itemize}
If $A$ admits a balanced dualizing complex $R$, then $R\cong \D(R\Gamma_A(A))$, where $\D$ is the Matlis dual.
\end{theorem}

\subsection{Positively graded AS-Gorenstein algebras}\label{notation-e-i}
\
\newline 
\indent Let $A$ be an $\mathbb{N}$-graded algebra, and $S=A/J$ where $J$ is the graded Jacobson radical. Here we follow the notations in \cite[Section 3]{LW}. Since $S$ is a finite-dimensional semisimple algebra, we may assume that $S\cong \mathop{\bigoplus}\limits_{i=1}^n M_{r_i}(D_i)$ where the $D_i$'s are division algebras. For any $1\leqslant i\leqslant n$, let $S_i$ be a simple $M_{r_i}(D_i)$-module. Then $\{ S_1,S_2,...,S_n\}$ represents all the isomorphic classes of simple $S$-modules, and all the isomorphic classes of graded simple $A$-modules up to shifting. There is a set of orthogonal primitive idempotents $e_1,e_2,...,e_n$ of $A_0$ (it may happen $e_1+e_2+...+e_n\neq 1$) such that $S_i\cong S\Bar{e_i}$ where $\Bar{e_i}$ is the image of $e_i$ in $A_0/J(A_0)$. Then clearly $Ae_i$ is the graded projective cover of $S_i$, and $e_i A$ is the graded projective cover of $S'_i=\Bar{e_i}S$. If $r_i=1$ for all $1 \leqslant i \leqslant n$, then $A$ is called basic.

\begin{lemma} \label{S-symmetric}
Let $A$ be an $\mathbb{N}$-graded algebra. Then
\begin{itemize}
    \item [(1)] $S$ is a symmetric algebra, that is, $\D(S)\cong S$ as $S^e$-modules;
    \item [(2)] $\D(S_i)\cong S^{\prime}_i$ for any $i$.
\end{itemize}
\end{lemma}
\begin{proof}
  (1) Any finite-dimensional semisimple algebra is symmetric \cite[16F]{Lam}.\\
   (2) By (1), $\D(S_i)=\gHom_k(S_i,k)\cong \gHom_k(S \otimes_S S_i, k) \cong \gHom_S(S_i,S)\cong S_i^{'}$.
\end{proof}
\begin{proposition} \label{inj-hull-S-i}
Let $A$ be an $\mathbb{N}$-graded algebra. Then
\begin{itemize}
    \item [(1)] $\D(A_A)$ is the injective hull of ${}_{A}S$.
    \item [(2)] $\D(e_iA)$ is the injective hull of $S_i$, and $\D(Ae_i)$ is the injective hull of $S^{\prime}_i$.
\end{itemize}
\end{proposition}
\begin{definition} \label{Defi-CM}
    Suppose $M$ is a finitely generated graded left $A$-module. We call $M$ is a graded Cohen-Macaulay module of dimension $s$ if $R^i\Gamma_A(M)=0$ for all $i\neq s$ and $R^s\Gamma_A(M)\neq 0$. We say $A$ is a graded CM-algebra of dimension $s$ if ${}_{A}A$ is a graded Cohen-Macaulay module of dimension $s$.
\end{definition}
$\mathbb{N}$-graded AS-Gorenstein algebras (sometimes called generalized AS-Gorenstein) were defined in several slightly different ways \cite{MV1,MV2,MM,RR,RRZ}. 
The following definition is given in \cite{LW}.
\begin{definition} \label{GAS-Gorenstein} 
   A noetherian $\mathbb{N}$-graded algebra $A$ is called an $\mathbb{N}$-graded Artin-Schelter Gorenstein algebra (or $\mathbb{N}$-graded AS-Gorenstein for short) of dimension $d$ if the following conditions hold.
\begin{itemize}
    \item[(1)] $\injdim({}_{A}A)=\injdim(A_A)=d < \infty$,
    \item[(2)] For every graded simple left $A$-module $M$, $\gExt_A^i(M,A)=0$ if $i\neq d$; for every graded simple right $A$-module $N$, $\gExt_{A^o}^i(N,A)=0$ if $i\neq d$,
    \item[(3)] $\gExt_A^d(-,A)$ induces a bijection between the isomorphism classes of graded simple $A$-modules and the isomorphism classes of graded simple $A^o$-modules, with the inverse $\gExt_{A^o}^d(-,A)$.
\end{itemize}
Furthermore, if $A$ has finite global dimension, then $A$ is called an $\mathbb{N}$-graded Artin-Schelter regular algebra (or $\mathbb{N}$-graded AS-regular for short).
\end{definition}

With the notation as above, if $A$ is an $\mathbb{N}$-graded AS-Gorenstein algebra of dimension $d$, then there is a permutation $\sigma \in \mathfrak{S}_n$ such that $\gExt^d_A(S_i,A)\cong S^{\prime}_{\sigma(i)}(\ell_i)$ for some $\ell_i\in \mathbb{Z}$ ($1\leqslant i\leqslant n$). So
\begin{align}\label{def-GAS-Goren}
   \gExt^d_A(S,A)\cong \mathop{\bigoplus}\limits_{i=1}^n (e_{\sigma(i)}S(\ell_i))^{r_i}.
\end{align}
With the notation as above, $\{ \ell_1,\ell_2,...,\ell_n\}$ is called the set of Gorenstein parameters of $A$,
If $A$ is $\mathbb{N}$-graded AS-Gorenstein of dimension $d$, then $A$ is a CM-algebra of dimension $d$. 

\begin{definition}\cite[Definition 4.3]{IKU}
  $\ell_{av}^{A}:=n^{-1}\sum_{i=1}^n \ell_i \in \mathbb{Q}$ is called the average Gorenstein parameter of $A$.
\end{definition}

The proofs of the following lemmas are direct from the definition of $\mathbb{N}$-graded AS-Gorenstein algebras.
\begin{lemma}\cite[Lemma 3.8]{LW}\label{def-GAS}
Suppose $A$ is an $\mathbb{N}$-graded AS-Gorenstein algebra of dimension $d$. With the notation as above,  $\gExt^d_{A^o}(S^{\prime}_i,A)\cong S_{\sigma^{-1}(i)}(\ell_{\sigma^{-1}(i)})$ for any $1\leqslant i\leqslant n$. So
\begin{align}\label{def-GAS-Goren coro}
   \gExt^d_{A^o}(S,A)\cong \mathop{\bigoplus}\limits_{i=1}^n (Se_{\sigma^{-1}(i)}(\ell_{\sigma^{-1}(i)}))^{r_i}.
\end{align}
\end{lemma}
For $\mathbb{N}$-graded algebras, Minamoto and Mori also gave a definition of AS-Gorenstein (AS-regular) property, which is called AS-Gorenstein (AS-regular) over $A_0$ \cite{MM}.  
\begin{definition} \cite[Definition 3.1]{MM}\label{def-AS-Goren over A_0}
    A noetherian $\mathbb{N}$-graded algebra $A$ is called AS-Gorenstein over $A_0$ (resp. AS-regular over $A_0$) of dimension $d$ with Gorenstein parameter $\ell$, if the following statements hold.
    \begin{itemize}
        \item [(1)] $\injdim{{}_AA}=d$ (resp. $\gldim A =d $);
        \item [(2)] there is an algebra automorphism $\nu$ over $A_0$ such that 
        $$R\gHom_A(A_0,A)\cong {}^1\D(A_0)^{\nu}(\ell)[-d]$$
        in $\mathbf{D}(\Gr A^e)$.
    \end{itemize}
\end{definition}

Obviously, if $A$ is AS-Gorenstein of dimension $d$ over $A_0$, then $A$ is $\mathbb{N}$-graded AS-Gorenstein of dimension $d$.

\begin{lemma}\cite[Lemma 2.8]{IKU}, \cite[Proposition 3.16]{LW}\label{graded Morita equi}
    Generalized AS-Gorenstein (regular) algebras of dimension $d$ are closed under graded Morita equivalences. 
\end{lemma}
\begin{lemma}\cite[Section 4.2]{IKU}\label{GP}
   Let $A$ be a ring-indecomposable, basic $\mathbb{N}$-graded AS-Gorenstein algebra of dimension $d$ with Gorenstein parameters $\{\ell_1^A,\ell_2^A,...,\ell_n^A\}$. Given integers $p_i$ ($1\leqslant i\leqslant n$), let $B=\gEnd_A(\mathop{\bigoplus}\limits_{i=1}^n Ae_{i}(p_i))$. Then $\ell_i^B=\ell_i^A-p_i+p_{\sigma(i)}$ for any $1\leqslant i\leqslant n$. In particular, $\ell_{av}^A=\ell_{av}^B$. 
\end{lemma}
Gorenstein parameters of $A$ can be adjusted to be almost identical under graded Morita equivalences \cite{IKU}. 

\begin{theorem}\cite[Theorem 4.7]{IKU}\label{AGP 1}
   Let $A$ be a ring-indecomposable basic $\mathbb{N}$-graded AS-Gorenstein algebra of dimension $d$ with Gorenstein parameters $\{\ell_1^A,\ell_2^A,...,\ell_n^A\}$. Then there exists a ring-indecomposable basic $\mathbb{N}$-graded AS-Gorenstein algebra $B$ of dimension $d$ such that the following holds.
   \begin{itemize}
       \item [(1)] $B$ is graded Morita equivalent to $A$.
       \item [(2)] $|\ell_i^B-\ell_{av}^B| < 1$ holds for each $1\leqslant i\leqslant n$.
   \end{itemize}
As a consequence, if $\ell_{av}^A \in \mathbb{Z}$, then $\ell_i^B=\ell_{av}^B$ for each $i$. 
\end{theorem}

\begin{theorem}\cite[Theorem 4.12]{LW}\label{Bdc}
    If $A$ is an $\mathbb{N}$-graded AS-Gorenstein algebra of dimension $d$, then $A$ admits a balanced dualizing complex given by $\D(R^d\Gamma_A(A))[d]$.

    In fact, the balanced dualizing complex $R$ of $A$ has the form 
    \begin{align}   \label{GAS-Goren} 
     R=\D(R^d\Gamma_A(A))[d]\cong \mathop{\bigoplus}\limits_{i=1}^n (Ae_{\sigma(i)}(-\ell_i))^{r_i}[d]\cong \mathop{\bigoplus}\limits_{i=1}^n (e_iA(-\ell_i))^{r_{\sigma(i)}}[d]
     \end{align}
     where  $\{\ell_1,\ell_2,\cdots,\ell_n\}$ is the set of Gorenstein parameters of $A$.
\end{theorem} 

The following fact indicates that if $A$ satisfies the $\chi$-condition and $\chi^{o}$-condition, then the condition (3) in Definition \ref{GAS-Gorenstein} could be removed.
\begin{proposition}
    Let $A$ be a noetherian $\mathbb{N}$-graded algebra satisfying the $\chi$-condition and $\chi^{o}$-condition. If the following conditions hold:
    \begin{itemize}
        \item [(1)] $\injdim({}_{A}A)=\injdim(A_A)=d < \infty$,
    \item[(2)] For every graded simple left $A$-module $M$, $\gExt_A^i(M,A)=0$ if $i\neq d$;
    For every graded simple right $A$-module $N$, $\gExt_{A^o}^i(N,A)=0$ if $i\neq d$.
    \end{itemize}
    Then $A$ is an $\mathbb{N}$-graded AS-Gorenstein algebra.
\end{proposition}
\begin{proof}
    For every graded simple $A$-module $M$, $\gExt^i_{A^o}(\gExt^j_A(M,A),A)=0$ for $(i,j)\neq (d,d)$ and $\gExt^d_{A^o}(\gExt^d_A(M,A),A)\cong M$ by \cite[Theorem 1]{MV1}. We claim that for every graded simple $A$-module $M$, $\gExt^d_A(M,A)$ is a graded simple $A^o$-module.

       Since $A$ satisfies the $\chi$-condition, $\gExt^d_A(M,A)$ is a finite-dimensional module. Let $L=\gExt^d_A(M,A)$, and $M^{\prime}$ be a graded simple $A^o$-submodule of $L$. If $L_A$ is not graded simple. Then the exact sequence
       $$0\to M^{\prime} \to L \to L^{\prime} \to 0$$
   with $ L^{\prime} \neq 0$ induces an exact sequence of $A$-modules
  $$
       0\to \gExt^d_{A^o}(L^{\prime},A) \to M \to \gExt^d_{A^o}(M^{\prime},A) \to 0.
 $$
 It follows that either $\gExt^d_{A^o}(L^{\prime},A)=0$ or $\gExt^d_{A^o}(L^{\prime},A)\cong M$.
 
 Suppose $\gExt^d_{A^o}(L^{\prime},A)\cong M$. Then $\gExt^d_{A^o}(M^{\prime},A)=0$, which is a contradiction.
 
  Suppose $\gExt^d_{A^o}(L^{\prime},A)=0$. Let $M^{\prime \prime}$ be a graded simple $A^o$-submodule of $L^{\prime}$. Then the exact sequence
 $$
    0\to M^{\prime \prime} \to L^{\prime} \to L^{\prime \prime} \to 0,
 $$
 induces an exact sequence of $A$-modules
 $$
     0\to \gExt^d_{A^o}(L^{\prime \prime},A) \to 0 \to \gExt^d_{A^o}(M^{\prime \prime},A) \to 0.
$$
Hence $\gExt^d_{A^o}(M^{\prime \prime},A)=0$, which is also a contradiction. 
 
Therefore $\gExt^d_A(M,A)$ is a graded simple $A^o$-module.

 Similarly, $\gExt^d_{A^o}(N,A)$ is a graded simple $A$-module for any graded simple $A^o$-module $N$. Hence $A$ is an $\mathbb{N}$-graded AS-Gorenstein algebra.
\end{proof}

\section{homological regularities}

In this section, we study various numerical regularities for $\mathbb{N}$-graded algebras, including  Castelnuovo-Mumford regularities, Tor-regularities (Ext-regularities), Ex-regularities, and Artin-Schelter regularities, and their companions---corresponding lowercase character named regularities. We explore their interrelationship, with particular emphasis on the connections between these numerical regularities. Some basic facts are collected in Table $1$.

\subsection{Definitions of various numerical regularities} In this subsection
we recall the definitions of various aforementioned regularities, and define the corresponding lowercase characters named regularities. A weighted version of these numerical regularities and relevant results over $\mathbb{N}$-graded algebras are discussed in a separate paper \cite{WY}.

Now let us recall the definition of Castelnuovo-Mumford regularity for noncommutative $\mathbb{N}$-graded algebras. 
\begin{definition}\label{CM-reg-def}
    Let $A$ be an $\mathbb{N}$-graded algebra. 
    \begin{itemize}
        \item [(1)] \cite{Jo3,Jo4,DW,KWZ1} The Castelnuovo-Mumford regularity ({\it CM-regularity} for short) of $X\in \mathbf{D}^+(\Gr A)$ is defined to be
$$\CMreg(X)=\sdeg (R\Gamma_A(X)) =\sup\{ i + j  \mid R^i\Gamma_A(X)_j\neq 0, i,j\in \mathbb{Z} \}. $$
        \item [(2)] The {\it cm-regularity} of $X\in \mathbf{D}^+(\Gr A)$ is defined to be
$$\cmreg(X)=-\ideg (R\Gamma_A(X)) =-\inf\{ i + j  \mid R^i\Gamma_A(X)_j\neq 0, i,j\in \mathbb{Z} \}. $$
    \end{itemize}
\end{definition}

\begin{lemma}\label{CM-reg-finite}
Suppose $A$ is a noetherian $\mathbb{N}$-graded algebra with a balanced dualizing complex. Then
\begin{itemize}
    \item [(1)]  $\CMreg(X) < \infty$ for all $0\neq X \in \mathbf{D}^b(\gr A)$.
  \item[(2)] $\CMreg({}_{A}A)=\CMreg(A_A)$, which is denoted by $\CMreg(A)$ later in the paper \cite[Observation 2.3]{Jo4}.
\end{itemize}
\end{lemma}
\begin{proof}
    (1) Let $R$ be a balanced dualizing complex of $A$. Then $R\gHom_A(-,R)$ and $R\gHom_{A^o}(-,R)$ define a duality between $\mathbf{D}^b(\gr A)$ and $\mathbf{D}^b(\gr A^o)$. Hence $0\neq X\in \mathbf{D}^b(\gr A)$ if and only if $0\neq R\gHom_A(X,R)\in \mathbf{D}^b(\gr A^o)$. By the local duality theorem  (see Theorem \ref{LD}), $0\neq R\gHom_A(X,R)\cong \D(R\Gamma_A(X))\in \mathbf{D}^b(\gr A^o)$. It follows that $$\CMreg(X) = \sdeg(R\Gamma_A(X)) = -\ideg \big(\D(R\Gamma_A(X)) \big) \neq \pm \infty.$$
    (2) Since $A$ has a balanced dualizing complex $R$, by Lemma \ref{local-cohom}, $R\Gamma_A(A)\cong R\Gamma_{A^o}(A)$. Hence $\CMreg({}_{A}A)=\CMreg(A_A)$. 
\end{proof}

In fact, if $A$ does not have balanced dualizing complex, $\CMreg(A)$ may be infinite \cite[Example 5.1]{KWZ1}.
It follows from \cite[Lemma 5.6]{KWZ1} that $\CMreg(A)$ can take any integers. Unlike $\CMreg(X)$, $\cmreg(X)$  may be $\infty$ in many cases.  If $A_0$ is semisimple and $A$ has a balanced dualizing complex, then by Corollary \ref{cmreg-finite}, either $\cmreg(A)=0$ or $\cmreg(A)=+\infty$.

\begin{definition}  \label{torreg-def} 
Let $A$ be an $\mathbb{N}$-graded algebra and $S=A/J$. 
\begin{itemize}
    \item [(1)] \cite{Jo3, Jo4, DW, KWZ1} The {\it Tor-regularity} of $X\in \mathbf{D}^-(\Gr A)$ is defined to be 
$$\Torreg(X)=\sdeg (S \,\, {}^L\!\otimes_A X).$$
    \item [(2)] The {\it tor-regularity} of $X\in \mathbf{D}^-(\Gr A)$ is defined to be 
$$\torreg(X)=-\ideg (S \,\, {}^L\!\otimes_A X)$$
\end{itemize}
\end{definition}

By definition, $\Torreg(X)=\sup \{ -i + j \mid \Tor^A_i(S,X)_j\neq 0, i,j\in \mathbb{Z} \},$ and $\torreg(X)=-\inf \{ -i + j \mid \Tor^A_i(S,X)_j\neq 0, i,j\in \mathbb{Z} \}$.

 If $A$ is an $\mathbb{N}$-graded algebra with $A_0$ semisimple, then $\Torreg({}_A A_0) \geqslant 0$ and $\torreg({}_A A_0)=0$.
 
\begin{definition} \cite[Definition 2.14.1]{BGS}
Let $A$ be an $\mathbb{N}$-graded ring with $A_0$ semisimple. Suppose $M$ is a graded $A$-module. We say $M$ has a linear projective resolution if $M$ has a graded projective resolution $P^{\bullet}$ such that $P^{-i}$ is generated by elements of degree $i$, 
 or equivalently,
 $\Tor^A_i(A_0,M)_j=0$ for all $j\neq i$. If $M$ has a linear projective resolution, then $M$ is called a Koszul module. If $A_0$ viewed as a graded $A$-module is a Koszul module, then $A$ is called a Koszul ring.
\end{definition}

 So, if $A_0$ is semisimple, then $\Torreg({}_A A_0)= 0$ if and only if $A$ is a Koszul ring, and $\Torreg({}_A A_0)$ can be regarded as a measure of how far that $A$ is from being a Koszul ring. Definition \ref{torreg-def} has a right version. 
 In general, $\Torreg({}_A S)=\Torreg(S_A)\geqslant 0$ and  $\torreg({}_A S)=\torreg(S_A)\geqslant 0$ by the definition, which are denoted by $\Torreg(S)$ and $\torreg(S)$ respectively. In general, that is, if $A_0$ is not necessarily semisimple, then for any $M \in \Gr A$ with $\ideg(M)< \infty$,  $\Torreg(M)=\torreg(M)=0$ if and only if $M$ has a linear projective resolution. 

Note that $\Torreg(S)$ can take all positive integers by taking different graded algebras \cite[Lemma 5.6]{KWZ1}. If $A_0$ is not semisimple, then $\torreg(S)$ can also take any positive integers. For example, let $A=A_0$ be not semisimple with $\gldim A=1$. Then $\torreg(S_A)=1$. For any given integer $n$, let $B=A^{\otimes{n}}$, then $\torreg(S_B)=n$.
 
The Ext-regularity in the noncommutative case is first introduced in \cite{Jo4}. Here is the definition.
\begin{definition} \label{Ext-reg-def}
Let $A$ be an $\mathbb{N}$-graded algebra and $S=A/J$. 
\begin{itemize}
    \item [(1)] \cite{Jo4, DW, KWZ1} The {\it Ext-regularity} of $X\in \mathbf{D}^-(\Gr A)$ is defined to be $$\Extreg(X)=-\ideg (R\gHom_A(X, S)).$$ 
   \item [(2)] The {\it ext-regularity} of $X\in \mathbf{D}^-(\Gr A)$ is defined to be
   $$\extreg(X)=\sdeg (R\gHom_A(X, S)).$$
\end{itemize}
\end{definition}

In fact, for any $X\in \mathbf{D}^-(\Gr A)$, the Tor-regularity  of $X$ is the same as the Ext-regularity of $X$, and the tor-regularity of $X$ is the same as the ext-regularity of $X$ as claimed in the following lemma (see also \cite[Remark 4.5]{DW}). 

\begin{lemma}\label{Extreg-Torreg}
Suppose $A$ is an $\mathbb{N}$-graded algebra and $X\in \mathbf{D}^-(\Gr A)$. Then $$\Extreg(X)=\Torreg(X), \, \textrm{ and }\, \extreg(X)=\torreg(X).$$
\end{lemma}
\begin{proof} 
    By Lemma \ref{S-symmetric}, $\D(S)\cong S$ as $S$-$S$ as bimodules. Hence 
     $$\D(S \, {}^L\!\otimes_A X) \cong R\gHom_A(X, \D(S)) \cong R\gHom_A(X, S).$$
    The conclusion follows.
\end{proof}
 
There are more concrete descriptions of $\Torreg(X)\, (=\Extreg(X))$ and $\torreg(X) \,\\(=\extreg(X))$ by using the generating degrees of the minimal graded projective resolution of $0\neq X \in \mathbf{D}^-(\gr A)$. Suppose $P^{\bullet}$ is a minimal graded projective resolution of $X$. If $P^{-m}\neq 0$, we may assume that $P^{-m}$ has the following form (see \S\ref{notation-e-i})
$$P^{-m}= \mathop{\bigoplus}\limits_i \big(\mathop{\bigoplus}\limits_j Ae_i(-s_m^{i,j})\big).$$
We fix the notation that 
\begin{align}\label{u-l-degree-of-generators}
\uxm{m}{X}:=\sup \{ s_m^{i,j} \mid i, j\},\, \textrm{ and } \, \lxm{m}{X}:=\inf \{ s_m^{i,j} \mid i, j\}
\end{align}
which are the maximum degree and the minimal degree of the elements in the minimal generating set of $P^{-m}$ respectively. In fact, $\lxm{m}{X}=\ideg(P^{-m}).$

\begin{lemma}\label{notations of X}  Let $A$ be a noetherian $\mathbb{N}$-graded algebra. For any $0\neq X \in \mathbf{D}^-(\gr A)$,
$$\Torreg(X)=\Extreg(X)=\sup\{\uxm{m}{X}-m\mid P^{-m}\neq 0\},$$
$$\torreg(X)=\extreg(X)=-\inf\{\lxm{m}{X}-m \mid P^{-m}\neq 0\}.$$ 
\end{lemma}

\begin{proof} It follows from the minimality of the projective resolution that
$$\Tor_m^A(S,X)=S\otimes_A P^{-m}=\mathop{\bigoplus}\limits_i \big(\mathop{\bigoplus}\limits_j Se_i(-s_m^{i,j})\big),$$ 
$$\gExt_A^m(X,S)=\gHom_A(P^{-m},S)=\mathop{\bigoplus}\limits_i \big(\mathop{\bigoplus}\limits_j e_iS(s_m^{i,j})\big).$$ 
Therefore, $$\Torreg(X)=\sup\{\uxm{m}{X}-m\mid P^{-m}\neq 0\}=\Extreg(X),$$
$$\torreg(X)=-\inf\{\lxm{m}{X}-m \mid P^{-m}\neq 0\}=\extreg(X).$$
\end{proof}

Nguyen \cite{Ngu} introduced another regularity \cite[Definition 3.1]{Ngu}, which is denoted by $\Exreg(X)$ in this paper, for bounded-below complexes $X$ over a commutative 
standard graded $k$-algebra $A$. 
It was proved in \cite[Theorem 3.9]{Ngu} that  $\Exreg(X)$ coincides with $\CMreg(X)$ for any $X\in \mathbf{D}^b(\gr A)$. This is proved to be true in the noncommutative case that $A$ is a noetherian $\mathbb{N}$-graded algebra with $A_0$ semisimple and $A$ has a balanced dualizing complex in Theorem \ref{ex}.
\begin{definition} \label{exreg-def}
    Let $A$ be a $\mathbb{N}$-graded algebra. 
    \begin{itemize}
        \item [(1)]  \cite{Ngu} The {\it Ex-regularity} of $X\in \mathbf{D}^+(\Gr A)$ is defined to be 
    $$\Exreg(X)=\sdeg (R\gHom_A(S, X)).$$ 
        \item [(2)] The {\it ex-regularity} of $X\in \mathbf{D}^+(\Gr A)$ is defined to be 
    $$\exreg(X)=-\ideg (R\gHom_A(S, X)).$$ 
    \end{itemize}
\end{definition}

Note that by Theorem \ref{ex} and \cite[Lemma 5.6]{KWZ1}, $\Exreg(A)$ runs over all integers. 

Suppose $B$ is an AS-regular algebra of dimension 3 with Gorenstein parameter 4 generated in degree $1$. Then $\exreg(B)=1$. Let $C=k[x]/(x^2)$. Then $\exreg(C)=-1$. Given any integer $n$, there exist nonnegative integers $p,q$ such that $n=q-p$. Consider the algebra $A=C^{\otimes{p}}\otimes B^{\otimes q}$. Then $\exreg(A)=n$. Hence $\exreg(A)$ runs over all integers.

Suppose $A$ has a balanced dualizing complex. Then $\D(R\Gamma_A(-))$ gives a duality between $\mathbf{D}_{fg}(\Gr A)$ and $\mathbf{D}_{fg}(\Gr A^o)$. So Ex-regularity and ext-regularity, respectively, ex-regularity and Ext-regularity,
are dual concepts if $A$ has a balanced dualizing complex. 
\begin{lemma}\label{dual-ext-ex} Let $A$ be a noetherian $\mathbb{N}$-graded algebra with a balanced dualizing complex $R$. Then 
\begin{itemize}
 \item [(1)] for any $X\in \mathbf{D}^+_{fg}(\Gr A)$, 
 $$\Exreg(X)=\extreg(\D(R\Gamma_A(X))), \textrm{ and } \exreg(X) = \Extreg(\D(R\Gamma_A(X))),$$
\item [(2)] for any $X\in \mathbf{D}^-_{fg}(\Gr A)$,
$$\extreg(X)=\Exreg(\D(R\Gamma_A(X))), \textrm{ and } \Extreg(X)=\exreg(\D(R\Gamma_A(X))).$$
\end{itemize}

\end{lemma}
\begin{proof}
The conclusion follows from 
$$ R\gHom_A(X,S)\cong R\gHom_{A^o}(S,\D(R\Gamma_A(X))), \textrm{ and }$$
$$ R\gHom_A(S,X)\cong R\gHom_{A^o}(\D(R\Gamma_A(X)),S).$$
\end{proof}

 The Matlis dual gives a duality between $\mathbf{D}_{lf}(\Gr A)$ and $\mathbf{D}_{lf}(\Gr A^o)$, which induces some dual relations between Ex-regularity and ext-regularity (tor-regularity), respectively, ex-regularity and Ext-regularity (Tor-regularity).
\begin{lemma}\label{dual-ext-ex-by0Matlis} 
Let $A$ be an $\mathbb{N}$-graded algebra such that $S$ is finite-dimensional. Then
 \begin{itemize}
 \item [(1)] for any $X\in \mathbf{D}^+_{lf}(\Gr A)$,
\begin{align}\label{Ex-Ext and ex-ext-3}
  \Exreg(X)=\extreg(\D(X)), \textrm{ and } \exreg(X) = \Extreg(\D(X)), 
\end{align}
\item [(2)] for any $X\in \mathbf{D}^-_{lf}(\Gr A)$,
\begin{align}\label{Ex-Ext and ex-ext-4}
  \extreg(X)=\Exreg(\D(X)), \textrm{ and } \Extreg(X)=\exreg(\D(X)).
\end{align}
\end{itemize}
\end{lemma}
\begin{proof}
The conclusion follows from 
$$R\gHom_{A}(S, X) \cong  R\gHom_{A^o}(\D(X),S) \textrm{ for any } X\in \mathbf{D}^+_{lf}(\Gr A), \textrm{ and }$$
$$R\gHom_{A}(X, S) \cong R\gHom_{A^o}(S,\D(X))\textrm{ for any } X\in \mathbf{D}^-_{lf}(\Gr A).$$
\end{proof}

The following lemma shows that $\Exreg(X)\geqslant \CMreg(X)$ for any $0\neq X\in \mathbf{D}^+(\Gr A)$ in general. The equality holds if $A$ has a balanced dualizing complex and $A_0$ is semisimple as showed in Theorem \ref{ex}, which generalizes \cite[Theorem 3.9]{Ngu}.
\begin{lemma}\label{lem-CM-ex}
Let $A$ be a noetherian $\mathbb{N}$-graded algebra. 
For any $0\neq X\in \mathbf{D}^+(\Gr A)$,
$$\Exreg(X)\geqslant \CMreg(X).$$
\end{lemma}
\begin{proof}
    Let $I^{\bullet}$ be a minimal graded injective resolution of $X$. For any $m\in \mathbb{Z}$ such that $\gExt_A^m(S,X)\neq 0$, $\gExt_A^m(S,X)\cong \gHom_A(S,I^m)\cong \soc I^m$, we may assume that (see \S\ref{notation-e-i})
$$\soc I^m \cong \mathop{\bigoplus}\limits_i \big(\mathop{\bigoplus}\limits_j Se_i(-\ell_m^{i,j})\big).$$
    Therefore, $\sdeg(\gExt^m_A(S,X)) = \mathop{\sup}\{\ell_m^{i,j}\mid i,j\}$, and 
    $$\Exreg(X)=\mathop{\sup}\{m+\mathop{\sup}\{\ell_m^{i,j}\mid i,j\}\mid \gExt^m_A(S,X)\neq 0, m\in \mathbb{Z}\}.$$
    
    By Lemma \ref{inj-decom}, $\Gamma_A(I^m)$ is the injective hull of $\soc I^m$, and by Proposition \ref{inj-hull-S-i},
    $$\Gamma_A(I^m)\cong \mathop{\bigoplus}\limits_i \big(\mathop{\bigoplus}\limits_j \D(e_iA)(-\ell_m^{i,j})\big).$$ 
    Since $R^m\Gamma_A(X)$ is a subquotient of $\Gamma_A(I^m)$, it follows that
    \begin{align*}
    \sdeg(R^m\Gamma_A(X))\leqslant & \sdeg(\Gamma_A(I^m))\\
    = &\mathop{\sup}\{\ell_m^{i,j}\mid i,j\}\\
    =& \sdeg(\gExt^m_A(S,X)).
    \end{align*}
    Hence 
    \begin{align*}
       \CMreg(X)=&\sup \{m+\sdeg(R^m\Gamma_A(X)) \mid R^m\Gamma_A(X)\neq 0, m\in \mathbb{Z}\}\\
       \leqslant& \sup \{m+\sdeg(R^m\Gamma_A(X)) \mid \gExt^m_A(S,X)\neq 0, m\in \mathbb{Z}\} \\
       \leqslant& \sup \{m+\sdeg(\gExt^m_A(S,X)) \mid \gExt^m_A(S,X)\neq 0, m\in \mathbb{Z}\}\\
       =&\Exreg(X).
    \end{align*}
\end{proof}
It may happen that $\Exreg(X) > \CMreg(X)$. For example, if $A_0$ is not semisimple, then $\Exreg(S) > \CMreg(S)=0$, as showed in the following. 
   Since $S=A/J=A_0/J(A_0)$ is finite-dimensional, $S$ is a torsion module. Thus $R\Gamma_A(S)\cong S$ and $\CMreg(S)=0$. 
   On the other hand,
   let $P^{\bullet}$ be a minimal graded projective resolution of ${}_AS$. If $A_0$ is not semisimple, then the minimal homogeneous generating set of $P^{-1}$
   contains an element $x$ of degree zero. Hence $x \in \gHom_A(P^{-1},S) \cong \gExt_A^1(S,S)$. Therefore $$\Exreg(S)\geqslant 1+\sdeg(\gExt_A^1(S,S))\geqslant 1+0=1>0.$$ 

The following lemma concerns the relation between the numerical regularities of modules in a short exact sequence, which is proved by using the corresponding induced long exact sequences (see \cite[Lemma 3.7]{KWZ2}), where $$\reg \in \{\CMreg, \cmreg, \Torreg, \torreg, \Extreg, \extreg, \exreg, \Exreg\}.$$

\begin{lemma}\label{exact sequence}
Let $A$ be an $\mathbb{N}$-graded algebra. For any distinguished triangle $X \to Y \to Z \to X[1]$ in $\mathbf{D}^b(\Gr A)$,
\begin{itemize}
    \item [(1)] $\reg(Y) \leqslant \max\{\reg(X),\reg(Z)\}$.
    \item [(2)] $\reg(X)\leqslant \max\{\reg(Y),\reg(Z)+1\}.$
    \item [(3)] $\reg(Z)\leqslant \max\{\reg(X)-1,\reg(Y)\}.$
\end{itemize}
\end{lemma}

\subsection{More relationship between the regularities} \label{relationship between the regularities}
As showed in Lemma \ref{lem-CM-ex}, 
$\CMreg(X) \leqslant \Exreg(X)$ for any $0\neq X\in \mathbf{D}^+(\Gr A)$.
If $A_0$ is semisimple and $A$ has a balanced dualizing complex, then we show in Theorem \ref{ex} that
$\CMreg(X) = \Exreg(X)$ for any  $0\neq X \in \mathbf{D}^+_{fg}(\Gr A)$, which is a generalized form of \cite[Theorem 3.9]{Ngu} in the noncommutative graded case.
To prove Theorem \ref{ex}, we first generalize \cite[Lemma 5.2]{DW} in  Proposition \ref{CM and Ext}, which says that if $A_0$ is semisimple, then for any $X\in \mathbf{D}^-(\gr A)$, 
$$\extreg(X)= -\ideg(X).$$
Proposition \ref{CM and Ext} is one of the key facts we developed to compare the regularities.
\begin{lemma}\label{CM and Ext inequality}
 Let $A$ be a noetherian $\mathbb{N}$-graded algebra. For any $0\neq X\in \mathbf{D}^-(\gr A)$, $$\extreg(X)\geqslant -\ideg(X).$$
\end{lemma}
\begin{proof}
     Without loss of generality, we may assume that $X^n=0$ for all $n \geqslant 1$ and $\sup(X)=0$. 
 It follows from Lemma \ref{notations of X} that 
 $$\extreg(X)=\sdeg(R\gHom_{A}(X,S))=\mathop{\sup} \{m-\lxm{m}{X} \mid P^{-m} \neq 0\}.$$
On the other hand,
 \begin{align*}
    -\ideg(X)&=-\inf \{-i+ \ideg (H^{-i}(X)) \mid H^{-i}(X)\neq 0\}\\
    &\leqslant -\inf\{-i+\ideg(P^{-i})\mid  H^{-i}(X)\neq 0\}\\
    &\leqslant -\inf\{-i+\ideg(P^{-i})\mid  P^{-i}\neq 0\}\\
    &=-\inf\{-i+\lxm{i}{X}\mid P^{-i}\neq 0\}\\
    &= \mathop{\sup} \{i-\lxm{i}{X} \mid P^{-i} \neq 0\}.
    \end{align*}
Therefore,  $-\ideg(X) \leqslant \extreg(X)$.
\end{proof}

\begin{proposition}\label{CM and Ext}
Let $A$ be a noetherian $\mathbb{N}$-graded algebra with $A_0$ semisimple. Then, for any $0\neq X \in \mathbf{D}^-(\gr A),$
$$\extreg(X)=-\ideg(X).$$
\end{proposition}
\begin{proof}
 Without loss of generality, we may assume that $X^n=0$ for all $n \geqslant 1$ and $\sup(X)=0$. 
 It follows from Lemma \ref{notations of X} that 
 $$\extreg(X)=\sdeg(R\gHom_{A}(X,A_0))=-\mathop{\inf} \{\lxm{m}{X}-m \mid P^{-m} \neq 0\}.$$
 
 Let $C=\{ -i \mid \ker d^{-i} \cap P^{-i}_{\lxm{i}{X}} \neq 0,i \in \mathbb{N}\}$, where $P^{-i}_{\lxm{i}{X}}$ is the degree $\lxm{i}{X}$ part of the graded module $P^{-i}$. 
In fact, $-i \in C$ if and only if  $\ideg(H^{-i}(X))=\lxm{i}{X}$.
 
Since $P^{\bullet}$ is a minimal graded projective resolution of $X$ and $A_0$ is semisimple, $\im d^{-i-1} \subseteq A_{\geqslant 1}P^{-i}$. Hence $\ideg(\im d^{-i-1}) > \lxm{i}{X}$. If $-i \in C$, then $\ideg(\im d^{-i-1}) > \lxm{i}{X}=\ideg(\ker d^{-i})$, and
$\ideg(H^{-i}(X))=\lxm{i}{X}$. Conversely, if $\ideg(H^{-i}(X))=\lxm{i}{X}$, then $\ideg(\ker d^{-i})=\lxm{i}{X}$, and $-i \in C$.
Hence, $C=\{-i  \mid \ideg (H^{-i}(X)) = \lxm{i}{X}, i \in \mathbb{N}\}$, and
$$0 \in C \subseteq \{-i  \mid H^{-i}(X) \neq 0,i \in \mathbb{N}\} \subseteq \{-i  \mid P^{-i} \neq 0,i \in \mathbb{N}\} .$$
  
Suppose $C \subsetneqq \{-i  \mid P^{-i} \neq 0,i \in \mathbb{N}\} $. If $P^{-m} \neq 0$ and $-m \notin C$\, ($m > 0$), then $\ker d^{-m} \bigcap P^{-m}_{\lxm{m}{X}} = 0$. Let $x \in P^{-m} $ with $\deg(x)=\lxm{m}{X}$ such that $d^{-m}(x) \neq 0$. It follows from $d^{-m}(x) \in A_{\geqslant 1}P^{-m+1}$ that $P^{-m+1} \neq 0$ and 
$$\lxm{m}{X}=\deg(d^{-m}(x))\geqslant \lxm{m-1}{X}+1.$$

Since $0 \in C$, 
there is an element $-r \in C$ such that $-k \notin C$ for all $k$ satisfying that $-m \leqslant -k < -r$. Then 
$$ \lxm{m}{X} \geqslant \lxm{m-1}{X} + 1 \geqslant \lxm{m-2}{X} +2 \geqslant \cdots  \geqslant \lxm{r}{X}+ m-r,$$
and so, $\lxm{m}{X} -m \geqslant \lxm{r}{X} - r.$ 
It follows that
\begin{align} \label{inf-inequality}
   \inf \{\lxm{m}{X}-m\mid P^{-m}\neq 0, -m\notin C\} \geqslant \inf \{\lxm{r}{X}-r\mid -r \in C\}.
\end{align}
Therefore
\begin{align} \label{exreg-id}
\mathop{\inf} \{\lxm{m}{X}-m \mid P^{-m}\neq 0, m\in \mathbb{N} \}
    =\mathop{\inf} \{\lxm{r}{X}-r \mid -r \in C\}.
\end{align}

 Suppose $C = \{-i  \mid P^{-i} \neq 0,i \in \mathbb{N}\}$. Then \eqref{exreg-id} holds trivially.

Therefore, \eqref{exreg-id} holds always.

It follows from $C\subseteq\{-i \mid H^{-i}(X) \neq 0,i\in \mathbb{N}\}$ that
\begin{align*}
    -\ideg(X)&=-\inf \{-i+ \ideg (H^{-i}(X)) \mid H^{-i}(X)\neq 0\}\\
    &\geqslant-\inf \{-i+ \ideg (H^{-i}(X)) \mid  -i\in C\}\\
&=-\inf \{-i+\lxm{i}{X} \mid -i \in C \}  \\
    &=\extreg(X) \quad (\textrm{by } \eqref{exreg-id}).
\end{align*}
It follows from Lemma \ref{CM and Ext inequality} that $\extreg(X)=-\ideg(X)$.
\end{proof}

\begin{remark}\label{remark of CM and Ext} 
   It follows from the proof of Proposition \ref{CM and Ext} that for any fixed integer $c\geqslant 0$,  
   \begin{align*}
   &-\inf \{ -i+ \ideg(H^{-i}(X)) \mid 0 \leqslant i\leqslant c\}\\
   =&\sup\{i + \sdeg(\gExt_A^i(X,A_0)) \mid 0 \leqslant i\leqslant c\}.
   \end{align*}
This generalizes \cite[Corollary 5.2]{Jo3} and \cite[Theorem 5.2]{Sch}.
\end{remark}

Suppose $A_0$ is semisimple and $A$ has a balanced dualizing complex. The following theorem says that the CM-regularity
is the same as the ex-regularity for any $0\neq X \in \mathbf{D}^+_{fg}(\Gr A)$, generalizing  \cite[Theorem 3.9]{Ngu}.  Proposition \ref{CM and Ext} plays a key role in the proof. 

\begin{theorem}\label{ex}
Suppose $A$ is a noetherian $\mathbb{N}$-graded algebra with $A_0$ semisimple and $A$ has a balanced dualizing complex $R$. 
Then, for any $0\neq X \in \mathbf{D}^+_{fg}(\Gr A)$, 
$$\CMreg(X)=\Exreg(X).$$
\end{theorem}
\begin{proof}
Let $Z=\D(R\Gamma_A(X))\in \mathbf{D}^-(\gr A^o)$.
Then,  by Definition \ref{CM-reg-def}
$$\CMreg(X)=\sdeg (R\Gamma_A(X)) = -\ideg(Z).$$
 
 Since $R$ is a balanced dualizing complex of $A$, $R\gHom_A(A_0,X)\cong R\gHom_{A^o}(Z,A_0)$.
 By Definition \ref{exreg-def},
$$\Exreg(X)=\sdeg (R\gHom_A(A_0, X))=\sdeg (R\gHom_{A^o}(Z, A_0))=\extreg(Z).$$
It follows from  the right version of Proposition \ref{CM and Ext} that
$$\CMreg(X)=-\ideg(Z)=\extreg(Z)=\Exreg(X).$$
\end{proof}

It follows from Lemma \ref{dual-ext-ex} (resp. \eqref{Ex-Ext and ex-ext-3}) that $\CMreg(X)=\extreg(\D(R\Gamma_A(X)))$ (resp. $\CMreg(X)=\extreg(\D(X)$).
If $X\in \mathbf{D}^b(\gr A)$, then it follows from Theorem \ref{ex} and Lemma \ref{CM-reg-finite} that $\Exreg(X)=\CMreg(X)<\infty$.

If $A_0$ is not semisimple, Theorem \ref{ex} may not be true. For example, if $A=A_0$ and $\pdim(S)=+\infty$, then $\Exreg(S)=+\infty$ and $\CMreg(S)=0$. This also shows that $\Exreg(X)$ may be infinite.

It follows from Lemma \ref{CM-reg-finite} and Theorem \ref{ex} that if $A$ has a balanced dualizing complex and $A_0$ is semisimple then $\Exreg({}_{A}A)=\Exreg(A_A)$, which is sometimes denoted by $\Exreg(A)$. If $A$ has a balanced dualizing complex $R$, then it follows from Lemma \ref{dual-ext-ex} (1) that $\exreg({}_{A}A)=\Torreg(R_A)$. Therefore, $\exreg({}_{A}A)=\exreg(A_A)$ if and only if $\Torreg({}_{A}R)=\Torreg(R_A)$. 

\begin{lemma} \label{Exreg-is-greater}
  Let $A$ be a noetherian $\mathbb{N}$-graded algebra. For any $0\neq X\in \mathbf{D}^-(\gr A)$, $$\Extreg(X) \geqslant -\ideg(R\gHom_A(X,A)).$$
\end{lemma}
\begin{proof}
It follows Lemma \ref{notations of X} that $\Extreg(X)=-\inf\{m-\uxm{m}{X} \mid P^{-m}\neq 0\}$.

On the other hand, 
\begin{align*}
   -\ideg(R\gHom_A(X,A))=&-\inf\{i+\ideg(\gExt_A^i(X,A))\mid \gExt_A^i(X,A)\neq 0\}\\
   &\leqslant -\inf\{i+\ideg(\gExt_A^i(X,A))\mid P^{-i}\neq 0\}\\
   &\leqslant -\inf\{i+\ideg(\gHom_A(P^{-i},A))\mid P^{-i}\neq 0\}\\
   &= -\inf\{i-\uxm{i}{X} \mid P^{-i}\neq 0\}.
 \end{align*} 
 Hence $-\ideg(R\gHom_A(X,A))\leqslant \Extreg(X)$.
\end{proof}

If $A_0$ is semisimple, then the inequality in Lemma \ref{Exreg-is-greater} is in fact an equality as proved in the following theorem.
Theorem \ref{Extreg} simplifies to \cite[Proposition 1.1]{Tr} within the context of commutative polynomial algebra, using the terminology therein. The result also generalizes \cite[Theorem 5.1]{Jo3}.

\begin{theorem}\label{Extreg}
Suppose $A$ is a noetherian $\mathbb{N}$-graded algebra with $A_0$ semisimple. If $0\neq X\in \mathbf{D}^b(\gr A)$ with $\pdim X < \infty$, then $$\Extreg(X)=-\ideg(R\gHom_A(X,A)).$$
\end{theorem}
\begin{proof}
     Let  $Z=R\gHom_A(X,A)$.  
    Since $\pdim(X) < \infty$, 
   $$R\gHom_A(X,A_0)\cong R\gHom_A(X,A)\, {}^L\!\otimes_A A_0\cong Z \, {}^L\!\otimes_A A_0.$$   
Then,  by the definitions
\begin{align*}
    \Extreg(X)&=-\ideg(R\gHom_A(X,A_0))\\
    &=-\ideg(Z \,\, {}^L\!\otimes_A A_0)\\
    &=\torreg(Z) \\
    &=\extreg(Z) \quad (\textrm{by Lemma } \ref{Extreg-Torreg}).
\end{align*}
It follows from the right version of Proposition \ref{CM and Ext} that
$$\extreg(Z)=-\ideg(Z)=-\ideg(R\gHom_A(X,A)).$$
Hence $\Extreg(X)=-\ideg(R\gHom_A(X,A)).$
\end{proof}

Note that the condition $\pdim(X) < \infty$ in Theorem \ref{Extreg} is necessary. For example, if $A=k[x]/(x^2)$, then $\pdim(k)=+\infty$ and $\Extreg(k)=0$. Since $\gExt_A^i(k, A)=0$, $i\neq 0$ and $\gExt_A^0(k, A)=k(-1)$, 
$$-\mathop{\inf} \{i+j \mid \gExt_A^i(k, A)_j \neq 0, i,j\in \mathbb{Z}\}=-1.$$ 

Suppose $A$ is a noetherian $\mathbb{N}$-graded algebra with $A_0$ semisimple. Let  $0\neq X \in \mathbf{D}^b(\gr A)$.  Table $1$ in the following presents the definitions and the relations between the various regularities, that is, the supremum degrees and the infimum degrees of the complexes of interest.
\renewcommand\arraystretch{2}
\begin{table}[htbp]
\begin{center}
	\begin{tabular}{ | m{2.6cm}<{}|m{4.2cm}<{\centering} | m{4.4cm}<{\centering} |}
		\hline
			\diagbox{$\mathbf{Ri}$\qquad}{$\mathbf{Cj}$\qquad} &   $\sdeg$   &   $-\ideg$   \\ \hline
		 $A_0 \, {}^L\!\otimes_A X$ & \makecell{$\Torreg(X)$ \\ $ = \Extreg(X)$ \\ $= \exreg(\D(X))$ \\ $\stackrel{\exists \text{ bdc}}{=} \exreg(\D(R\Gamma_A(X)))$\\
   $\stackrel{\exists \text{ bdc}}{=} \Extreg(R\Gamma_A(X))$}   & \makecell{$\torreg(X)$ \\ $=\extreg(X)$ \\ $=\Exreg(\D(X))$ \\ $=-\ideg(X)$ \\$\stackrel{\exists \text{ bdc}}{=} \Exreg(\D(R\Gamma_A(X)))$\\ $\stackrel{\exists \text{ bdc}}{=} \extreg(R\Gamma_A(X))$}  \\ \hline
		 $R\gHom_A(X,A_0)$    & \makecell{$\extreg(X)$ \\ $=\torreg(X)$ \\ $=\Exreg(\D(X))$ \\ $=-\ideg(X)$ \\$\stackrel{\exists \text{ bdc}}{=} \Exreg(\D(R\Gamma_A(X)))$\\ $\stackrel{\exists \text{ bdc}}{=} \extreg(R\Gamma_A(X))$}  & \makecell{$\Extreg(X)$ \\ $ = \Torreg(X)$ \\ $= \exreg(\D(X))$ \\ $\stackrel{\exists \text{ bdc}}{=} \exreg(\D(R\Gamma_A(X)))$\\
   $\stackrel{\exists \text{ bdc}}{=} \Extreg(R\Gamma_A(X))$}\\ \hline
		 $R\gHom_A(A_0,X)$	 &\makecell{$\Exreg(X)$ \\ $=\extreg(\D(X))$ \\ $\stackrel{\exists \text{ bdc}}{=}\CMreg(X)$ \\ $\stackrel{\exists \text{ bdc}}{=}\extreg(\D(R\Gamma_A(X)))$}    &   \makecell{$\exreg(X)$ \\ $=\Extreg(\D(X))$ \\ $\stackrel{\exists \text{ bdc}}{=}\Extreg(\D(R\Gamma_A(X)))$}   \\ \hline
 $R\Gamma_A(X)$	& \makecell{$\CMreg(X)$ \\  $\stackrel{\exists \text{ bdc}}{=}\Exreg(X)$\\ $\stackrel{\exists \text{ bdc}}{=}\extreg(\D(X))$ \\ $\stackrel{\exists \text{ bdc}}{=}\extreg(\D(R\Gamma_A(X)))$}  & \makecell{$\cmreg(X)$ \\ $\stackrel{X\in \mathbf{D}^b(\fd A)} {=\!=} \extreg(X)$} \\ \hline
 $R\gHom_A(X,A)$	 &  & $\stackrel{\text{pdim X $< \infty$}}{=\!=}\Extreg(X)$   \\ 
		\hline
	\end{tabular}
\end{center}
\caption{}
\end{table}

Let us give a short explanation. In the table, \text{``bdc"} means that under the condition that $A$ has a balanced dualizing complex. We use RiCj to represent the position at the intersection of the $i$th row  and $j$th column. 

R1C1 and R2C2 hold by definitions, Lemma \ref{Extreg-Torreg}, \eqref{Ex-Ext and ex-ext-4}, Lemma \ref{dual-ext-ex} (2) and Lemma \ref{iso about local cohomological} (by taking $Y=S$) if $A$ has a balanced dualizing complex.

R2C1 and R1C2 hold by definitions, Lemma \ref{Extreg-Torreg}, \eqref{Ex-Ext and ex-ext-4}, Proposition \ref{CM and Ext}, Lemma \ref{dual-ext-ex} (2) and Lemma \ref{iso about local cohomological} (by taking $Y=S$) when $A$ has a balanced dualizing complex.

R3C1 holds by the definition, \eqref{Ex-Ext and ex-ext-3}, Theorem \ref{ex} and Lemma  \ref{dual-ext-ex} (1) when $A$ has a balanced dualizing complex.

R4C1 holds by the definition, Theorem \ref{ex}, \eqref{Ex-Ext and ex-ext-3} and Lemma  \ref{dual-ext-ex} (1) when $A$ has a balanced dualizing complex.

R3C2 holds by Definition \ref{exreg-def}, \eqref{Ex-Ext and ex-ext-3} and Lemma  \ref{dual-ext-ex} (1) when $A$ has a balanced dualizing complex.

R4C2 is Definition \ref{CM-reg-def}, and if $X\in \mathbf{D}^b(\fd A)$, then $R\Gamma_A(X)\cong X$ and $\cmreg(X)=-\ideg(X)=\extreg(X)$.

R5C2 holds by Theorem \ref{Extreg}.

R5C1 means $\sdeg(R\gHom_A(X,A))$.

\section{Some interconnections among the regularities} 
\subsection{Supremum and infimum degrees of Hom and Tensor complexes} \label{degree-of-Hom-Tensor}
In this subsection, we establish several inequalities about the supremum or infimum degrees of $R\gHom_A(X,Y)$ and $Y\, {}^L\!\otimes_A X$ for some $X, Y \in \mathbf{D}(\Gr A)$ in the following propositions. Note by our abused terminology in this paper that $A$ is a noetherian $\mathbb{N}$-graded algebra means that $A$ is a noetherian locally finite $\mathbb{N}$-graded algebra, and so $\mathbf{D}(\gr A)$ is a full subcategory of $\mathbf{D}_{lf}(\Gr A)$.
Cohomological spectral sequences are used in the proofs of the following propositions.
We refer to \cite[Chapter 5]{Wei} for the notation of cohomological spectral sequences.

\begin{proposition} \label{sdeg-ideg inequality 1}
    Let $A$ be a noetherian $\mathbb{N}$-graded algebra. Then for any $0\neq X\in \mathbf{D}^-(\gr A)$, $0\neq Y\in \mathbf{D}^+(\Gr A)$, 
   \begin{align}  \label{inequality-1}
   -\ideg(R\gHom_A(X,Y))\leqslant \Extreg(X)-\ideg(Y).
   \end{align}
    \end{proposition}
    
\begin{proof}
If $\Extreg(X) = +\infty$ or $-\ideg(Y)=+\infty$, then the inequality is trivially true. 

  Suppose $\Extreg(X) = r < \infty$ and $-\ideg(Y)=p < \infty$.

By Lemma \ref{notations of X}, $r=\sup\{\uxm{m}{X}-m\mid P^{-m}\neq 0\}$.
It follows that $\uxm{m}{X}\leqslant r+m$ for all $m \in \mathbb{Z}$ such that $P^{-m} \neq 0$.

Since $p=-\ideg(Y)=-\inf\{n+\ideg(H^n(Y))\mid H^n(Y)\neq 0, n\in \mathbb{Z}\}$, 
$\ideg(H^n(Y))\geqslant -p-n$ for any $n\in \mathbb{Z}$ such that $H^n(Y)\neq 0$.

Now consider the double complex $C^{\bullet \bullet}$ in the first quadrant given by 
$$C^{m,n}=\gHom_A(P^{-m},Y^n).$$
Then there is a convergent spectral sequence 
$${}^I\!E_2^{m,n}=\gExt_A^m(X,H^n(Y))\Rightarrow \gExt_A^{m+n}(X,Y).$$

For any $P^{-m}=\mathop{\bigoplus}\limits_{i} \big(\mathop{\bigoplus}\limits_{j} Ae_i(-s_m^{i,j})\big) \neq 0$ and $H^n(Y)\neq 0$, 
\begin{align*}
    &\ideg(\gHom_A(P^{-m},H^n(Y))) \\
    =& \ideg\big(\mathop{\bigoplus}\limits_{i} (\mathop{\bigoplus}\limits_{j} e_iH^n(Y)(s_m^{i,j}))\big) = \inf \{\ideg(e_iH^n(Y)) - s_m^{i,j}\mid i,j\}\\
    \geqslant & \inf \{\ideg(H^n(Y)) - s_m^{i,j}\mid i,j\} =\ideg(H^n(Y)) - \sup\{s_m^{i,j}\mid i,j\}  \\
    =& -p-n-\uxm{m}{X} \geqslant -p-n-m-r.
\end{align*} 
As $\gExt_A^m(X,H^n(Y))$ is a subquotient of $\gHom_A(P^{-m},H^n(Y))$, so 
$$\ideg(\gExt_A^m(X,H^n(Y)))\geqslant -p-n-m-r.$$
It follows from the convergence of the spectral sequence that 
$$\ideg(\gExt_A^{m+n}(X,Y))\geqslant -p-n-m-r.$$ 
Hence
\begin{align*}  
&-\ideg(R\gHom_A(X,Y))\\
=&-\inf\{i+\ideg(\gExt_A^{i}(X,Y))\mid \gExt_A^{i}(X,Y)\neq 0\}\\
\leqslant & p+r
\end{align*} 
and the conclusion holds.
\end{proof}

Similarly, the following two propositions are proved by taking the minimal graded injective resolution in the second variable
 and considering the cohomological spectral sequences of the double complex induced by Hom.

\begin{proposition} \label{sdeg-ideg inequality 1-1}
    Let $A$ be a noetherian $\mathbb{N}$-graded algebra. Then the following hold. 
 \begin{itemize}
        \item [(1)] For any $0\neq X\in \mathbf{D}^-(\Gr A)$ with 
         torsion cohomologies,
        $0\neq Y\in \mathbf{D}^+(\Gr A)$,  
\begin{align}  \label{inequality-3-1}
   \sdeg(R\gHom_A(X,Y))\leqslant -\ideg(X)+\Exreg(Y).
   \end{align}
   \item [(2)] If $A$ has a balanced dualizing complex, then for any $0\neq X \in \mathbf{D}^-_{fg}(\Gr A)$
   and $0\neq Y \in \mathbf{D}^+_{fg}(\Gr A)$,  
\begin{align} \label{inequality-3-2}
   \sdeg(R\gHom_A(X,Y))\leqslant \cmreg(X)+\Exreg(Y).
   \end{align}
\end{itemize}     
\end{proposition}

\begin{proof} (1)
  We may assume that $-\ideg(X)=p <\infty$ and $\Exreg(Y)=q < \infty$. Then, for any $-n \in \mathbb{Z}$ with $H^{-n}(X)\neq 0$, $\ideg(H^{-n}(X))\geqslant -p+n$. 
  
 Let $I^\bullet$ be the minimal graded injective resolution of $Y$. Then $$\gExt_A^m(S,Y)\cong \gHom_A(S,I^m)\cong \soc I^m.$$ 
 Suppose that $\soc I^m\cong \mathop{\bigoplus}\limits_i \big(\mathop{\bigoplus}\limits_j Se_i(-s_m^{i,j})\big)$. 
  It follows from Proposition \ref{inj-hull-S-i} that 
   $$E(\soc I^m)\cong \mathop{\bigoplus}\limits_i \big(\mathop{\bigoplus}\limits_j \D(e_iA)(-s_m^{i,j})\big).$$
   Then,  
    $\sdeg(E(\soc I^m)) =\max\{s_m^{i,j} \} = \sdeg(\gExt_A^m(S,Y)) \leqslant q-m. $
    
The double complex $C^{\bullet \bullet}$ in the first quadrant given by 
$$C^{m,n}=\gHom_A(X^{-m},I^n),$$
induces a convergent spectral sequence (arising from the second filtration)
$${}^{I\!I}\!E_2^{m,n}=\gExt_A^m(H^{-n}(X),Y)\Rightarrow \gExt_A^{m+n}(X,Y).$$
Since, by assumption,  
$H^{-n}(X)$ is torsion for all $n$,
$$\gHom_A(H^{-n}(X), I^m) = \gHom_A(H^{-n}(X), E(\soc I^m)).$$
Then, for any $m,n \in \mathbb{Z}$ with $I^m\neq 0$ and $H^{-n}(X)\neq 0$, 
\begin{align*}
    \sdeg(\gHom_A(H^{-n}(X), I^m))=& \sdeg(\gHom_A(H^{-n}(X), E(\soc I^m))) \\
    \leqslant& q-m+p-n.
\end{align*} 
Hence, $\sdeg(\gExt_A^m(H^{-n}(X),Y))\leqslant q-m+p-n$.

It follows from the convergence of the spectral sequence that
$$\sdeg(\gExt_A^{m+n}(X,Y))\leqslant q-m+p-n.$$ 
Therefore
\begin{align*}  
&\sdeg(R\gHom_A(X,Y))\\
=&\sup\{i+\sdeg(\gExt_A^{i}(X,Y))\mid \gExt_A^{i}(X,Y)\neq 0\}\\
\leqslant & p+q
\end{align*} 
and (1) is proved.

(2) Now, suppose $A$ has a balanced dualizing complex and $0\neq X \in \mathbf{D}^-_{fg}(\Gr A)$,  $0\neq Y \in \mathbf{D}^+_{fg}(\Gr A)$. Then 
by Lemma \ref{iso about local cohomological} or \cite[Proposition 1.1]{Jo4},
$R\gHom_A(X,Y)\cong R\gHom_A(R\Gamma_A(X),Y).$ 

Since $A$ has a balanced dualizing complex, $R\Gamma_A(X) \in \mathbf{D}^-(\Gr A)$ and $R^i\Gamma_A(X)$ is torsion for all $i$,
it follows from (1) that 
\begin{align*}
 &\sdeg(R\gHom_A(X,Y))\\
 =&\sdeg(R\gHom_A(R\Gamma_A(X),Y)) \leqslant -\ideg(R\Gamma_A(X))+\Exreg(Y)\\
 =&\cmreg(X)+ \Exreg(Y).
 \end{align*} 
\end{proof}

\begin{lemma}\label{iso about local cohomological} Suppose $A$ has a balanced dualizing complex.
If $X  \in \mathbf{D}^-(\Gr A), Y \in \mathbf{D}^+_{fg}(\Gr A)$, or $Y = \D(Y')$ for some $Y' \in \mathbf{D}^-(\Gr A^o)$  with torsion cohomologies, then 
$$R\gHom_A(X,Y)\cong R\gHom_A(R\Gamma_A(X),Y).$$
\end{lemma}
\begin{proof} The proof is almost the same as \cite[Proposition 1.1]{Jo4}. It follows from $A$ has finite cohomological dimension and $\D(Y)$ has torsion cohomologies that
\begin{align*}
R\gHom_A(R\Gamma_{A^o}(A), Y)\cong & R\gHom_A(R\Gamma_{A^o}(A), \D(\D(Y)))\\
\cong& \D\big(\D(Y)\, {}^L\!\otimes_A R\Gamma_{A^o}(A) \big)  \quad (\D(Y) \in \mathbf{D}^-(\Gr A) )\\
\cong &\D\big(R\Gamma_{A^o}(\D(Y))\big) \quad (\textrm{by \cite[Proposition 2.1]{Jo1}} \\
\cong & \D(\D(Y)) \cong Y \quad (\textrm {by \cite[Lemma 4.4]{VdB}}.
\end{align*} 
Hence
\begin{align*}
R\gHom_A(R\Gamma_A(X), Y)\cong & R\gHom_A(R\Gamma_{A}(A) {}^L\!\otimes_A  X, Y)\big) \\
\cong & R\gHom_{A}(X, R\gHom_A(R\Gamma_A(A), Y)) \\
\cong & R\gHom_{A}(X, R\gHom_A(R\Gamma_{A^o}(A), Y)) \\
\cong & R\gHom_A(X,Y).
\end{align*} 
\end{proof}

\begin{proposition} \label{sdeg-ideg inequality 1-4}
Let $A$ be a noetherian $\mathbb{N}$-graded algebra. Then the following hold.
\begin{itemize}
        \item [(1)] For any $0\neq X\in \mathbf{D}^-(\Gr A)$
        and $0\neq Y\in \mathbf{D}^-(\gr A^o)$,
\begin{align}\label{inequality-2-1}
-\ideg(R\gHom_A(X,\D(Y)))\leqslant \sdeg(X)+\Extreg(Y).
\end{align}
\item [(2)] If $A$ has a balanced dualizing complex, then for any $0\neq X \in \mathbf{D}^-_{fg}(\Gr A)$ and $0\neq Y \in \mathbf{D}^-(\gr A^o)$ with torsion cohomologies,
\begin{align}  \label{inequality-2-2}
   -\ideg(R\gHom_A(X,\D(Y)))\leqslant \CMreg(X)+\Extreg(Y).
   \end{align}
\end{itemize}   
\end{proposition}
\begin{proof} (1)
  We may assume that $\sdeg(X)=p <\infty$ and $\Extreg(Y)=q < \infty$. Then, for any $-n \in \mathbb{Z}$ with $H^{-n}(X)\neq 0$, $\sdeg(H^{-n}(X))\leqslant p+n$. 

 Let $Q^\bullet \to Y$ be a minimal graded projective resolution of $Y \in \mathbf{D}^-(\gr A^o)$. 
 If  $0\neq Q^{-m}=\mathop{\bigoplus}\limits_i \big(\mathop{\bigoplus}\limits_j e_iA(-s_m^{i,j})\big)$, then $\ideg \gExt_{A^o}^m(Y,S)=\inf\{-s_m^{i,j}\}=-\uxm{m}{Y}$.  By Lemma \ref{notations of X},
 $\Extreg(Y)= \sup\{\uxm{m}{Y}-m\mid Q^{-m}\neq 0\}.$
 
 Note $\D(Y) \to \D(Q^\bullet)$ is a minimal graded injective resolution of $\D(Y)$. 
Therefore, 
  \begin{align*}
   \gExt_{A^o}^m(Y,S) \cong \gExt_{A}^m(S,\D(Y)) = \gHom_{A}(S, \D(Q^{-m})).
\end{align*}
For any $m \in \mathbb{Z}$ with $\gExt_{A^o}^m(Y, S)\neq 0$,  $\ideg (\gExt_{A^o}^m(Y,S)) =- \uxm{m}{Y} \geqslant -q - m$.

The double complex $C^{\bullet \bullet}$ in the first quadrant given by 
$$C^{m,n}=\gHom_A(X^{-m},\D(Q^{-n})),$$
induces a convergent spectral sequence (arising from the second filtration)
$${}^{I\!I}\!E_2^{m,n}=\gExt_A^m(H^{-n}(X),\D(Y))\Rightarrow \gExt_A^{m+n}(X,\D(Y)).$$
For any $m,n \in \mathbb{Z}$ with $Q^{-m}\neq 0$ and $H^{-n}(X)\neq 0$, 
\begin{align*}
&\gHom_A(H^{-n}(X), \D(Q^{-m})) \\
=& \gHom_A\big(H^{-n}(X),  \mathop{\bigoplus}\limits_i \big(\mathop{\bigoplus}\limits_j \D(e_iA)(s_m^{i,j}))\big)\\
 =& \mathop{\bigoplus}\limits_i \big(\mathop{\bigoplus} \limits_j \big(\gHom_A(H^{-n}(X),  \D(e_iA))(s_m^{i,j})\big)\big)\\
 =& \mathop{\bigoplus}\limits_i \big(\mathop{\bigoplus} \limits_j \D \big(e_iA \otimes_A H^{-n}(X)\big)(s_m^{i,j})\big).\\
\end{align*}
Hence
\begin{align*}
    &\ideg(\gHom_A(H^{-n}(X), \D(Q^{-m})))\\
    =&\inf \{ -\sdeg(e_iH^{-n}(X)) - s_m^{i,j} \mid i, j\} \\
    \geqslant &\inf \{- \sdeg(H^{-n}(X)) - s_m^{i,j} \mid i, j\} \\
     = &-\sdeg(\gHom_A(H^{-n}(X)) - \sup\{s_m^{i,j} \mid i, j\} \\
      = &-\sdeg(\gHom_A(H^{-n}(X)) - \uxm{m}{Y} \\
    \geqslant& -p-n-q-m,
\end{align*} 
and $\ideg(\gExt_A^m(H^{-n}(X), \D(Y)))\geqslant -p-n-q-m$.

It follows from the convergence of the spectral sequence that
$$\ideg(\gExt_A^{m+n}(X,\D(Y)))\geqslant -p-n-q-m.$$ 
Therefore
\begin{align*}  
&-\ideg(R\gHom_A(X,\D(Y)))\\
=&-\inf\{i+\ideg(\gExt_A^{i}(X,\D(Y)))\mid \gExt_A^{i}(X,\D(Y))\neq 0\}\\
\leqslant & p+q
\end{align*} 
and (1) holds.

 (2) Now, suppose $A$ has a balanced dualizing complex, $0\neq X \in \mathbf{D}^-_{fg}(\Gr A)$ and $0\neq Y \in \mathbf{D}^-(\gr A^o)$ with torsion cohomologies.
By Lemma \ref{iso about local cohomological},  
$$R\gHom_A(X,\D(Y))\cong R\gHom_A(R\Gamma_A(X),\D(Y)).$$ 
Since $R\Gamma_A(X) \in \mathbf{D}^-(\Gr A)$, it follows from 
(1) that 
\begin{align*}
 &-\ideg(R\gHom_A(X,\D(Y)))\\
 =&-\ideg(R\gHom_A(R\Gamma_A(X),\D(Y))) \leqslant \sdeg(R\Gamma_A(X))+\Extreg(Y)\\
 =&\CMreg(X)+ \Extreg(Y).
 \end{align*} 
\end{proof}

By taking $Y=S$ in Proposition \ref{sdeg-ideg inequality 1-4} (2), it follows  that $\Extreg(X)\leqslant \CMreg(X)+ \Extreg(S)$, which generalizes \cite[Theorem 2.5]{Jo4}.

\begin{proposition}\label{sdeg-ideg inequality 2}
    Let $A$ be a noetherian $\mathbb{N}$-graded algebra. Then the following statements hold.
    \begin{itemize}
        \item [(1)] For any $0\neq X\in \mathbf{D}^-(\gr A)$, $0\neq Y\in \mathbf{D}^-(\Gr A^{o})$,
    \begin{align}\label{inequality-2}
    -\ideg(Y\, {}^L\!\otimes_A X) \leqslant \extreg(X)-\ideg(Y).
    \end{align}
        \item [(2)] For any $0\neq X\in \mathbf{D}^-(\Gr A)$, $0\neq Y\in \mathbf{D}^-(\gr A^{o})$,
    \begin{align} \label{inequality-3}
    -\ideg(Y\, {}^L\!\otimes_A X) \leqslant \extreg(Y)-\ideg(X).
    \end{align}
    \end{itemize}
\end{proposition}

\begin{proof}
   (1) If $\extreg(X)=+\infty$ or $-\ideg(Y)=+\infty$, the inequality holds trivially.
   
   Let $\extreg(X)=p < \infty $ and $-\ideg(Y)=q < \infty$. 
 
It follows from Lemma \ref{notations of X} that $p=\sup\{m-\lxm{m}{X}\mid P^{-m}\neq 0\}$. Thus, for any $m\in \mathbb{Z}$ such that $P^{-m}\neq 0$, $p\geqslant m-\lxm{m}{X}$. 

By definition, $-q \leqslant -n + \ideg(H^{-n}(Y))$ for all  $n\in \mathbb{Z}$ such that $H^{-n}(Y)\neq 0$.

Now consider the double complex $C^{\bullet \bullet}$ in the third quadrant given by 
$$C^{-m,-n}= Y^{-n}\otimes_A P^{-m}.$$
Then there is a convergent spectral sequence
$${}^I\!E_2^{-m,-n}=\Tor_m^A(H^{-n}(Y),X)\Rightarrow \Tor_{m+n}^A(Y,X).$$

For any $m,n \in \mathbb{Z}$ with $P^{-m}\neq 0$ and $H^{-n}(Y)\neq 0$,
 $$\ideg(H^{-n}(Y)\otimes_A P^{-m})\geqslant n-q+\lxm{m}{X}\geqslant n-q+m-p.$$
  As $\Tor_m^A(H^{-n}(Y),X)$ is a subquotient of $H^{-n}(Y)\otimes_A P^{-m}$, 
  $$\ideg(\Tor_m^A(H^{-n}(Y),X))\geqslant \ideg(H^{-n}(Y)\otimes_A P^{-m} )\geqslant n-q+m-p.$$
By the convergence of the spectral sequence,  $$\ideg(\Tor_{m+n}^A(Y,X))\geqslant -q+n+m-p.$$  Hence
\begin{align*}
   & -\ideg(Y\, {}^L\!\otimes_A X)\\
   =&-\inf\{-m-n + \ideg(\Tor_{m+n}^A(Y,X))\mid \Tor_{m+n}^A(Y,X)\neq 0\}\\
   \leqslant& p+q.
\end{align*}

(2) The proof is similar to (1), replacing the double complex in (1) by $C^{-n,-m}=Q^{-n} \otimes_A X^{-m}$ where $Q^{\bullet}$ is a minimal graded projective resolution of $Y$.
\end{proof}

 There are more inequalities as given in the following corollary for the supremum  or infimum degrees of $R\gHom_A(X,Y)$ and $Y\, {}^L\!\otimes_A X$ by using the isomorphism $\D(Y\, {}^L\!\otimes_A X)\cong R\gHom_A(X,\D(Y))$.

\begin{corollary}\label{sdeg-ideg inequality coro}
    Let $A$ be a noetherian $\mathbb{N}$-graded algebra. Then the following hold.
    \begin{itemize}
        \item [(1)] For any $0\neq X\in \mathbf{D}^-(\gr A)$, $0\neq Y\in \mathbf{D}^-(\Gr A^o)$,  
     $$ \sdeg(Y\, {}^L\!\otimes_A X )\leqslant \Extreg(X)+\sdeg(Y).$$
           
        \item [(2)] For any $0\neq X\in \mathbf{D}^-(\Gr A)$, $0\neq Y \in \mathbf{D}^-(\gr A^o)$, 
         $$ \sdeg(Y\, {}^L\!\otimes_A X )\leqslant \Extreg(Y)+\sdeg(X).$$

        \item [(3)] For any $0\neq X\in \mathbf{D}^-(\Gr A)$ with torsion cohomologies, $0\neq Y \in \mathbf{D}^-_{lf}(\Gr A^o)$,   
         $$-\ideg(Y\, {}^L\!\otimes_A X )\leqslant -\ideg(X)+\extreg(Y).$$

        \item [(4)] For any $0\neq X\in \mathbf{D}^-_{lf}(\Gr A)$, $0\neq Y \in \mathbf{D}^-(\Gr A^o)$ with torsion cohomologies,   
         $$-\ideg(Y\, {}^L\!\otimes_A X )\leqslant -\ideg(Y)+\extreg(X).$$
         
    \item [(5)] If $A$ has a balanced dualizing complex, then for any $0\neq X \in \mathbf{D}^-_{fg}(\Gr A)$ and $0\neq Y \in \mathbf{D}^+_{fg}(\Gr A)$,  
$$-\ideg(\D(Y)\, {}^L\!\otimes_A X )\leqslant \cmreg(X)+\Exreg(Y).$$

   \item [(6)] If $A$ has a balanced dualizing complex, then for any $0\neq X \in \mathbf{D}^-_{fg}(\Gr A)$ and $0\neq Y \in \mathbf{D}^-(\gr A^o)$ with torsion cohomologies,
   $$ \sdeg(Y\, {}^L\!\otimes_A X )\leqslant \CMreg(X)+\Extreg(Y).$$

        \item [(7)] For any $0\neq X\in \mathbf{D}^-(\gr A)$, $0\neq Y \in \mathbf{D}^-(\Gr A^o)$, 
        $$\sdeg(R\gHom_A(X,\D(Y)))\leqslant \extreg(X)-\ideg(Y).$$
            
        \item [(8)] For any $0\neq X\in \mathbf{D}^-(\Gr A)$, $0\neq Y \in \mathbf{D}^-(\gr A^o)$, 
        $$\sdeg(R\gHom_A(X,\D(Y)))\leqslant -\ideg(X)+\extreg(Y).$$
    \end{itemize}
\end{corollary}
    
\begin{proof}
(1) It follows from  $\D(Y\, {}^L\!\otimes_A X)\cong R\gHom_A(X,\D(Y))$ and \eqref{inequality-1} that
         \begin{align*}
             \sdeg(Y\, {}^L\!\otimes_A X)=&-\ideg(\D(Y\, {}^L\!\otimes_A X))\\
            =&-\ideg(R\gHom_A(X,\D(Y))) \\
            \leqslant & \Extreg(X)-\ideg(\D(Y))\\
            =& \Extreg(X)+\sdeg(Y).
         \end{align*}
       
(2) Similar to the proof of (1) by using the right version of \eqref{inequality-1}, or 
\begin{align*}
     \sdeg(Y\, {}^L\!\otimes_A X)=&-\ideg(\D(Y\, {}^L\!\otimes_A X))\\
            =&-\ideg(R\gHom_A(X,\D(Y))) \\
            \leqslant & \Extreg(Y)+\sdeg(X) \quad (\textrm{by \eqref{inequality-2-1}})
         \end{align*}
         
 (3) This is proved by using Proposition \ref{sdeg-ideg inequality 1-1}. 
 \begin{align*}
              -\ideg(Y\, {}^L\!\otimes_A X)=&\sdeg(\D(Y\, {}^L\!\otimes_A X))\\
             =&\sdeg(R\gHom_A(X,\D(Y))) \\
             \leqslant & -\ideg(X)+\Exreg(\D(Y)) \quad (\textrm{by \eqref{inequality-3-1}})\\
             =& -\ideg(X) + \extreg(Y) \quad (\textrm{by \eqref{Ex-Ext and ex-ext-4}}).
         \end{align*}

  (4) Similar to the proof of (3) by using the right version of  \eqref{inequality-3-1}.
  
(5) It follows from $\D(\D(Y)\, {}^L\!\otimes_A X) \cong R\gHom_A(X,\D(\D(Y)))$ that
\begin{align*}
    -\ideg(\D(Y)\, {}^L\!\otimes_A X )
    =&\sdeg(R\gHom_A(X,\D(\D(Y))))\\
    =&\sdeg(R\gHom_A(X,Y))\\
    \leqslant& \cmreg(X)+\Exreg(Y) \quad (\textrm{by \eqref{inequality-3-2}}).
\end{align*}

(6) By \eqref{inequality-2-2},
$$\sdeg(Y\, {}^L\!\otimes_A X )=\sdeg(\D(Y\, {}^L\!\otimes_A X))
    \leqslant \CMreg(X)+\Extreg(Y).$$

(7) By \eqref{inequality-2},
$$\sdeg(R\gHom_A(X,\D(Y)))=-\ideg(Y\, {}^L\!\otimes_A X)
             \leqslant \extreg(X)-\ideg(Y).$$
         
 (8) By \eqref{inequality-3},
  $$\sdeg(R\gHom_A(X,\D(Y)))=-\ideg(Y \, {}^L\!\otimes_A X) 
             \leqslant \extreg(Y)-\ideg(X).$$      
\end{proof}

\subsection{Some equalities}
In this subsection we explore the conditions under which the inequalities \eqref{inequality-1}, \eqref{inequality-3-1}, \eqref{inequality-3-2}, \eqref{inequality-2-1}, \eqref{inequality-2-2}, \eqref{inequality-2} and \eqref{inequality-3} are equalities in Propositions \ref{sdeg-ideg inequality 1}, \ref{sdeg-ideg inequality 1-1}, \ref{sdeg-ideg inequality 1-4} and  \ref{sdeg-ideg inequality 2}. We assume that $A_0$ is semisimple, and prove two technical lemmas first. 
For any $0\neq X\in \mathbf{D}^-(\gr A)$, if $p:=-\ideg(X)$ and $P^{\bullet}$ is the minimal graded projective resolution of $X$, then Lemma \ref{pre for equality 1} says that at least one generator in a minimal generating subset  of some $P^{-\alpha}$ is a $(-\alpha)$-cocycle of degree $\alpha - p$.  Lemma \ref{pre for equality 2} says that if $p=-\ideg(X) =-\ideg(e_k X)$ for some $1\leqslant k \leqslant n$ then there is some  $(-\alpha)$-cocycle $y$ of the minimal degree in $P^{-\alpha}$ such that $e_k y \neq 0$. 
Lemmas \ref{pre for equality 1} and \ref{pre for equality 2} play a key role in proving Propositions \ref{sdeg-ideg equality 1} and \ref{sdeg-ideg equality 2}. For the notation see \S \ref{notation-e-i}.
\begin{lemma}\label{pre for equality 1}
Suppose $A$ is a noetherian $\mathbb{N}$-graded algebra with $A_0$ semisimple. Let $0\neq X\in \mathbf{D}^-(\gr A)$ and $P^{\bullet}$ be the minimal graded projective resolution of $X$.

If $p=-\ideg(X)$, then there  is some $P^{-\alpha}$  such that $P^{-\alpha}$ has an indecomposable direct summand generated by an element of degree $(\alpha-p)$ which is contained in $\ker d^{-\alpha}_{P^{\bullet}}$.  
\end{lemma}
\begin{proof} 
  It follows from Proposition \ref{CM and Ext} and Lemma \ref{notations of X} that 
  $$p=-\ideg(X)=\extreg(X)=\sup\{m-\lxm{m}{X}\mid P^{-m}\neq 0\}.$$ 
  Without loss of generality, we may assume that $X^i=0$ for all $i>0$ and $\sup(X)=0$. 
Then either $0-\lxm{0}{X}=p$ or $0-\lxm{0}{X}<p$ and  $\alpha-\lxm{\alpha}{X}=p$ for some $\alpha > 0$. 

Suppose $0-\lxm{0}{X}=p$. It follows from $P^0=\ker d^{0}_{P^{\bullet}}$ that $Ae_k(p-0)\subseteq P^0=\ker d^{0}_{P^{\bullet}}$ for some $1\leqslant k \leqslant n$. In this case, $Ae_k(p-0)$ is an indecomposable direct summand of $P^0$ generated by an element of degree $(0-p)$.

Suppose $0-\lxm{0}{X}<p$ and $\alpha-\lxm{\alpha}{X}=p$ for some $\alpha > 0$.  We may assume that $\alpha$ is the minimal positive integer such that $\alpha-\lxm{\alpha}{X}=p$, 
 that is, $\alpha'-\lxm{\alpha'}{X}<p$ for any $0 \leqslant \alpha' < \alpha$. 
In any minimal generating set of $P^{-\alpha}$, there is an element of degree  $\lxm{\alpha}{X}=\alpha-p$. So,
there is some $k$ with $1\leqslant k \leqslant n$ such that $Ae_k(p-\alpha)$ is a direct summand of $P^{-\alpha}$, i.e.,
$$P^{-\alpha}=\mathop{\bigoplus}\limits_{i} \big(\mathop{\bigoplus}\limits_{j} Ae_i(-s_{\alpha}^{i,j})\big)=Ae_k(p-\alpha) \oplus P'$$
for some $P'$. If $Ae_k(p-\alpha)\nsubseteq \ker d^{-\alpha}_{P^{\bullet}}$, then it follows from $d^{-\alpha}_{P^{\bullet}}(Ae_k(p-\alpha))\subseteq A_{\geqslant1}P^{-\alpha+1}$ and $\alpha -1 < \alpha$ that 
$$\alpha-p\geqslant \lxm{\alpha-1}{X}+1>\alpha-1-p+1=\alpha-p,$$
which is impossible.
Therefore, $Ae_k(p-\alpha) \subseteq \ker d^{-\alpha}_{P^{\bullet}}$.
\end{proof}

\begin{lemma}\label{pre for equality 2} Suppose $A$ is a noetherian $\mathbb{N}$-graded algebra with $A_0$ semisimple. Let $0\neq X\in \mathbf{D}^-(\gr A)$ and $P^{\bullet}$ be the minimal graded projective resolution of $X$.

If $p=-\ideg(X)=-\ideg(e_kX)$ for some  $1\leqslant k \leqslant n$, then there exists $\beta \in \mathbb{Z}$ with $H^{-\beta}(e_kP^{\bullet})\neq 0$ such that 
$$\ideg(\ker d_{P^{\bullet}}^{-\beta})=\beta-p\, \textrm{ and } \, e_k(\ker d_{P^{\bullet}}^{-\beta})_{\beta-p} \neq 0.$$
\end{lemma}
\begin{proof}
     In general, for any $1 \leqslant k \leqslant n$,
     \begin{align*}
         &-\ideg(e_kX)\\
         =& -\inf \{-m+ \ideg(H^{-m}(e_kP^{\bullet}))\mid  H^{-m}(e_kP^{\bullet})\neq 0, m\in \mathbb{Z} \}\\
         \leqslant & -\inf \{ -m+ \ideg(\ker d_{e_kP^{\bullet}}^{-m}) \mid  H^{-m}(e_kP^{\bullet})\neq 0, m\in \mathbb{Z}\}\\
         \leqslant & -\inf \{ -m+ \ideg(P^{-m})\mid H^{-m}(e_kP^{\bullet})\neq 0, m\in \mathbb{Z} \}\\
         \leqslant & -\inf\{-m+ \ideg(P^{-m})\mid P^{-m}\neq 0, m\in \mathbb{Z}\}\\
         =& -\inf\{-m+ \lxm{m}{X} \mid P^{-m}\neq 0, m\in \mathbb{Z}\}\\
         =&\extreg(X) \quad (\textrm{by Lemma } \ref{notations of X})\\
         =&-\ideg(X) \quad (\textrm{by Proposition } \ref{CM and Ext}).  
     \end{align*}
 If $p=-\ideg(X)=-\ideg(e_kX)$ for some  $1\leqslant k \leqslant n$, then all the inequalities above are in fact equalities.
 Therefore, there exists some $\beta \in \mathbb{Z}$ with $ H^{-\beta}(e_kP^{\bullet})\neq 0$ such that
     $$p=\beta-\ideg( \ker d_{e_kP^{\bullet}}^{-\beta})=\beta-\ideg(e_k \ker d_{P^{\bullet}}^{-\beta}).$$
     Consequently,
     \begin{align*}
      \beta-p&=\ideg(e_k \ker d_{P^{\bullet}}^{-\beta})\\ &\geqslant \ideg(\ker d_{P^{\bullet}}^{-\beta}) \\
      &\geqslant \ideg(P^{-\beta})=\lxm{\beta}{X}\geqslant \beta-p.
     \end{align*}
Hence $\ideg(e_k \ker d_{P^{\bullet}}^{-\beta})=\ideg(\ker d_{P^{\bullet}}^{-\beta})=\beta-p$, and $e_k (\ker d_{P^{\bullet}}^{-\beta})_{\beta-p} \neq 0$. 
\end{proof}

Now we are ready to give some criteria for when \eqref{inequality-1}-\eqref{inequality-3} in \S \ref{degree-of-Hom-Tensor} are equalities under the condition that $A_0$ is semisimple.  Note that if $A_0$ is semisimple,  then $-\ideg(X) = \extreg(X)$ by Proposition \ref{CM and Ext}. If $A$ is connected graded, then the conditions (2) in all the following propositions and corollaries in this subsection are trivially true, so  the identities in (1) are always true in the connected graded case, which are also new.

Propositions \ref{sdeg-ideg equality 1} and \ref{sdeg-ideg equality 2} give equivalent conditions for \eqref{inequality-2} and \eqref{inequality-3} being equalities if $A_0$ is semisimple.

\begin{proposition}\label{sdeg-ideg equality 1}
   Suppose $A$ is a noetherian $\mathbb{N}$-graded algebra with $A_0$ semisimple. For any $0\neq Y \in \mathbf{D}^b(\gr {A^o})$, the following are equivalent.
 \begin{itemize}
     \item [(1)] $\ideg(Y\, {}^L\!\otimes_A X) =  \ideg(X)+\ideg(Y)$ for all $0\neq X\in \mathbf{D}^b(\gr A)$.
     \item [(2)] $\ideg(Ye_i)=\ideg(Y)$ for all $1\leqslant i\leqslant n$.
 \end{itemize} 
\end{proposition}
\begin{proof}
    (1) $\Rightarrow$ (2) It follows by taking $X=Ae_i$.

    (2) $\Rightarrow$ (1)
Without loss of generality, let $X^i=0, Y^i=0$ for all $i>0$ and $\sup(X)=\sup(Y)=0$. 
Since $0\neq X\in \mathbf{D}^b(\gr A)$ and $0 \neq Y \in \mathbf{D}^b(\gr {A^o})$, $-\ideg(X) \neq  \pm \infty$ and $-\ideg(Y) \neq \pm \infty$. Let $p_1= -\ideg(X)$ and $p_2=-\ideg(Y)$. 
Let $P^{\bullet}$ and $Q^{\bullet}$ be the minimal graded projective resolutions of $X$ and $Y$ respectively.

Note that $Y\, {}^L\!\otimes_A X\cong Q^{\bullet} \otimes_{A} P^{\bullet}$ 
    where $(Q^{\bullet} \otimes_{A} P^{\bullet})^m= \mathop{\bigoplus}\limits_{q}Q^{q} \otimes_{A} P^{m-q},$ and its $m$th-differential is the morphism $$\partial^m(x \otimes y)=d_{Q^{\bullet}}(x) \otimes y +(-1)^{q}x\otimes d_{P^{\bullet}}(y).$$
It follows from Proposition \ref{sdeg-ideg inequality 2} (1) that 
     \begin{align*}
         &-\ideg(Y\, {}^L\!\otimes_A X)\\
         &=-\inf\{m+\ideg(H^m(Q^{\bullet}\otimes_{A} P^{\bullet}))\mid \Tor_{-m}^A(Y,X)\neq 0\}\\
         &\leqslant \extreg(X)-\ideg(Y)= p_1+p_2.
     \end{align*}
Hence, for any $m\in \mathbb{Z}$ such that $\Tor_{-m}^A(Y,X)\neq 0$, $$\ideg(H^{m}(Q^{\bullet}\otimes_{A} P^{\bullet})) \geqslant -p_1-p_2-m.$$
To prove (1), it suffices to show that there is some $\alpha \in \mathbb{Z}$ with
$\Tor_{-\alpha}^A(Y,X)\neq 0$ such that
$$\ideg(H^{\alpha}(Q^{\bullet}\otimes_{A} P^{\bullet}))\leqslant -p_1-p_2-\alpha.$$

On one hand, it follows from Lemma \ref{pre for equality 1} that there is some $\alpha_1 \in \mathbb{Z}$ such that $P^{-{\alpha}_1}$ has an indecomposable direct summand generated by an element of degree $(\alpha_1-p_1)$ contained in $\ker d^{-\alpha_1}_{P^{\bullet}}$.  
So, there is some
$k$ with $1\leqslant k\leqslant n$ such that $Ae_k(p_1-\alpha_1)\subseteq \ker d^{-\alpha_1}_{P^{\bullet}}$.

On the other hand, by the assumption in (2) and the right version of Lemma \ref{pre for equality 2}, there exists $\alpha_2 \in \mathbb{Z}$ with $H^{-\alpha_2}(Q^{\bullet}e_k)\neq 0$ such that 
$$\ideg(\ker d_{Q^{\bullet}}^{-{\alpha_2}})=\alpha_2-p_2\, \textrm{ and } \, (\ker d_{Q^{\bullet}}^{-{\alpha_2}})_{\alpha_2-p_2}\,e_k \neq 0.$$
There is some $0\neq y \in (\ker d_{Q^{\bullet}}^{-\alpha_2})_{\alpha_2-p_2} \subseteq Q^{-\alpha_2}$ such that 
$$ 0 \neq y \otimes_A e_k(p_1-\alpha_1) \subseteq Q^{-\alpha_2} \otimes_A P^{-\alpha_1}.$$

Then $0\neq y \otimes_A e_k(p_1-\alpha_1) \in \ker \partial^{\alpha}$, where $\alpha=-\alpha_2-\alpha_1$. Hence $$\ideg(\ker \partial^{\alpha})\leqslant \alpha_2-p_2 + \alpha_1 -p_1 = -\alpha-p_1-p_2.$$ 
Since $P^{\bullet}$ and $Q^{\bullet}$ are minimal graded projective resolutions, 
$$\im \partial^{\alpha-1} \subseteq \mathop{\bigoplus}\limits_{q}(Q^{q+1}A_{\geqslant 1} \otimes_{A} P^{\alpha-1-q}+Q^{q} \otimes_{A} A_{\geqslant 1}P^{\alpha-q}) \subseteq (Q^{\bullet}\otimes_{A} P^{\bullet})^{\alpha}.$$ 
It follows that 
\begin{align*}
&\ideg(\im \partial^{\alpha-1}) \\
\geqslant& \min\{\lxm{-q-1}{Y}+1 + \lxm{-\alpha+1+q}{X}, \lxm{-q}{Y}+1 + \lxm{-\alpha+q}{X}  \}\\
\geqslant& \min\{-q-1-p_2+1 -\alpha+1+q-p_1, -q-p_2+1 -\alpha+q-p_1 \}\\
=& -p_1-p_2-\alpha+1.
\end{align*}
Hence
    $\ideg (H^{\alpha}(Q^{\bullet}\otimes_{A} P^{\bullet}))\leqslant -p_1-p_2-\alpha.$
    
Therefore, 
$-\ideg(Y\, {}^L\!\otimes_A X)=p_1+p_2=-\ideg(X)-\ideg(Y),$ and the proof of (2) $\Rightarrow$ (1) is finished.
\end{proof}

\begin{proposition}\label{sdeg-ideg equality 2}
   Suppose $A$ is a noetherian $\mathbb{N}$-graded algebra with $A_0$ semisimple. 
   For any $0\neq X\in \mathbf{D}^b(\gr A)$, the following are equivalent.
 \begin{itemize}
     \item [(1)] $\ideg(Y\, {}^L\!\otimes_A X) =  \ideg(X)+\ideg(Y)$ for all $0\neq Y \in \mathbf{D}^b(\gr {A^o})$.
     \item [(2)] $\ideg(e_iX)=\ideg(X)$ for all $1\leqslant i\leqslant n$.
 \end{itemize} 
\end{proposition}
\begin{proof} (1) $\Rightarrow$ (2) It follows by taking $Y=e_iA$  in (1).

(2) $\Rightarrow$ (1) Similar to the proof of (2) $\Rightarrow$ (1) in Proposition \ref{sdeg-ideg equality 1}.
\end{proof}

Propositions \ref{sdeg-ideg equality 3} and \ref{sdeg-ideg equality 4} give criteria for when \eqref{inequality-1} is an equality under the condition that $A_0$ is semisimple. 

\begin{proposition}\label{sdeg-ideg equality 3}
 Suppose $A$ is a noetherian $\mathbb{N}$-graded algebra with $A_0$ semisimple. For any $0\neq Y \in \mathbf{D}^b(\gr {A})$, the following are equivalent.
 \begin{itemize}
    \item [(1)] $-\ideg(R\gHom_A(X,Y))=\Extreg(X)-\ideg(Y)$ for all $0\neq X\in \mathbf{D}^b(\gr A)$ with $\pdim(X)<\infty$.
    \item [(2)] 
    $\ideg(e_iY)=\ideg(Y)$ for all $1\leqslant i\leqslant n$.
 \end{itemize}
\end{proposition} 
\begin{proof}
    (1) $\Rightarrow$ (2) It follows by taking $X=Ae_i$  in (1) that
    $$-\ideg(R\gHom_A(Ae_i,Y))=\Extreg(Ae_i)-\ideg(Y).$$
    Since $R\gHom_A(Ae_i,Y)\cong e_iY$ as graded $k$-spaces and $\Extreg(Ae_i)=0$, $$\ideg(e_iY)=0+\ideg(Y)=\ideg(Y).$$ 

    (2) $\Rightarrow$ (1) Let $Z=R\gHom_A(X,A)$. Since $\pdim(X)<\infty$, $0\neq Z \in \mathbf{D}^b(\gr A^o)$ and 
      \begin{align}\label{isomorphism}
          R\gHom_A(X,Y) \cong R\gHom_A(X,A) \, {}^L\!\otimes_A Y \cong Z \, {}^L\!\otimes_A Y. 
      \end{align}
      It follows from Proposition \ref{sdeg-ideg equality 2} that 
      $$ \ideg(Z \, {}^L\!\otimes_A Y)=\ideg(Y)+\ideg(Z).$$
      By Theorem \ref{Extreg}, 
      $\Extreg(X)=-\ideg(R\gHom_A(X,A))=-\ideg(Z).$

      Thus, 
      \begin{align*}
          -\ideg(R\gHom_A(X,Y))=&-\ideg(Z \, {}^L\!\otimes_A Y)\\
          =&-\ideg(Y)-\ideg(Z)\\
          =&-\ideg(Y)+\Extreg(X).
      \end{align*}
The proof of (2) $\Rightarrow$ (1) is finished.
\end{proof}

\begin{proposition}\label{sdeg-ideg equality 4}
 Suppose $A$ is a noetherian $\mathbb{N}$-graded algebra with $A_0$ semisimple. For any $0\neq X \in \mathbf{D}^b(\gr A)$ with $\pdim(X)<\infty$, the following are equivalent.
 \begin{itemize}
     \item [(1)] $-\ideg(R\gHom_A(X,Y))=\Extreg(X)-\ideg(Y)$ for all $0\neq Y\in \mathbf{D}^b(\gr A)$.
    \item [(2)] $\Extreg(X)=-\ideg(R\gHom_A(X,Ae_i))$ for all $1\leqslant i\leqslant n$. 
 \end{itemize}
\end{proposition} 

\begin{proof} (1) $\Rightarrow$ (2) It follows by taking $Y=Ae_i$.

(2) $\Rightarrow$ (1) Let $Z=R\gHom_A(X,A)$. Since $\pdim(X)<\infty$, $0\neq Z \in \mathbf{D}^b(\gr A^o)$. It follows from Theorem \ref{Extreg} and $R\gHom_A(X,Ae_i)\cong Z \otimes_A Ae_i\cong Ze_i$ that 
$$ -\ideg(Z)=\Extreg(X)\stackrel{(2)}{=}-\ideg(R\gHom_A(X,Ae_i))=-\ideg(Ze_i).$$ 
 Hence by Proposition \ref{sdeg-ideg equality 1} and \eqref{isomorphism}, 
 \begin{align*}
      -\ideg(R\gHom_A(X,Y))=&-\ideg(Z \, {}^L\!\otimes_A Y)\\
      =&-\ideg(Y)-\ideg(Z)\\
      =&-\ideg(Y)+\Extreg(X).
 \end{align*}
 The proof of (2) $\Rightarrow$ (1) is finished.
\end{proof}

Propositions \ref{sdeg-ideg equality 7} and \ref{sdeg-ideg equality 8} tell us when \eqref{inequality-3-1} is an equality.

\begin{proposition}\label{sdeg-ideg equality 7}
    Suppose $A$ is a noetherian $\mathbb{N}$-graded algebra with $A_0$ semisimple. For any $0\neq Y \in \mathbf{D}^b(\Gr {A})$ with $\D(Y) \in \mathbf{D}^b(\gr {A^o})$, the following are equivalent.
 \begin{itemize}
     \item [(1)] $\sdeg(R\gHom_A(X,Y)) = -\ideg(X)+\Exreg(Y)$ for all $0\neq X\in \mathbf{D}^b(\gr A)$.
     \item [(2)]  $\sdeg(e_iY)=\sdeg(Y)$ for all $1\leqslant i\leqslant n$.
 \end{itemize} 
\end{proposition}
 \begin{proof}  Since $\D(Y) \in \mathbf{D}^b(\gr {A^o})$, $Y\in \mathbf{D}^b_{lf}(\Gr {A})$.
 
     (1) $\Rightarrow$ (2) Taking $X=Ae_i$, then $$\sdeg(e_iY)=\Exreg(Y)=\extreg(\D(Y))=-\ideg(\D(Y))=\sdeg(Y).$$

     (2) $\Rightarrow$ (1) Since $-\ideg(\D(Y)e_i)=-\ideg(\D(e_iY))=\sdeg(e_iY)\stackrel{(2)}{=}\sdeg(Y)=-\ideg(\D(Y))$, it follows from Proposition \ref{sdeg-ideg equality 1} that 
       \begin{align*}
       -\ideg(\D(Y)\, {}^L\!\otimes_A X) =& -\ideg(X)-\ideg(\D(Y)).
       \end{align*} 
      
    Since $\D(\D(Y)\, {}^L\!\otimes_A X)
     \cong R\gHom_A(X,Y)$, 
     \begin{align*}
         \sdeg(R\gHom_A(X,Y))=&-\ideg(\D(Y)\, {}^L\!\otimes_A X)\\
         =&-\ideg(X)-\ideg(\D(Y))\\
         =&-\ideg(X)+\extreg(\D(Y))\\
         =&-\ideg(X)+\Exreg(Y).
     \end{align*}
\end{proof}

\begin{proposition}\label{sdeg-ideg equality 8}
    Suppose $A$ is a noetherian $\mathbb{N}$-graded algebra with $A_0$ semisimple.  For any $0\neq X \in \mathbf{D}^b(\gr A)$, the following are equivalent.
 \begin{itemize}
      \item [(1)] $\sdeg(R\gHom_A(X,Y)) = -\ideg(X)+\Exreg(Y)$ for all $0\neq Y \in \mathbf{D}^b(\Gr {A})$ with $\D(Y) \in \mathbf{D}^b(\gr {A^o})$.
     \item [(2)] $\ideg(e_iX)=\ideg(X)$ for all $1 \leqslant i\leqslant n$.
 \end{itemize} 
\end{proposition}
 \begin{proof}
     (1) $\Rightarrow$ (2)  By  taking $Y=\D(e_iA)$, it follows from $\Exreg(Y)=\extreg(\D(Y))=0$.

     (2) $\Rightarrow$ (1) Since $\D(Y) \in \mathbf{D}^b(\gr {A^o})$, $Y \in \mathbf{D}^b_{lf}(\Gr {A})$ and $\D(\D(Y)\cong Y$.
      Then $\D(\D(Y)\, {}^L\!\otimes_A X)\cong R\gHom_A(X,\D(\D(Y))) \cong R\gHom_A(X,Y)$. 
     Hence
     \begin{align*}
         \sdeg(R\gHom_A(X,Y))=&-\ideg(\D(Y)\, {}^L\!\otimes_A X)\\
         =&-\ideg(X)-\ideg(\D(Y))  \quad (\textrm{by Proposition } \ref{sdeg-ideg equality 2})\\
         =&-\ideg(X)+\extreg(\D(Y))\\
         =&-\ideg(X)+\Exreg(Y).
     \end{align*}
\end{proof}

Propositions \ref{sdeg-ideg equality 5} and \ref{sdeg-ideg equality 6} concern when \eqref{inequality-2-1} is an equality.
\begin{proposition}\label{sdeg-ideg equality 5}
    Suppose $A$ is a noetherian $\mathbb{N}$-graded algebra with $A_0$ semisimple. For any $0\neq X \in \mathbf{D}^b(\Gr {A})$ with $\D(X) \in \mathbf{D}^b(\gr {A^o})$, the following are equivalent.
 \begin{itemize}
     \item [(1)] $-\ideg(R\gHom_A(X,\D(Y))) = \sdeg(X)+\Extreg(Y)$ for all $0\neq Y\in \mathbf{D}^b(\gr A^o)$ with $\pdim Y < \infty$.
     \item [(2)]  $\sdeg(e_iX)=\sdeg(X)$ for all $1\leqslant i\leqslant n$.
 \end{itemize} 
\end{proposition}
 \begin{proof}
    (1) $\Rightarrow$ (2) It follows by taking $Y=e_iA$.

    (2) $\Rightarrow$ (1) Since $-\ideg(\D(X)e_i)=-\ideg(\D(e_iX))=\sdeg(e_iX)\stackrel{(2)}{=}\sdeg(X)=-\ideg(\D(X))$, it follows from the right version of Proposition \ref{sdeg-ideg equality 3} that 
   \begin{align*}
   -\ideg(R\gHom_A(X,\D(Y)))=&-\ideg(R\gHom_{A^o}(Y,\D(X)))\\
   =&\Extreg(Y)-\ideg(\D(X))\\
   =&\Extreg(Y)+\sdeg(X). 
   \end{align*} 
\end{proof}

\begin{proposition}\label{sdeg-ideg equality 6}
    Suppose $A$ is a noetherian $\mathbb{N}$-graded algebra with $A_0$ semisimple. For any $0\neq Y\in \mathbf{D}^b(\gr A^o)$ with $\pdim Y < \infty$, the following are equivalent.
 \begin{itemize}
     \item [(1)] $-\ideg(R\gHom_A(X,\D(Y))) = \Extreg(Y)+\sdeg(X)$ for all $0\neq X \in \mathbf{D}^b(\Gr {A})$ with $\D(X) \in \mathbf{D}^b(\gr {A^o})$.
      \item [(2)] $\Extreg(Y)=-\ideg(R\gHom_{A^o}(Y,e_iA))$ for all $1\leqslant i\leqslant n$.
 \end{itemize} 
\end{proposition}
 \begin{proof}
    (1) $\Rightarrow$ (2) It follows by  taking $X=\D(e_iA)$.

    (2) $\Rightarrow$ (1) 
    It follows from the right version of Proposition \ref{sdeg-ideg equality 4} that 
    \begin{align*}
    -\ideg(R\gHom_A(X,\D(Y)))=&-\ideg(R\gHom_{A^o}(Y,\D(X)))\\
    =&\Extreg(Y)-\ideg(\D(X))\\
    =&\Extreg(Y)+\sdeg(X). 
    \end{align*}  
\end{proof}

Corollaries \ref{sdeg-ideg equality 9} and \ref{sdeg-ideg equality 10} concern when \eqref{inequality-3-2} and \eqref{inequality-2-2} are equalities respectively.

\begin{corollary}\label{sdeg-ideg equality 9}
  Suppose $A$ is a noetherian $\mathbb{N}$-graded algebra with $A_0$ semisimple, and $A$ has a balanced dualizing complex. For any $X \in \mathbf{D}^b(\fd A)$, the following are equivalent. 
    \begin{itemize}
      \item [(1)] $\sdeg(R\gHom_A(X,Y))=\cmreg(X)+\Exreg(Y)$ for all $0\neq Y \in \mathbf{D}^b(\Gr {A})$ with $\D(Y) \in \mathbf{D}^b(\gr {A^o})$.
     \item [(2)] $\ideg(e_iX)=\ideg(X)$ for all $1 \leqslant i\leqslant n$.
  \end{itemize}
\end{corollary}
\begin{proof}
    Since $X \in \mathbf{D}^b(\fd A)$, $R\Gamma_A(X)\cong X$ and $\cmreg(X)=-\ideg(X)$. It follows from Proposition \ref{sdeg-ideg equality 8} that the conclusion holds.
\end{proof}

\begin{corollary}\label{sdeg-ideg equality 10}
  Suppose $A$ is a noetherian $\mathbb{N}$-graded algebra with $A_0$ semisimple and $A$ has a balanced dualizing complex. For any $Y\in \mathbf{D}^b(\fd A^o)$ with $\pdim Y <\infty$, the following are equivalent.
  \begin{itemize}
      \item [(1)] $-\ideg(R\gHom_A(X,\D(Y)))=\CMreg(X)+\Extreg(Y)$ for all $0\neq X \in \mathbf{D}^b(\gr A)$.
      \item [(2)] $\Extreg(Y)=-\ideg(R\gHom_{A^o}(Y,e_iA))$ for all $1\leqslant i\leqslant n$.
  \end{itemize}
\end{corollary}
\begin{proof}
     By Lemma \ref{iso about local cohomological}, $R\gHom_A(X,\D(Y))\cong R\gHom_A(R\Gamma_A(X),\D(Y))$. Since $X \in \mathbf{D}^b(\gr A)$ and $A$ has a balanced dualizing complex, $\D(R\Gamma_A(X))\in \mathbf{D}^b(\gr A^o)$. 
     
     (1) $\Rightarrow$ (2) it follows by the proof of (1) $\Rightarrow$ (2) of Proposition \ref{sdeg-ideg equality 6} and taking $R\Gamma_A(X)=\D(e_iA)$.

     (2) $\Rightarrow$ (1) it follows from the proof of (2) $\Rightarrow$ (1) of Proposition \ref{sdeg-ideg equality 6} that
     \begin{align*}
         -\ideg(R\gHom_A(X,\D(Y)))=& -\ideg(R\gHom_A(R\Gamma_A(X),\D(Y)))\\
         =& \sdeg(R\Gamma_A(X))+\Extreg(Y)\\
         =&\CMreg(X)+\Extreg(Y).
     \end{align*}    
\end{proof}

\section{Numerical Artin-Schelter regularities}
\subsection{Numerical Artin-Schelter regularity}
Following \cite{KWZ1}, a notion of numerical Artin-Schelter regularity is given first for noetherian $\mathbb{N}$-graded algebras in this subsection.  Then we generalize \cite[Theorems 2.5 and 2.6]{Jo4} to $\mathbb{N}$-graded algebras with balanced dualizing complexes, and develop more similar (or more stronger) relations between the regularities (say, see Theorem \ref{relation between CM and Tor} (4) and Proposition \ref{relation between CM and Tor-p}). For any noetherian $\mathbb{N}$-graded algebra $A$ with a balanced dualizing complex, we prove that the numerical AS-regularity $\ASreg(A)=0$ if and only if
 $\CMreg(X)= \Torreg(X)+\CMreg(A)$ for all $0\neq X \in \mathbf{D}^b(\gr A)$; if and only if
 $\Torreg(X)= \CMreg(X)+\Torreg(S)$ for all $0\neq X \in \mathbf{D}^b(\gr A)$. 

\begin{definition} (see \cite[Definition 0.6]{KWZ1})
 The Artin-Schelter regularity ({\it AS regularity} for short) of a noetherian $\mathbb{N}$-graded algebra $A$ is defined to be $$\ASreg(A) =\CMreg(A) + \Torreg(S).$$
\end{definition}
If $A$ is a noetherian $\mathbb{N}$-graded algebra with $A_0$ semisimple and $A$ has a balanced dualizing complex, then $\ASreg(A)=\Exreg(A)+\Torreg(S)$ by Theorem \ref{ex}.

Note that $\ASreg(A)$ runs over all positive integers \cite[Lemma 5.6]{KWZ1}.

Theorem \ref{relation between CM and Tor} (1) and (2) in the following are  generalization of \cite[Theorems 2.5 and 2.6]{Jo4}. The proofs here rely on Proposition  \ref{sdeg-ideg inequality 1}. 
\begin{theorem}\label{relation between CM and Tor}
 Suppose $A$ is a noetherian $\mathbb{N}$-graded algebra with a balanced dualizing complex $R$. Then
   \begin{itemize}
       \item [(1)] $\CMreg(X)\leqslant \Torreg(X)+\CMreg(A)$ for all $0\neq X \in \mathbf{D}^-(\gr A)$.
       \item [(2)] $\Torreg(X)\leqslant \CMreg(X)+\Torreg(S)$ for all $0\neq X \in \mathbf{D}^-(\gr A)$.
       \item [(3)] $\ASreg(A)\geqslant 0$.
       \item [(4)] $\CMreg(A)+\exreg({}_AA)\geqslant 0$, and $\CMreg(A)+\exreg(A_A)\geqslant 0$. 
   \end{itemize}
\end{theorem}
\begin{proof}
(1) By Definition \ref{CM-reg-def},
\begin{equation}\label{CMreg-identity}
 \CMreg(X)=\sdeg(R\Gamma_A(X))=-\ideg(\D(R\Gamma_A(X))). 
\end{equation}
In particular, $\CMreg(A)=-\ideg(R)$.
Hence  
\begin{align*}
    \CMreg(X) =&-\ideg(\D(R\Gamma_A(X)))\\
    =&-\ideg(R\gHom_A(X,R))  \quad (\textrm{by Theorem } \ref{LD})\\
    \leqslant& \Extreg(X)-\ideg(R)   \quad (\textrm{by Proposition } \ref{sdeg-ideg inequality 1})\\
    =& \Torreg(X)+\CMreg(A)  \quad (\textrm{by Lemma } \ref{Extreg-Torreg}). 
\end{align*}

     (2) Since $R\gHom_A(X,S)\cong R\gHom_{A^o}(S,\D(R\Gamma_A(X)))$, 
    \begin{align*}
        \Torreg(X)=&-\ideg(R\gHom_A(X,S)) \quad (\textrm{by Lemma } \ref{Extreg-Torreg})\\
        =&-\ideg(R\gHom_{A^o}(S,\D(R\Gamma_A(X))))\\
    \leqslant & \Extreg(S)-\ideg(\D(R\Gamma_A(X))) \quad (\textrm{by Proposition } \ref{sdeg-ideg inequality 1})\\
    =&\Torreg(S)+\CMreg(X) \quad (\textrm{by } \eqref{CMreg-identity}). 
    \end{align*}

   (3) It follows from (1) by taking $X=A$ and $\CMreg(A) < \infty$ that $\Torreg(A) \geqslant 0$.  Then, by (2) $$\ASreg(A)=\CMreg(A) + \Torreg(S) \geqslant \Torreg(A) \geqslant 0.$$ 

   (4) By Definitions \ref{Defi-CM} and \ref{defi-bdc}, $$\CMreg(R_A) = \sdeg(R\Gamma_{A^o}(R_A)) = \sdeg(D(A))=0.$$
   Note that $\Torreg(R_A)=\exreg({}_AA)$ and $\Torreg({}_AR)=\exreg(A_A)$ by Lemma \ref{dual-ext-ex} (1) and its right version.
   It follows from the right module version of (1) that 
  $$0=\CMreg(R_A)\leqslant \CMreg(A)+\Torreg(R_A)=\CMreg(A)+\exreg({}_AA),$$
  and $0 \leqslant \CMreg(A)+\exreg(A_A)$.
\end{proof}

Since $\CMreg(R) =0$, it follows from (2) that $\exreg(A_A)=\Torreg({}_{A}R)\leqslant \Torreg(S)$ and $\exreg({}_AA)=\Torreg(R_A)\leqslant \Torreg(S)$. So, (4) is stronger than (3) in Theorem \ref{relation between CM and Tor}.

By using Theorem \ref{relation between CM and Tor}, we show in the following proposition that sufficiently high truncations of $ M \in \gr A$ have linear projective resolutions if $A_0$ is semisimple and $\Torreg(S) < \infty$. In fact, this fact was proved in \cite{EG} for polynomial algebras $A$,  in \cite{AE} for  Koszul commutative graded algebras $A$, in \cite[Theorem 3.1]{Jo4} for Koszul connected graded algebras $A$ and in \cite[Theorem 3.13]{KWZ2} for connected graded algebras $A$. 

\begin{proposition}\label{linear resolution}
Suppose that $A$ is a noetherian $\mathbb{N}$-graded algebra with a balanced dualizing complex $R$. If $\Torreg(S)\leqslant p$, then for any $0\neq M \in \gr A$ with $\CMreg(M) \leqslant r$,
$\Torreg(M_{\geqslant r}(r+p))\leqslant 0$.

If, further, $A_0$ is semisimple, then $\Torreg(M_{\geqslant r}(r+p)) = 0$, and $M_{\geqslant r}(r+p)$ has a linear projective resolution.
\end{proposition}
\begin{proof}
    Since $\Torreg(M_{\geqslant r}(r+p))=\Torreg(M_{\geqslant r})-r-p$, it remains to show that $\Torreg(M_{\geqslant r})\leq r+p$. It follows from Theorem \ref{relation between CM and Tor} (2) that  
    $$\Torreg(M_{\geqslant r})\leqslant \CMreg(M_{\geqslant r})+\Torreg(S)\leqslant \CMreg(M_{\geqslant r})+p.$$  
   It suffices to prove that
   $$\CMreg(M_{\geqslant r})\leq r.$$
    
It follows from the short exact sequence
    \begin{align*}
        0 \to M_{\geqslant r} \to M \to M/M_{\geqslant r} \to 0
    \end{align*}
 and  Lemma \ref{exact sequence} that $\CMreg(M_{\geqslant r})\leqslant \max \{ \CMreg(M), \CMreg(M/M_{\geqslant r})+1 \}.$
    
     Since $M/M_{\geqslant r}$ is a torsion $A$-module, 
    $$ \CMreg(M/M_{\geqslant r})=\sup \{ j \mid  (M/M_{\geqslant r})_j \neq 0 \} \leqslant r-1.$$
     Thus $\CMreg(M/M_{\geqslant r})+1 \leqslant r$, and so
      $\CMreg(M_{\geqslant r}) \leqslant r.$     
\end{proof}

\begin{proposition} \label{relation between CM and Tor-p} Suppose $A$ is a noetherian $\mathbb{N}$-graded algebra with $A_0$ semisimple and $A$ has a balanced dualizing complex $R$. Then, for any $0\neq X \in \mathbf{D}^b(\gr A)$ with $\pdim X< \infty$,
\begin{align} \label{stronger-inequality}
\Torreg(X) \leqslant \CMreg(X)+\exreg({}_AA).
\end{align}
\end{proposition}

\begin{proof} Since $\pdim X < \infty$, $\Torreg(X)=-\ideg(R\gHom_A(X,A))$ by Theorem \ref{Extreg}. It follows from $R\gHom_A(X,A)\cong R\gHom_{A^o}(R,\D(R\Gamma_A(X)))$ that
    \begin{align*}
        \Torreg(X) 
        =&-\ideg(R\gHom_{A^o}(R,\D(R\Gamma_A(X))))\\
        \leqslant& \Torreg(R_A)-\ideg(\D(R\Gamma_A(X))) \quad (\textrm{by Proposition } \ref{sdeg-ideg inequality 1})\\
        =&\exreg({}_AA)+\CMreg(X) \quad (\textrm{by Lemma } \ref{dual-ext-ex} \,(1)) \textrm{ and  } \eqref{CMreg-identity}).
    \end{align*}
    Hence $\Torreg(X) \leqslant \CMreg(X)+\exreg({}_AA).$
\end{proof}

\begin{remark}\label{Torreg(R)=Torreg(A_0)}
    The inequality \eqref{stronger-inequality} is stronger than 
    the inequality in Theorem \ref{relation between CM and Tor} (2).
    If, furthermore, $A_0$ is semisimple and $\injdim {}_{A}A =\injdim A_A < \infty$, then $\pdim R_A < \infty$ and $\pdim {}_AR < \infty$. 
    It follows from the proof of  Proposition \ref{relation between CM and Tor-p} and its right version that
    $\Torreg(R_A)=\Torreg({}_{A}R)$, 
    and  $\exreg({}_AA)=\exreg(A_A)$.
    
    If $A_0$ is semisimple and $\gldim A < \infty$, then, by Theorem \ref{relation between CM and Tor} (2) and Proposition \ref{relation between CM and Tor-p}, $\Torreg(R)=\exreg(A)=\Torreg(A_0)$.
If  $\gldim A = \infty$, then $\Torreg(R)=\Torreg(A_0)$ may not be true. For example, if $A=k[x]/(x^2)$, then $A$ is a Koszul AS-Gorenstein algebra of dimension 0 with Gorenstein parameter $-1$, where $\Torreg(k)=0$ and $\Torreg(R)=\exreg(A)=-1$.
\end{remark}

\noindent
{\bf Question}: When is $\Torreg({}_{A}R)=\Torreg(R_A)$, or  $\exreg(A_A)=\exreg({}_AA)$?

In the connected graded case, the fact that the inequality in (1) or (2) in Theorem \ref{relation between CM and Tor} is always an equality is closed to the Artin-Schelter regular property of the algebra (\cite[Theorem 4.1]{R}, \cite[Theorem 5.4]{DW} and \cite[Corollary 3.4]{KWZ1}).
The fact that $\ASreg(A)=0$ characterizes when the inequality in (1) or (2) in Theorem \ref{relation between CM and Tor} is always an equality, as given in the following corollary.

\begin{corollary}\label{ASreg=0}
    Suppose $A$ is a noetherian $\mathbb{N}$-graded algebra with a balanced dualizing complex. Then the following are equivalent.
    \begin{itemize}
        \item [(1)] $\ASreg(A)=0$.
        
        \item [(2)] $\CMreg(X)=\Torreg(X)+\CMreg(A)$ for all $0\neq X \in \mathbf{D}^b(\gr A)$.
        \item [(3)] $\Torreg(X)=\CMreg(X)+\Torreg(S)$ for all $0\neq X \in \mathbf{D}^b(\gr A)$.
        \item [(4)] There exists some $0\neq X \in \mathbf{D}^b(\gr A)$ such that $$\CMreg(X)=\Torreg(X)+\CMreg(A) \, \textrm { and } \, \Torreg(X)=\CMreg(X)+\Torreg(S).$$
    \end{itemize}
\end{corollary}
\begin{proof}
    (1) $\Rightarrow$ (2) and (3) For any $0\neq X \in \mathbf{D}^b(\gr A)$, By Theorem \ref{relation between CM and Tor} (1) and (2), 
    \begin{align*}
    \CMreg(X)&\leqslant \Torreg(X)+\CMreg(A) \\
    &\leqslant \CMreg(X)+\Torreg(S)+\CMreg(A)\\
    &=\CMreg(X)+ \ASreg(A).
    \end{align*} 
   It follows from $\ASreg(A)=0$ and $\CMreg(X)\neq \pm \infty$ that
    $$\CMreg(X)=\Torreg(X)+\CMreg(A)\, \textrm{ and } \, \Torreg(X)=\CMreg(X)+\Torreg(S).$$

    (2) $\Rightarrow$ (1) It follows from $\CMreg(S)=0$ by taking $X=S$. 

    (3) $\Rightarrow$ (1) It follows from $\Torreg(A)=0$ by taking $X=A$. 

    (1) $\Rightarrow$ (4) If $X=S$, then $$0=\ASreg(A)=\Torreg(S)+\CMreg(A)=\CMreg(S)=0$$ and $$\Torreg(S)=0+\Torreg(S)=\CMreg(S)+\Torreg(S).$$

    (4) $\Rightarrow$ (1) Since $\CMreg(X)<\infty$ by Lemma \ref{CM-reg-finite}, it follows from  
    $$\CMreg(X)=\Torreg(X)+\CMreg(A)=\CMreg(X)+\Torreg(S)+\CMreg(A)$$ 
     that $\ASreg(A)=\Torreg(S)+\CMreg(A)=0$.
\end{proof}

\begin{corollary}\label{CMreg+Torreg=0}
    Suppose $A$ is a noetherian $\mathbb{N}$-graded algebra with $A_0$ semisimple and $A$ has a balanced dualizing complex. Then the following are equivalent.
    \begin{itemize}
        \item [(1)] $\CMreg(A)+\exreg({}_AA)=0$.
         \item [(2)] $\Torreg(X)=\CMreg(X)+\exreg({}_AA)$ for all $0\neq X \in \mathbf{D}^b(\gr A)$ with $\pdim X < \infty$.
        \item [(3)] There exists some $0\neq X \in \mathbf{D}^b(\gr A)$ with $\pdim X < \infty$ such that $$\CMreg(X)=\Torreg(X)+\CMreg(A) \, \textrm { and } \, \Torreg(X)=\CMreg(X)+\exreg({}_AA).$$
    \end{itemize}
\end{corollary}
\begin{proof}
   (1) $\Rightarrow$ (2)  For any $0\neq X \in \mathbf{D}^b(\gr A)$ with $\pdim X < \infty$, by  Theorem \ref{relation between CM and Tor} (1) and Proposition \ref{relation between CM and Tor-p}, 
    \begin{align*}
    \CMreg(X)&\leqslant \Torreg(X)+\CMreg(A) \\
    &\leqslant \CMreg(X)+\exreg({}_AA)+\CMreg(A).
    \end{align*} 
    It follows from  $\exreg({}_AA)+\CMreg(A)=0$ that
    $$\CMreg(X)=\Torreg(X)+\CMreg(A)\, \textrm{ and } \, \Torreg(X)=\CMreg(X)+\exreg({}_AA).$$

    (2) $\Rightarrow$ (1) It follows by taking $X=A$.  
 
   (1) $\Rightarrow$ (3) Take $X=A$.

   (3) $\Rightarrow$ (1) Since $\CMreg(X)<\infty$ by Lemma \ref{CM-reg-finite}, it follows from  
    $$\CMreg(X)=\Torreg(X)+\CMreg(A)=\CMreg(X)+\exreg({}_AA)+\CMreg(A)$$ 
    that $\exreg({}_AA)+\CMreg(A)=0$.
\end{proof}

By the proof (1) $\Rightarrow$ (2) in Corollary \ref{CMreg+Torreg=0}, (1) implies $(2')$: for any $0\neq X \in \mathbf{D}^b(\gr A)$ with $\pdim X < \infty$, 
 $\CMreg(X)=\Torreg(X)+\CMreg(A).$   

If, further, $\injdim {}_{A}A=\injdim A_A < \infty$, then $(2') \Rightarrow$ (1). In fact, by Remark \ref{Torreg(R)=Torreg(A_0)}), $\exreg({}_AA)=\exreg(A_A)$ and $\Torreg(R_A)=\Torreg({}_{A}R)$.  
 If we take $X={}_{A}R$ in $(2')$, then 
     $$0=\CMreg({}_{A}R)=\Torreg({}_{A}R)+\CMreg(A)=\exreg(A_A)+\CMreg(A).$$

For the dual version of $(2') \Leftrightarrow$ (1), see Corollary \ref{CMreg+Torreg=0 4}.

\begin{lemma}\label{CM+Torreg=0}
    Let $A$ be an $\mathbb{N}$-graded AS-Gorenstein algebra of dimension $d$ with Gorenstein parameters $\{\ell_1,\ell_2,\cdots,\ell_n\}$. Then the following are equivalent.
    \begin{itemize}
        \item [(1)]  $\Exreg(A)+\exreg(A)=0$.
        \item [(2)] $\ell_1=\ell_2=\cdots=\ell_n$.
     \end{itemize}
\end{lemma}
 
\begin{proof} Since $A$ is an $\mathbb{N}$-graded AS-Gorenstein algebra of dimension $d$ with Gorenstein parameters $\{\ell_1,\ell_2,\cdots,\ell_n\}$, $R\gHom_A(S,A)\cong \mathop{\bigoplus}\limits_{i=1}^n (e_{\sigma(i)}S(\ell_i))^{r_i} [-d]$. 
   Thus  
   $$\exreg(A)=-d+\max\{\ell_i\mid 1\leqslant i\leqslant n\}, \textrm{ and } \Exreg(A)=d-\min\{\ell_i\mid 1\leqslant i\leqslant n\}.$$
   It follows that $\Exreg(A)+\exreg(A)=0$ if and only if $\ell_1=\ell_2=\cdots=\ell_n$,
   that is, (1) $\Leftrightarrow$ (2).
\end{proof}
Note that if $A$ is an {\it AS-Gorenstein algebra over $A_0$} of dimension $d$ in the sense of Minamoto and Mori \cite[Definition 3.1]{MM} (see Definition \ref{def-AS-Goren over A_0}), then $\Exreg(A)+\exreg(A)=0$.

   If, furthermore, $A_0$ is semisimple, then  $\CMreg(A)=\Exreg(A)$ by Theorem \ref{ex}. Hence $\Exreg(A)+\exreg(A)=0$ if and only if $\CMreg(A)+\exreg(A)=0$.\\

 \noindent
 {\bf Question}: Suppose $A$ is a noetherian $\mathbb{N}$-graded algebra with $A_0$ semisimple and $A$ has a balanced dualizing complex. If $\CMreg(A)+\exreg({}_AA)=0$ and $\exreg({}_AA)=\exreg(A_A)$,  is $A$ an $\mathbb{N}$-graded AS-Gorenstein algebra 
 with the  Gorenstein parameters being the same?
 
\begin{proposition}\label{Extreg+CMreg}
    Suppose $A$ is a noetherian $\mathbb{N}$-graded algebra with a balanced dualizing complex. Then 
    \begin{itemize}
        \item [(1)] $-\ideg(X)\leqslant \exreg(X)+\CMreg(A)$ for all $0\neq X \in \mathbf{D}^+_{fg}(\Gr A)$.
        \item [(2)] $\exreg(X) \leqslant -\ideg(X)+\Extreg(S)$ for all $0\neq X \in \mathbf{D}^+_{fg}(\Gr A)$.
        \item [(3)] If $A_0$ is semisimple, then $\exreg(X) \leqslant -\ideg(X)+\exreg(A_A)$ for all $0\neq X \in \mathbf{D}^b(\gr A)$ with $\injdim X < \infty$.
    \end{itemize}
\end{proposition}

\begin{proof}
 (1) Since $X \cong \D(R\Gamma_{A^o}(\D(R\Gamma_A(X))))$, 
 \begin{align*}
     -\ideg(X)=&\CMreg(\D(R\Gamma_A(X)))\\
      \leqslant& \Torreg(\D(R\Gamma_A(X))) + \CMreg(A) \quad (\textrm{by Theorem } \ref{relation between CM and Tor} (1))\\
     =& \exreg(X) + \CMreg(A)\quad (\textrm{by Lemma \ref{dual-ext-ex} (1)}).
 \end{align*}
 
(2) By Lemma \ref{dual-ext-ex} (1) and Theorem \ref{relation between CM and Tor} (2), 
 \begin{align*}
     \exreg(X)=&\Torreg(\D(R\Gamma_A(X)))\\
     \leqslant& \CMreg(\D(R\Gamma_A(X)))+\Extreg(S)\\
     =&-\ideg(X)+\Extreg(S).
 \end{align*}
  
 (3) For any $0\neq X \in \mathbf{D}^b(\gr A)$ with $\injdim X < \infty$, 
 $\pdim(\D(R\Gamma_A(X)))< \infty$.
By the right version of Proposition \ref{relation between CM and Tor-p}, 
 \begin{align*}
\exreg(X)=&\Extreg(\D(R\Gamma_A(X)))\\
\leqslant &\CMreg(\D(R\Gamma_A(X))) + \exreg(A_A)\\
=&-\ideg(X) + \exreg(A_A).
 \end{align*}
\end{proof}

In fact, Proposition \ref{Extreg+CMreg} (1) and (2) are directly from the duality between $\mathbf{D}_{fg}(\Gr A)$ and $\mathbf{D}_{fg}(\Gr A^o)$ and the right version of Theorem \ref{relation between CM and Tor} (1) and (2). Proposition \ref{Extreg+CMreg} (3) is the dual version of Proposition \ref{relation between CM and Tor-p}. Similarly, Corollaries \ref{ASreg=0 2} and \ref{CMreg+Torreg=0 3} in the following are nothing new but dual versions of Corollaries \ref{ASreg=0} and \ref{CMreg+Torreg=0}.

\begin{corollary}\label{ASreg=0 2}
    Suppose $A$ is a noetherian $\mathbb{N}$-graded algebra with a balanced dualizing complex $R$. Then the following are equivalent.
    \begin{itemize}
        \item [(1)] $\ASreg(A)= 0$.
        \item [(2)] $\exreg(X) = -\ideg(X)+\Extreg(S)$ for all $0\neq X \in \mathbf{D}^b(\gr A)$.
        \item [(3)] $-\ideg(X)= \exreg(X)+\CMreg(A)$ for all $0\neq X \in \mathbf{D}^b(\gr A)$. 
        \item [(4)] There exists $0\neq X \in \mathbf{D}^b(\gr A)$ such that 
         $$-\ideg(X)= \exreg(X)+\CMreg(A) \, \textrm { and } \, \\
         \exreg(X) = -\ideg(X)+\Extreg(S).$$
    \end{itemize}
\end{corollary}
\begin{proof}
    (1) $\Rightarrow$ (2) and (3) For any $0\neq X \in \mathbf{D}^b(\gr A)$, by  Proposition \ref{Extreg+CMreg} (1) and (2), 
    \begin{align*}
    -\ideg(X)&\leqslant  \exreg(X) +\CMreg(A) \\
    &\leqslant -\ideg(X)+\Extreg(S)+\CMreg(A)\\
    &=-\ideg(X)+ \ASreg(A).
    \end{align*} 
   It follows from $\ASreg(A)=0$ that
    $$-\ideg(X)=  \exreg(X)+\CMreg(A),\, \textrm{ and } \,$$ 
    $$ \exreg(X)=-\ideg(X)+\Extreg(S).$$

    (2) $\Rightarrow$ (1) Note that
    $$\exreg(R) = -\ideg(R\gHom_A(S,R))=-\ideg(R\gHom_{A^o}(A,S))=0.$$
By taking $X=R$ in (2), then $$0=-\ideg(R)+\Extreg(S)=\CMreg(A)+\Extreg(S)=\ASreg(A).$$

    (3) $\Rightarrow$ (1) By taking $X=S$ in (3), then 
    $$0=-\ideg(S)=\exreg(S)+\CMreg(A)=\Extreg(S)+\CMreg(A)=\ASreg(A).$$

    (1) $\Rightarrow$ (4) Take $X=S$.

    (4) $\Rightarrow$ (1) Since $-\ideg(X)<\infty$ and  
    \begin{align*}
        -\ideg(X)=&\exreg(X) +\CMreg(A)\\
        =&-\ideg(X)+\Extreg(S)+\CMreg(A),
    \end{align*}
    It follows that $\ASreg(A)=\Extreg(S)+\CMreg(A)=0$.
\end{proof}

\begin{corollary}\label{CMreg+Torreg=0 3}
    Suppose $A$ is a noetherian $\mathbb{N}$-graded algebra with $A_0$ semisimple and $A$ has a balanced dualizing complex $R$. Then the following are equivalent.
    \begin{itemize}
        \item [(1)] $\CMreg(A)+\exreg(A_A)=0$.
        \item [(2)] $\exreg(X)=-\ideg(X)+\exreg(A_A)$ for all $0\neq X \in \mathbf{D}^b(\gr A)$ with $\injdim X < \infty$.
        \item [(3)] There exists $0\neq X \in \mathbf{D}^b(\gr A)$ with $\injdim X < \infty$ such that $$\exreg(X)=-\ideg(X)+\exreg(A_A) \, \textrm { and } \, 
        -\ideg(X)=\exreg(X)+\CMreg(A).$$
    \end{itemize}
\end{corollary}
\begin{proof}
   (1) $\Rightarrow$ (2)  For any $0\neq X \in \mathbf{D}^b(\gr A)$ with $\injdim X < \infty$, by Proposition \ref{Extreg+CMreg} (1) and (3), 
    \begin{align*}
    -\ideg(X)&\leqslant \exreg(X)+\CMreg(A) \\
    &\leqslant -\ideg(X)+\exreg(A_A)+\CMreg(A).
    \end{align*} 
    It follows from  $\exreg(A_A)+\CMreg(A)=0$ that
    $$-\ideg(X)=\exreg(X)+\CMreg(A), \, \exreg(X)=-\ideg(X)+\exreg(A_A).$$

    (2) $\Rightarrow$ (1) By taking $X=R$, then $$\exreg(R)=-\ideg(R)+\exreg(A_A)=\CMreg(A)+\exreg(A_A).$$ The conclusion follows form 
    $$\exreg(R)=-\ideg(R\gHom_A(S, R))=-\ideg(R\gHom_{A^o}(A, S))=0.$$ 

    (2) $\Rightarrow$ (3) By taking $X={}_{A}R$ in (2), then 
    $$ (0=) \, \exreg(R)=-\ideg(R)+\exreg(A_A), \, \textrm { and } \,$$
    $$-\ideg(R)=\CMreg(A)=0+\CMreg(A)=\exreg(R)+\CMreg(A).$$

   (3) $\Rightarrow$ (1) Since 
    $$-\ideg(X)=\exreg(X)+\CMreg(A)=-\ideg(X)+\exreg(A_A)+\CMreg(A),$$ 
    it follows from $-\ideg(X)<\infty$ that $\exreg(A_A)+\CMreg(A)=0$.
\end{proof}

\begin{corollary}\label{CMreg+Torreg=0 4}
    Suppose $A$ is a noetherian $\mathbb{N}$-graded algebra with $A_0$ semisimple and $A$ has a balanced dualizing complex. If $\injdim {}_{A}A=\injdim A_A < \infty$, then the following are equivalent.
    \begin{itemize}
        \item [(1)] $\CMreg(A)+\exreg(A)=0$.
        \item [(2)] $-\ideg(X)=\exreg(X)+\CMreg(A)$ for all $0\neq X \in \mathbf{D}^b(\gr A)$ with $\injdim X < \infty$.
    \end{itemize}
\end{corollary}
 \begin{proof}
     (1) $\Rightarrow$ (2) It is the same as (1) $\Rightarrow$ (2) in Corollary \ref{CMreg+Torreg=0 3}.

     (2) $\Rightarrow$ (1) By taking $X={}_{A}A$ in (2), then 
    $$ 0=-\ideg({}_{A}A) = \exreg({}_{A}A)+\CMreg(A) = \exreg(A)+\CMreg(A).$$
 \end{proof}

\subsection{Little numerical AS-regularity} In this subsection,
similar to the numerical AS-regularity given by \cite{KWZ1}, we define a lowercase character named regularity $\asreg(A)$, called little numerical AS-regularity of $A$. We give some relations between the little AS-regularity and other lowercase character named regularities. If $A$ has a balanced dualizing complex with $A_0$ semisimple, then we prove that $\asreg(A)=0$ if and only if $A$ is finite-dimensional.
\begin{definition}
    The  little numerical AS-regularity of a noetherian $\mathbb{N}$-graded algebra $A$ is defined to be 
    $$\asreg(A):=\cmreg(A) + \torreg(S).$$
\end{definition}
Obviously, $\asreg(A)=\cmreg(A)+\extreg(S)=\cmreg(A)+\Exreg(S)$.
\begin{proposition}\label{extreg+sdeg}
    Suppose $A$ is a noetherian $\mathbb{N}$-graded algebra with a balanced dualizing complex $R$. Then 
    \begin{itemize}
        \item [(1)] $\extreg(X)\leqslant \cmreg(X)+\extreg(S)$ for all $0\neq X \in \mathbf{D}^-(\gr A)$.
        \item [(2)] $\cmreg(X) \leqslant \extreg(X)+\cmreg(A)$ for all $0\neq X \in \mathbf{D}^-(\gr A)$.
        \item [(3)] $\asreg(A) \geqslant 0$.
    \end{itemize}
    \end{proposition}
\begin{proof}
    (1) Since $R\gHom_A(X,S)\cong R\gHom_{A^o}(S,\D(R\Gamma_A(X)))$, it follows that
    \begin{align*}
        \extreg(X)
        =&\sdeg(R\gHom_{A^o}(S,\D(R\Gamma_A(X))))\\
        \leqslant& \extreg(S)-\ideg(R\Gamma_A(X)) \quad (\textrm {by Corollary } \ref{sdeg-ideg inequality coro}\, (7))\\
        =& \extreg(S)+\cmreg(X).
    \end{align*}

    (2) Since $\D(R\Gamma_A(X))\cong R\gHom_A(X,R)$,
    \begin{align*}
        \cmreg(X)=&-\ideg(R\Gamma_A(X)) =\sdeg(\D(R\Gamma_A(X)))\\
        =&\sdeg(R\gHom_A(X,R)) \\
        \leqslant& \extreg(X)-\ideg(R\Gamma_A(A))\quad (\textrm {by Corollary } \ref{sdeg-ideg inequality coro}\, (7))\\
        =& \extreg(X)+\cmreg(A).
    \end{align*}

    (3) It follows by taking $X=S$ in (2) (or $X=A$ in (1)) that
    $$0=\cmreg(S) \leqslant  \extreg(S) + \cmreg(A) =\asreg(A).$$
\end{proof}

Note that $\extreg(S)\geqslant 0$. If $A$ has a balanced dualizing complex, then $\cmreg(A)\neq -\infty$. Then $\asreg(A)=0$ implies that both  $\extreg(S)< \infty$ and $\cmreg(A)< \infty$. Furthermore, if $A_0$ is semisimple, then $\extreg(S)=0$. Hence $\asreg(A)=0$ if and only if $\cmreg(A)=0$. 
\begin{corollary}\label{finite-dim alg}
   Suppose $A$ is a noetherian $\mathbb{N}$-graded algebra with a balanced dualizing complex $R$. Then $A$ is finite-dimensional if and only if $\cmreg(X)< \infty$ for all $0\neq X \in \mathbf{D}^b(\gr A)$.
\end{corollary}
\begin{proof}
    If $\cmreg(X)< \infty$ for all $0 \neq X \in \mathbf{D}^b(\gr A)$, then 
 $\sdeg \D(R\Gamma_A(X)) < \infty$ for all $0 \neq X \in \mathbf{D}^b(\gr A)$. This is equivalent to that
$\sdeg (Y) < \infty$ for all $0 \neq Y \in \mathbf{D}^b(\gr A^o)$, which is equivalent to that $A$ is finite-dimensional.
\end{proof}
 
\begin{corollary}\label{cmreg-finite}
 Suppose $A$ is a noetherian $\mathbb{N}$-graded algebra with a balanced dualizing complex $R$. If $A_0$ is semisimple, then the following are equivalent. 
\begin{itemize}
        \item [(1)] $\cmreg(A)<\infty$.
        \item [(2)] $\cmreg(X)< \infty$ for all $0\neq X \in \mathbf{D}^b(\gr A)$.
        \item [(3)] $A$ is finite-dimensional.
        \item [(4)] $\cmreg(A)=0$.
        \item [(5)] $\asreg(A)=0$.
        \item [(6)] $\asreg(A)< \infty$.
\end{itemize}     
\end{corollary}

\begin{proof}  (1) $\Leftrightarrow$ (2) Since $A_0$ is semisimple, it follows that 
 $\extreg(S)=0$. 
 Then  $\cmreg(A)\geqslant 0$ by Proposition \ref{extreg+sdeg} (3). Hence,  by Proposition \ref{extreg+sdeg} (2) and Proposition \ref{CM and Ext}, for any $0\neq X \in \mathbf{D}^b(\gr A)$,
$$\cmreg(X)\leqslant\extreg(X)+\cmreg(A)=-\ideg(X)+\cmreg(A).$$
 Therefore $\cmreg(A)<\infty$ if and only if $\cmreg(X)< \infty$ for all $0 \neq X \in \mathbf{D}^b(\gr A)$. 

 (2) $\Leftrightarrow$ (3) By Corollary \ref{finite-dim alg}.

(3) $\Rightarrow$ (4), (4) $\Rightarrow$ (1) and (4) $\Leftrightarrow$ (5) $\Leftrightarrow$ (6) are obvious.
\end{proof}
If $A_0$ is not semisimple, then (3) $\Rightarrow$ (5) may not be true. For example, if $A=A_0$ with $\pdim{S}=n$, then $\extreg(S)=n$, $\cmreg(A)=0$, and so $\asreg(A)=n$. This means $\asreg(A)$ may run over all positive integers. If $A_0$ is semisimple and $A$ has a balanced dualizing complex, then it follows from Corollary \ref{cmreg-finite} that $\asreg(A)=0$ or $\asreg(A)=+\infty$.

\noindent
{\bf Question}:  Suppose $A$ is a noetherian $\mathbb{N}$-graded algebra with a balanced dualizing complex. If $\asreg(A)=0$ and $A_0$ is not semisimple, is $A$ finite-dimensional?

\begin{corollary}\label{extreg+sdeg=0}
    Suppose $A$ is a noetherian $\mathbb{N}$-graded algebra with a balanced dualizing complex $R$. Then the following are equivalent.
    \begin{itemize}
        \item [(1)] $\asreg(A)=0$.
        \item [(2)] $\extreg(X)= \cmreg(X) + \extreg(S)$ for all $0\neq X \in \mathbf{D}^-(\gr A)$.
        \item [(3)] $\cmreg(X) = \extreg(X)+\cmreg(A)$ for all $0\neq X \in \mathbf{D}^-(\gr A)$.
    \end{itemize}
\end{corollary}
\begin{proof}
    (1) $\Rightarrow$ (2) By Proposition \ref{extreg+sdeg},
    \begin{align*}
        \extreg(X)\leqslant& \cmreg(X) +  \extreg(S)\\
        \leqslant& \extreg(X) + \cmreg(A) + \extreg(S)\\
        =&\extreg(X).
    \end{align*}
Hence  $\extreg(X) = \cmreg(X) +  \extreg(S).$

 (1) $\Rightarrow$ (3) Similarly, it follows from $$\cmreg(X) \leqslant \extreg(X) + \cmreg(A) \leqslant \cmreg(X) + \cmreg(A) + \extreg(S) \leqslant \cmreg(X).$$ 

    (2) $\Rightarrow$ (1) It follows by taking $X=A$.

    (3) $\Rightarrow$ (1) It follows by taking $X=S$.
\end{proof}

The following Proposition \ref{exreg+sdeg} and Corollary \ref{exreg+sdeg=0} are in fact the dual versions of Proposition \ref{extreg+sdeg} and Corollary \ref{extreg+sdeg=0}.

\begin{proposition}\label{exreg+sdeg}
    Suppose $A$ is a noetherian $\mathbb{N}$-graded algebra with a balanced dualizing complex $R$. Then 
    \begin{itemize}
        \item [(1)] $\Exreg(X)\leqslant \sdeg(X) + \Exreg(S)$ for all $0\neq X \in \mathbf{D}^+_{lf}(\Gr A)$.
        \item [(2)] $\sdeg(X)\leqslant \Exreg(X)+\cmreg(A)$ for all $0\neq X \in \mathbf{D}^+_{fg}(\Gr A)$.
    \end{itemize}
\end{proposition}
\begin{proof}
    (1) It follows from Corollary \ref{sdeg-ideg inequality coro} (7), 
   \begin{align*}
   \Exreg(X)=&\sdeg(R\gHom_A(S, X)) \\
   \leqslant &\extreg(S)-\ideg(\D(X))=\Exreg(S)+\sdeg(X).
\end{align*}

    (2) Since $X \cong \D(R\Gamma_{A^o}(\D(R\Gamma_A(X))))$ and $\D(R\Gamma_A(X)) \in \mathbf{D}^-_{fg}(\Gr A)$,
    \begin{align*}
        \sdeg(X)=&-\ideg(R\Gamma_{A^o}(\D(R\Gamma_A(X))))\\
        =&\cmreg(\D(R\Gamma_A(X))) \\
    \leqslant& \extreg(\D(R\Gamma_A(X)))+\cmreg(A) \quad (\textrm{by Proposition \ref{extreg+sdeg}} \, (2)) \\
    =& \Exreg(X)+\cmreg(A).
    \end{align*}
\end{proof}

If $A_0$ is semisimple, then $\Exreg(S)=0$ and $\Exreg(X)\leqslant \sdeg(X)$ for all $0\neq X \in \mathbf{D}^+_{lf}(\Gr A)$. Furthermore, if  $X \in \mathbf{D}^+_{lf}(\Gr A)$ with $\D(X) \in \mathbf{D}^-(\gr {A^o})$, then by Proposition \ref{CM and Ext} and $R\gHom_A(A_0,X)\cong R\gHom_{A^o}(\D(X),A_0)$, 
$$\sdeg(X)=-\ideg(\D(X))=\extreg(\D(X))=\Exreg(X).$$

\begin{corollary}\label{exreg+sdeg=0}
    Suppose $A$ is a noetherian $\mathbb{N}$-graded algebra with a balanced dualizing complex $R$. Then the following are equivalent.
    \begin{itemize}
        \item [(1)] $\asreg(A)=0$.
        \item [(2)] $\Exreg(X)=\sdeg(X) + \Exreg(S)$ for all $0\neq X \in \mathbf{D}^+_{fg}(\Gr A)$.
        \item [(3)] $\sdeg(X)= \Exreg(X)+\cmreg(A)$ for all $0\neq X \in \mathbf{D}^+_{fg}(\Gr A)$.
    \end{itemize}
\end{corollary}

\begin{proof}
   (1) $\Rightarrow$ (2)  By Proposition \ref{exreg+sdeg},
    \begin{align*}
        \Exreg(X)\leqslant& \sdeg(X) + \Exreg(S)\\
        \leqslant& \Exreg(X)+\cmreg(A)+\Exreg(S)\\
        =&\Exreg(X).
    \end{align*}
Hence $\Exreg(X)=\sdeg(X) + \Exreg(S)$.

(1) $\Rightarrow$ (3) Similarly, it follows from
 $\sdeg(X)\leqslant \Exreg(X)+\cmreg(A) \leqslant \sdeg(X) + \Exreg(S)+\cmreg(A) = \sdeg(X).$
 
    (2) $\Rightarrow$ (1) Since $\Exreg(R)=\extreg(\D(R\Gamma_A(R)))=\extreg(A)=0,$ it follows by taking $X=R$ in (2). 

    (3) $\Rightarrow$ (1) It follows by taking $X=S$ in (3). 
\end{proof}
\begin{proposition}\label{hom-ineq}
    Suppose $A$ is a noetherian $\mathbb{N}$-graded algebra with a balanced dualizing complex $R$. Then 
    \begin{itemize}
        \item [(1)] $\sdeg(R\gHom_A(X,A)) \leqslant \Exreg({}_AA)+\cmreg(X)$ for all $0\neq X \in \mathbf{D}^-(\gr A)$.
        \item [(2)] $-\ideg(R\gHom_A(X,A))\leqslant \exreg({}_AA)+\CMreg(X)$ for all $0\neq X \in \mathbf{D}^-(\gr A)$.
    \end{itemize}
\end{proposition}
\begin{proof}
   Since $R\gHom_A(X,A)\cong R\gHom_{A^o}(R,\D(R\Gamma_A(X)))$, by Corollary \ref{sdeg-ideg inequality coro} (7), 
    \begin{align*}
       \sdeg(R\gHom_A(X,A))=&\sdeg(R\gHom_{A^o}(R,\D(R\Gamma_A(X))))\\
       \leqslant& \extreg(R_A)-\ideg(R\Gamma_A(X))\\
       =& \Exreg({}_{A}A)+\cmreg(X).
    \end{align*}
Hence (1) holds. On the other hand, by Proposition \ref{sdeg-ideg inequality 1},
   \begin{align*}
       -\ideg(R\gHom_A(X,A))=&-\ideg(R\gHom_{A^o}(R,\D(R\Gamma_A(X))))\\
       \leqslant& \Extreg(R_A)-\ideg(\D(R\Gamma_A(X)))\\
       =& \exreg({}_{A}A)+\CMreg(X).
    \end{align*}
  Hence (2) holds.  
\end{proof}

\subsection{CM-regularity homogeneous and ex-regularity homogeneous} In this subsection, we introduce  CM-regularity  homogeneous  and ex-regularity homogeneous properties for noetherian $\mathbb{N}$-graded algebras. We examine the conditions under which
the inequalities stated in Theorem \ref{relation between CM and Tor} transform into equalities, and generalize \cite[Theorem 2.8]{KWZ1} to $\mathbb{N}$-graded algebras.

\begin{definition}\label{CM-homogeneous-def} Suppose $A$ is an $\mathbb{N}$-graded algebra.
 \begin{itemize}
     \item [(1)] $A$ is called {\it  left (resp. right) CM-regularity homogeneous} if $\CMreg(A)=\CMreg(Ae_i)$  (resp. $\CMreg(A)=\CMreg(e_iA)$) for all $1 \leqslant i \leqslant n$.
     \item [(2)] $A$ is called {\it left (resp. right) ex-regularity homogeneous}  if  $\exreg({}_{A}A)=\exreg(Ae_i)$ (resp. 
 $\exreg(A_A)=\exreg(e_iA)$) for all $1 \leqslant i \leqslant n$.
    \end{itemize}
    \end{definition}
\begin{remark}\label{ex-Tor homogeneous}
    If $A$ has a balanced dualizing complex, then for all $1 \leqslant i \leqslant n$, $\exreg({}_{A}A)=\exreg(Ae_i)$ (resp. $\exreg(A_A)=\exreg(e_iA)$) if and only if $\Torreg(R_A)=\Torreg(e_iR)$ (resp. $\Torreg({}_{A}R)=\Torreg(Re_i)$).
\end{remark}

Note that any connected graded algebra $A$ is CM-regularity homogeneous and ex-regularity homogeneous.  
\begin{remark}\label{left-right homogeneous}
 If $A$ is $\mathbb{N}$-graded AS-Gorenstein, then it follows from \eqref{def-GAS-Goren}, Lemma \ref{def-GAS} and \eqref{GAS-Goren} that $A$ is left CM-regularity homogeneous (resp. left ex-regularity homogeneous) if and only if $A$ is right CM-regularity homogeneous (resp. right ex-regularity homogeneous).

  Moreover, if $A$ is an $\mathbb{N}$-graded AS-Gorenstein algebra and $\exreg(A)+\Exreg(A)=0$ (i.e., the Gorenstein parameters of $A$ are the same), then $A$ is left (right) CM-regularity homogeneous and left (right) ex-regularity homogeneous.

\end{remark}

\noindent
{\bf Question}: Suppose $A$ is a noetherian $\mathbb{N}$-graded algebra. 
\begin{itemize}
    \item [(1)] If $A$ is left CM-regularity homogeneous, is $A$ necessarily right CM-regularity homogeneous?
    \item [(2)] If $A$ is left ex-regularity homogeneous, is $A$ necessarily right ex-regularity homogeneous?
\end{itemize} 

Theorem \ref{relation between CM and Tor 2} in the following gives a characterization when the inequality presented in Theorem \ref{relation between CM and Tor} (2) becomes an equality, which is an extension of \cite[Proposition 5.6]{DW} and \cite[Theorem 2.8]{KWZ1} and a direct corollary of Proposition \ref{sdeg-ideg equality 3}.

\begin{theorem}\label{relation between CM and Tor 2}
Suppose $A$ is a noetherian $\mathbb{N}$-graded algebra with $A_0$ semisimple and $A$ has a balanced dualizing complex $R$. Then the following  are equivalent.
   \begin{itemize}
       \item [(1)] $A$ is  left CM-regularity homogeneous.
       \item [(2)] $\CMreg(X)=\Torreg(X)+\CMreg(A)$ for all $0\neq X\in \mathbf{D}^b(\gr A)$ with finite projective dimension.
   \end{itemize}
\end{theorem}
\begin{proof}
By Definition \ref{CM-reg-def}, 
$$\CMreg(A) =-\ideg(\D(R\Gamma_A(A))) =-\ideg(R). $$
It follows from  Theorem \ref{LD} that, for all $1\leqslant i \leqslant n$,  
\begin{align*}
    \CMreg(Ae_i)
    =&-\ideg(\D(R\Gamma_A(Ae_i)))\\
    =&-\ideg(R\gHom_A(Ae_i,R))\\
    =&-\ideg(e_iR).
\end{align*}
Therefore, $\CMreg(A)=\CMreg(Ae_i)$ if and only if $-\ideg(R)=-\ideg(e_iR)$.

Similarly,
\begin{align*}
    \CMreg(X)
    =&-\ideg(\D(R\Gamma_A(X)))\\
    =&-\ideg(R\gHom_A(X,R)).
\end{align*} 
It follows from Proposition \ref{sdeg-ideg equality 3} that for any $X\in \mathbf{D}^b(\gr A)$ with finite projective dimension, 
$$\CMreg(X)=-\ideg(R\gHom_A(X,R))=\Extreg(X)-\ideg(R)$$
if and only if $\ideg(e_iR)=\ideg(R)$ for all $1\leqslant i\leqslant n$.
Hence 
$$\CMreg(X)=\Torreg(X)+\CMreg(A)$$ holds for all $0\neq X\in \mathbf{D}^b(\gr A)$ with finite projective dimension if and only if $A$ is  left CM-regularity homogeneous.
\end{proof}

The following is the dual version of the previous theorem. 
\begin{proposition}\label{right CM-reg homogenous}
    Suppose $A$ is a noetherian $\mathbb{N}$-graded algebra with $A_0$ semisimple and $A$ has a balanced dualizing complex $R$. Then the following are equivalent. 
    \begin{itemize}
        \item [(1)] $A$ is right CM-regularity homogeneous.
        \item [(2)] $-\ideg(X)=\exreg(X)+\CMreg(A)$ for all $X \in \mathbf{D}^b(\gr A)$ with $\injdim X < \infty$.
    \end{itemize}
\end{proposition}  
\begin{proof}
It follows from the right version of Theorem \ref{relation between CM and Tor 2} that $A$ is right CM-regularity homogeneous if and only if for all $X \in \mathbf{D}^b(\gr A)$ with $\injdim X < \infty$, $$\CMreg(\D(R\Gamma_A(X)))
    =\Torreg(\D(R\Gamma_A(X)))+\CMreg(A).$$
The conclusion follows from $$-\ideg(X)=\CMreg(\D(R\Gamma_A(X))) \textrm{ and } \exreg(X)=\Torreg(\D(R\Gamma_A(X))).$$
\end{proof}

\begin{proposition}\label{relation between CM and Tor 3}
    Suppose $A$ is a noetherian $\mathbb{N}$-graded algebra with $A_0$ semisimple and $A$ has a balanced dualizing complex $R$. If $\injdim {}_{A}A=\injdim A_A <\infty$, then the following statements are equivalent.
    \begin{itemize}
        \item [(1)] $A$ is right ex-regularity homogeneous. 
        \item [(2)] $\Extreg(X)=\CMreg(X)+\exreg(A)$ for all $0\neq X\in \mathbf{D}^b(\gr A)$ with $\pdim X < \infty$.
       \item [(3)] $-\ideg(R\gHom_A(X,A))=\CMreg(X)+\exreg(A)$ for all $0\neq X\in \mathbf{D}^b(\gr A)$.
    \end{itemize}
    \end{proposition}
    \begin{proof}
  (1) $\Rightarrow$ (3) It follows from $R\gHom_{A^o}(R,e_iA)\cong R\gHom_A(\D(R\Gamma_{A^o}(e_iA)), A)\cong R\gHom_A(Re_i,A)$ that 
        \begin{align*}
            -\ideg(R\gHom_{A^o}(R,e_iA))=&-\ideg(R\gHom_A(Re_i,A)) \\
            =&\Extreg(Re_i) \quad (\textrm{by Theorem \ref{Extreg}})\\
            =&\Torreg(Re_i) \quad (\textrm{by Lemma \ref{Extreg-Torreg}})\\
            =&\Torreg(R) \quad (\textrm{by (1) and Remark } \ref{ex-Tor homogeneous})\\
            =&\Extreg(R) \quad (\textrm{by Lemma \ref{Extreg-Torreg}}).
        \end{align*}
Since  $\pdim R_A < \infty$, 
it follows from the the right version of Proposition \ref{sdeg-ideg equality 4} that for any $0\neq X \in \mathbf{D}^b(\gr A)$, 
\begin{align*}
    -\ideg(R\gHom_A(X,A))
    =&-\ideg(R\gHom_{A^o}(R,\D(R\Gamma_A(X))))\\
    =&\Extreg(R_A)-\ideg(\D(R\Gamma_A(X)))\\
    =&\exreg(A)+\CMreg(X).
\end{align*}

  (3) $\Rightarrow$ (2) If $\pdim X < \infty$, then by Theorem \ref{Extreg}, 
  $$\Torreg(X)=-\ideg(R\gHom_A(X,A))=\CMreg(X)+\exreg(A).$$

  (2) $\Rightarrow$ (1)  It follows from $\D(R\Gamma_A(Re_i))\cong e_iA$ that  $\CMreg(Re_i)=0$ for any $1\leqslant i \leqslant n$. Hence by taking $X=Ae_i$ in (2),  
  $$\exreg(Ae_i)=\CMreg(Ae_i) + \exreg(A)=\exreg(A).$$ 
 It follows that $A$ is right ex-regularity homogeneous.
    \end{proof}
    
    The following corollary is the dual version of (1) $\Leftrightarrow$ (2) in Proposition \ref{relation between CM and Tor 3}. 
    
\begin{corollary}\label{the right Tor-regularity homogeneous}
       Suppose $A$ is a noetherian $\mathbb{N}$-graded algebra with $A_0$ semisimple and $A$ has a balanced dualizing complex $R$. If $\injdim {}_{A}A=\injdim A_A <\infty$, then the following statements are equivalent.
       \begin{itemize}
        \item [(1)] $\exreg(X)=-\ideg(X)+\exreg(A)$ for all $0\neq X\in \mathbf{D}^b(\gr A)$ with $\injdim X < \infty$.
       \item [(2)] $A$ is left ex-regularity homogeneous.
       \end{itemize}
    \end{corollary}
    \begin{proof}
    It follows from the right version of Proposition \ref{relation between CM and Tor 3} that $A$ is left ex-regularity homogeneous if and only if for all $0\neq X\in \mathbf{D}^b(\gr A)$ with $\injdim X < \infty$, 
   $$ \Extreg(\D(R\Gamma_A(X)))=\CMreg(\D(R\Gamma_A(X)))+\exreg(A).$$
The conclusion follows from  $$\exreg(X)=\Extreg(\D(R\Gamma_A(X))) \textrm{ and }
 -\ideg(X) = \CMreg(\D(R\Gamma_A(X))).$$
    \end{proof}
    \begin{corollary}\label{right ex-homo}
       Suppose $A$ is a noetherian $\mathbb{N}$-graded algebra with $A_0$ semisimple and $A$ has a balanced dualizing complex $R$. If $\gldim A < \infty$, then the following are equivalent
       \begin{itemize}
        \item [(1)] $\Torreg(X)=\CMreg(X)+\Torreg(A_0)$ for all $0\neq X\in \mathbf{D}^b(\gr A)$.
       \item [(2)] $A$ is right ex-regularity homogeneous. 
    \end{itemize}
    \end{corollary}
    \begin{proof}
        If $\gldim A < \infty$, then by Remark \ref{Torreg(R)=Torreg(A_0)}, $\Torreg(R)=\Torreg(A_0)$. Therefore by Proposition \ref{relation between CM and Tor 3}, the conclusion holds.
    \end{proof}
     The dual version of Corollary \ref{right ex-homo} is given in the following.
 \begin{corollary}\label{dual of ex-hom}
     Suppose $A$ is a noetherian $\mathbb{N}$-graded algebra with $A_0$ semisimple and $A$ has a balanced dualizing complex $R$. If $\gldim A < \infty$, then the following are equivalent.
     \begin{itemize}
        \item [(1)] $\exreg(X)=-\ideg(X)+\Extreg(A_0)$ for all $0\neq X\in \mathbf{D}^b(\gr A)$.
       \item [(2)] $A$ is left ex-regularity homogeneous.
       \end{itemize}
 \end{corollary}   
 \begin{proof}
 It follows from the right version of Corollary \ref{right ex-homo} that $A$ is left ex-regularity homogeneous if and only if for all $0\neq X\in \mathbf{D}^b(\gr A)$, 
 $$\Torreg(\D(R\Gamma_A(X)))=\CMreg(\D(R\Gamma_A(X)))+\Torreg(A_0).$$
 The conclusion follows from that 
 $$ \exreg(X)=\Torreg(\D(R\Gamma_A(X)))  \textrm{ and }
 -\ideg(X) = \CMreg(\D(R\Gamma_A(X))). $$
 \end{proof}
 
\section{Relationship between homological regularities and AS-regular property}
In this section, following the idea in\cite{KWZ1}, we study the relation between the homological regularities and the $\mathbb{N}$-graded AS-regular property for noetherian $\mathbb{N}$-graded algebras, which generalizes the results in the connected graded case.
 
A noetherian $\mathbb{N}$-graded algebra $A$ is said to satisfy the left {\it Auslander-Buchsbaum Formula}, if for any $0\neq X \in \mathbf{D}^b(\gr A)$ with $\pdim X < \infty$, $$\pdim X +\depth X =\depth A.$$ 

Note that any noetherian connected graded algebra satisfying  the $\chi$-condition  satisfies the left Auslander-Buchsbaum Formula \cite[Theorem 3.2]{Jo2}. For $\mathbb{N}$-graded algebras $A$ with a balanced dualizing complex, the following theorem characterizes when $A$ satisfies the Auslander-Buchsbaum formula.  

 \begin{theorem}\cite[Theorem 1.4]{LW}
     Let $A$ be a noetherian $\mathbb{N}$-graded algebra with a balanced dualizing complex. Then $A$ satisfies the left Auslander-Buchsbaum Formula if and only if  $\gExt_{A^o}^d(M,A)\neq 0$ for all graded simple $A^o$-module $M$, where $d=\depth{A}$.
 \end{theorem}
 
Given that the Auslander-Buchsbaum Formula is satisfied, the subsequent characterization of $\mathbb{N}$-graded AS-Gorenstein algebras extends the results presented in \cite[Theorem 3.6]{DW}.
\begin{proposition} \cite[Theorem 8.14]{LW} \label{ABF}
 Let $A$ be a noetherian $\mathbb{N}$-graded algebra with a balanced dualizing complex $R$. If $A$ satisfies the left and right Auslander-Buchsbaum Formula, then the following are equivalent.
  \begin{itemize}
      \item [(1)] $A$ is an $\mathbb{N}$-graded AS-Gorenstein algebra.
      \item [(2)] $\injdim({}_{A}A) < \infty$.
      \item [(3)] $\pdim({}_{A}R) < \infty$.
      \item [(4)] For any $X \in \mathbf{D}^b(\gr A)$, $\pdim(X) < \infty$ if and only if  $\injdim(X) < \infty$.
  \end{itemize}
\end{proposition}

\begin{lemma}\label{graded equivalence}
 Suppose $A$ is an $\mathbb{N}$-graded algebra. Let $P=\mathop{\bigoplus}\limits_{i=1}^n Ae_{i}(p_i)$ for some integers $p_1, p_2,\dots, p_n$, and $B=\gEnd_A(P)$. Then 
 \begin{itemize}
     \item [(1)] $A$ satisfies the right Auslander-Buchsbaum
Formula if and only if so is $B$. 
     \item [(2)] $A$ has a balanced dualizing complex if and only if so is $B$.
 \end{itemize}
\end{lemma}
\begin{proof}
(1) Since $P=\mathop{\bigoplus}\limits_{i=1}^n Ae_{i}(p_i)$ is a finitely generated graded projective generator, $A$ is graded Morita equivalent to $B$. For any $X\in \mathbf{D}^b(\gr 
 A^o)$, $\pdim{X} = \pdim X\otimes_A P$, and  $R\gHom_{A^o}(S,X)\cong R\gHom_{B^o}(S\otimes_A P,X\otimes_A P)$. Note that any graded simple $B^o$-modules is isomorphic to a direct summand of $S \otimes_A P$ up to a shift. It follows that  $\depth{X}=\depth{X\otimes_A P}$  and $\depth{A}=\depth{P_B}=\depth{B}$. Hence the conclusion holds.

 (2) See \cite[Proposition 3.16]{LW}.
\end{proof}

To prove the main result in this section, we modify \cite[Theorem 4.7]{IKU} in Theorem \ref{AGP} so that the algebra $B$ is graded Morita equivalent to $A$ and $B_0$ is semisimple. Note that the modules considered in \cite{IKU} are right graded modules, while the modules here are left graded modules.

 \begin{theorem}\label{AGP}
     Let $A$ be a ring-indecomposable basic $\mathbb{N}$-graded AS-Gorenstein algebra of dimension $d$ with the average Gorenstein parameter $\ell^A_{av}\in\mathbb{Z}$. If $A_0$ is semisimple,
  then there exists a ring-indecomposable basic $\mathbb{N}$-graded AS-Gorenstein algebra $B$ of dimension $d$ with all the Gorenstein parameters are $\ell^A_{av}$ such that 
  $B_0$ is semisimple and $B$ is graded Morita equivalent to $A$.
 \end{theorem}

\begin{proof}
Overall we follow the lines of the proof of \cite[Theorem 4.7]{IKU}, but with several modifications implemented. Let $I=\{1,2,\dots,n\}$ and define
\begin{align}
    m^A(i,j):=\min\{\ell\mid e_iA_{\ell}e_j\neq 0\}.
\end{align}
If there exists a ring-indecomposable basic $\mathbb{N}$-graded AS-Gorenstein algebra $B$ satisfying that $\ell^B_i=\ell^B_{av}=\ell^A_{av}$ for all $i$ and that $B$ is graded Morita equivalent to $A$ such that for any $i\neq j \in I$, $m^B(i,j)\geqslant 1$ and for any $i$, $m^B(i,i)=0$, then $B_0=\mathop{\bigoplus}\limits_{i=1}^n e_iB_0e_i$.

Since $A_0$ is semisimple, $e_iA_0e_i$ is a simple ring for each $i$. It follows from the ring isomorphisms
\begin{align*}
    e_iA_0e_i \stackrel{op}{\cong}& \Hom_{\Gr A}(Ae_i,Ae_i)\\
    \cong& \Hom_{\Gr B}(\gHom_A(P,Ae_i),\gHom_A(P,Ae_i))\\
    \cong& \Hom_{\Gr B}(Be_i,Be_i)\\
    \stackrel{op}{\cong}& e_iB_0e_i
\end{align*}
 that $B_0$ is semisimple.

Now it remains to show that there exists a ring-indecomposable basic $\mathbb{N}$-graded AS-Gorenstein algebra $B$ satisfying that $\ell^B_i=\ell^B_{av}=\ell^A_{av}$ for each $i$ and that $B$ is graded Morita equivalent to $A$ such that $m^B(i,j)\geqslant 1$ for any $i\neq j \in I$. It suffices to make the following modifications to the proof of \cite[Theorem 4.7]{IKU}. 
\begin{itemize}
    \item [(1)] Replace Definition A.1 with the following:
   Let $I$ be a finite set and $m:I^2\to \mathbb{Z}$ a map. The map $m$ is called positive if $m(i,j)\geqslant 1$ for all $i\neq j \in I$ and $m(i,i)\geq 0$ for all $i\in I$. 
  Let $\Sq(I)$ denote the set of sequences in $I$ of positive length with distinct adjacent elements, that is, for $\qqq=(i_0,\dotsc,i_{n-1})$ with $ i_k\in I$, $\qqq \in \Sq(I)$ if and only if $i_k\neq i_{k+1}$ for $k\in\mathbb{Z}/n\mathbb{Z}$. 
  Note that such a sequence must have its length $\abs{\qqq}\geqslant 2$.
  For $\qqq\in\Sq(I)$, we define
  \begin{equation*}
    m(\qqq):=\sum_{k\in\mathbb{Z}/n\mathbb{Z}} m(i_k,i_{k+1}).
  \end{equation*}
  Then $m$ is called \emph{$\Sigma$-positive} if $m(i,i)\geq 0$ 
  for all $i\in I$ and $m(\qqq)\geq \abs{\qqq}$ for all $\qqq\in\Sq(I)$. 
  \item [(2)] Replace Theorem A.2 with the following: Let $m:I^2\to \mathbb{Z}$ be a map. Then $m$ admits a positive conjugate if and only if $m$ is $\Sigma$-positive.
  \item [(3)] Modify $m_{\min}:=\min\{m(\qqq) \mid\qqq\in \Sq(I), \abs{\qqq}\geqslant 2\}$ to $m_{\min}:=\min\{m(\qqq)-\abs{\qqq} \mid \qqq \in \Sq(I)\}$ in the third line below Theorem A.2.
  \item [(4)] As for Lemma A.3, use the following statement instead:
   Let $m:I^2\to\mathbb{Z}$ be a map such that
  $m(i,i)\geqslant 0$ for $i\in I$.
  \begin{itemize}
    \item [(1$'$)]$m(\qqq)=(sm)(\qqq)$ holds for each $s:I\to\mathbb{Z}$
    and each $\qqq\in \Sq(I)$. Then, being $\Sigma$-positive is preserved under conjugation.
    \item[(2$'$)] If $m$ is $\Sigma$-positive, then $m_{\min}=m(\qqq)-\abs{\qqq}$ is valid for some multiplicity-free $\qqq\in\Sq(I)$.
    \item [(3$'$)] If $m(\qqq)\geqslant \abs{\qqq}$ holds
    for each multiplicity-free $\qqq\in\Sq(I)$,
    then $m$ is $\Sigma$-positive.
    \item [(4$'$)] $m$ is $\Sigma$-positive if and only if
    $\sum_{i\in I'}m(i,\tau i)\geqslant 
   \vert I' \vert$
    holds for each $\tau\in\Aut(I)$, where $I'=\{i\in I \mid \tau(i)\neq i\}$.
 \end{itemize}
 \item [(5)] Replace $s(i)=\sum_{j=0}^{i-1}m(j,j+1)$ for $i \in [1,n-1]$ with $s(i)=\sum_{j=0}^{i-1}m(j,j+1)-i$ for $i \in [1,n-1]$ in the first line of the proof of Lemma A.4.
\end{itemize}
\end{proof}
   
For the convenience, we make the following hypothesis.
 \begin{hypothesis}\label{Hy}
     $A$ is a noetherian (locally finite) $\mathbb{N}$-graded algebra with a balanced dualizing complex $R$, and $A$ satisfies the left and right Auslander-Buchsbaum Formula. 
 \end{hypothesis}
Now, we are ready to prove Theorems \ref{AS regular 11}, \ref{AS regular 22}, \ref{AS regular 33}, and their corollaries, generalizing \cite[Theorems 3.2, 0.8]{KWZ1}. 
For any noetherian connected graded algebra $A$ with a balanced dualizing complex, \cite[Theorem 3.2, 0.8]{KWZ1} says that $A$ is AS-regular if and only if that  $A$ is Cohen–Macaulay and $\ASreg(A) = 0$; if and only if $\ASreg(A) = 0$.

 \begin{theorem}\label{AS regular 11}
Suppose that $A$ satisfies Hypothesis \ref{Hy}, $A$ is basic, ring-indecomposable and $A_0$ is semisimple. Then the following are equivalent.
\begin{itemize}
    \item [(1)] $A$ is $\mathbb{N}$-graded AS-regular of dimension $d$ with the average Gorenstein parameter $\ell_{av}^A \in \mathbb{Z}$.
    \item [(2)] For some integers $p_1, p_2,\dots, p_n$,  $B: =\gEnd_A(\mathop{\bigoplus}\limits_{i=1}^n Ae_{i}(p_i))$ is an $\mathbb{N}$-graded CM-algebra of dimension $d$ with $B_0$ semisimple  such that $\ASreg(B)=0$.
\end{itemize}
\end{theorem}
\begin{proof}
    (1) $\Rightarrow$ (2) By Theorem \ref{AGP} and Lemma \ref{graded Morita equi} there is a ring-indecomposable basic $\mathbb{N}$-graded AS-regular algebra $B=\gEnd_A(\mathop{\bigoplus}\limits_{i=1}^n Ae_{i}(p_i))$ of dimension $d$ with Gorenstein parameters $\{\ell_{av}^A,\ell_{av}^A,...,\ell_{av}^A\}$ such that $B_0$ is semisimple. It follows from Theorem \ref{Bdc} and Definition \ref{GAS-Gorenstein} that $B$ has a balanced dualizing complex $R_B=B[d](-\ell_{av}^A)$. Hence $B$ is a CM-algebra of dimension $d$ (see Definition \ref{Defi-CM}), and $\CMreg(B)=-\ideg(R)=d-\ell_{av}^A$.

  By Theorem \ref{Extreg},  $\Extreg(X)=-\inf \{i+j \mid \gExt^i_B(X,B)_j \neq 0 \}$ for any $X \in \mathbf{D}^b(\gr B)$.
    Thus 
    \begin{align*}
        \CMreg(X)=&-\ideg \big(D(R\Gamma(X))\big)=-\ideg \big(\RHom(X, R_B)\big)\\
        =&-\inf \{ i+j \mid \gExt_B^i(X,R_B)_j \neq 0 \}\\
        =&-\inf \{ i+j \mid \gExt_B^i(X, B[d](-\ell_{av}^A))_j \neq 0 \}\\
        =&\Extreg(X)+d-\ell_{av}^A=\Torreg(X)+\CMreg(B).
    \end{align*}
    It follows from Corollary \ref{ASreg=0} that $\ASreg(B)=0$.

(2) $\Rightarrow$ (1)
Suppose first that $\gldim B =d_1< \infty$. It follows from Lemma \ref{graded equivalence} that $B$ satisfies the left and right Auslander-Buchsbaum
Formula and has a balanced dualizing complex. By Proposition \ref{ABF}, $B$ is an $\mathbb{N}$-graded AS-regular algebra of dimension $d_1$, say, with Gorenstein parameters $\{\ell_1,\ell_2,...,\ell_n\}$. 
Since $B$ satisfies the left Auslander-Buchsbaum Formula, $$d_1=\injdim{B_B}=\pdim{{}_BR}=\depth{B}-\depth{R}=\depth{B}=d.$$
Hence 
\begin{align*}
    \Exreg(B)=&d+\sdeg (\gExt^d_{B}(B_0,B))\\
    =& d+ \sdeg \big(\mathop{\bigoplus}\limits_{i=1}^n e_{\sigma (i)}B_0(\ell_i)\big) \quad (\textrm {by } \eqref{def-GAS-Goren})\\
    =& d-\min\{\ell_i\}.
\end{align*}
By Theorem \ref{ex}, $\CMreg(B)=\Exreg(B)= d-\min\{\ell_i\}.$
It follows from $\ASreg(B)=0$ that 
$\Torreg(B_0) =- \CMreg(B) = - d + \min\{\ell_i\}.$

On the other hand, let $P^{\bullet}$ be the minimal graded projective resolution of ${}_{B}B_0$. It follows from the AS-regular property of $B$ that $P^{-d}=\mathop{\bigoplus}\limits_{i=1}^n Be_{\sigma(i)}(-\ell_{i})$.
Then $\sdeg(B_0\otimes_B P^{-d})=\max\{\ell_i\}$, and 
\begin{align*}
    \Torreg(B_0) =& \sup\{-i + \sdeg \Tor^B_i(B_0, B_0)\}\\
    \geqslant &-d + \sdeg(B_0\otimes_B P^{-d})=-d + \max\{\ell_i\}.
\end{align*} 
It follows that $- d + \min\{\ell_i\} \geqslant -d + \max\{\ell_i\}$, and so, 
$\min\{\ell_i\}=\max\{\ell_i\}=\ell_{av}^B$.
By Theorem \ref{AGP}, 
 $\ell_{av}^A=\ell_{av}^B \in \mathbb{Z}$. Since $A$ is graded Morita equivalent to $B$,
 $A$ is $\mathbb{N}$-graded AS-regular.

 To finish the proof, it suffices to prove that $\gldim B < \infty$. Suppose $\gldim B = \infty$. Then $\pdim (B_0)=+\infty$. Let
 \begin{align*}
    ...\to P^{-j} \mathop{\to}^{\partial^{-j}} P^{-j+1}\to ... \to P^{-1} \to B \to B_0 \to 0
 \end{align*}
 be a minimal graded projective resolution of $B_0$ as left $B$-module.  Then,
 \begin{align*}
    ...\to P^{-j-2} \to P^{-j-1} \mathop{\to}^{\partial^{-j-1}} \ker \partial^{-j} \to 0
 \end{align*}
 is a minimal graded projective resolution  and $\pdim (\ker \partial^{-j})=+\infty$ for all $j \geqslant 0$.

 Since $B$ is a CM-algebra of dimension $d$, $\D(R\Gamma_B(B))\cong \D(R^d\Gamma_B(B))[d]:=\omega[d]$. 
 We claim that $\injdim_B \omega=d$ and $d=\cd(\Gamma_B)$. 
 
 By Theorem \ref{LD}, $\D(R\Gamma_B(X)) \cong \RHom_B(X,\omega[d])$ for all $X \in \mathbf{D}^b(\gr B)$.
 Then, for any $i> 0$ and $M \in \gr B$, $\gExt^i_B(M,\omega[d])\cong \D(R^{-i}\Gamma_B(M))=0$. Note that for any torsion module $0 \neq N$, $\gExt^0_B(N,\omega[d])\cong \D(R^0\Gamma_B(N))\cong \D(N) \neq 0$. It follows that $\injdim(\omega[d])=0$, and so $\injdim \omega=d$. 
 
 Since 
 $\D(R^q\Gamma_B(M))\cong \gExt^{-q}_B(M,\omega[d]) = \gExt^{-q+d}_B(X,\omega)=0$ for any $q> d$ and $M \in \Gr B$, $d=\cd(\Gamma_B)$.

Since $\RHom_B(B_0,\omega) =\RHom_B(B_0,\omega[d])[-d] \cong \D(R\Gamma_B(B_0))[-d]\cong B_0[-d]$,
 \begin{align*}
    \D(R^i\Gamma_B(\ker \partial^{-j} ))&\cong \gExt^{-i}_B(\ker \partial^{-j}, \omega[d])\\
    &\cong \gExt^{d-i}_B(\ker \partial^{-j},\omega)\\
    &\cong \gExt^{d-i+j+1}_B(B_0,\omega).
 \end{align*}
 Therefore, $\D(R^i\Gamma_B(\ker \partial^{-j} ))= 0$ for any $i< d$ and $j \geqslant d-1$.
It follows that for all $j \geqslant d-1$, $\CMreg(\ker \partial^{-j})= d + \sdeg(R^d\Gamma_B(\ker \partial^{-j}))$. 

If $j\geqslant d$, then the short exact sequence
    $0\to \ker \partial^{-j} \to P^{-j} \mathop{\to}^{\partial^{-j}} \ker \partial^{-j+1} \to 0$
 induces an exact sequence
\begin{align} \label{exact-seq-RGamma}
   0 \to R^d\Gamma_B(\ker \partial^{-j}) \to R^d\Gamma_B(P^{-j}) \to R^d\Gamma_B(\ker \partial^{-j+1}) \to 0.
\end{align}
So
\begin{equation}\label{CMreg-kernel}
\begin{aligned}
    \CMreg(\ker \partial^{-j})=&d+ \sdeg(R^d\Gamma_B(\ker \partial^{-j}))\\
    \leqslant &d + \sdeg( R^d\Gamma_B(P^{-j}))  \quad (\textrm{by } \eqref{exact-seq-RGamma})\\
    =&\CMreg(P^{-j})  \quad (B \textrm{ is Cohen-Macaulay})
\end{aligned}
\end{equation}
and 
\begin{equation}\label{sup-degree-B-tensor-P}
\begin{aligned}
    &\sdeg(B_0\otimes_B P^{-j-1})\\
    =&-0 + \sdeg(B_0\otimes_B P^{-j-1})
    \leqslant \Torreg(\ker \partial^{-j})\\
    \leqslant &\CMreg(\ker \partial^{-j})+\Torreg(B_0)  \quad (\textrm{by Theorem} \ref{relation between CM and Tor} \,(2))\\
    \leqslant& \CMreg(P^{-j})+\Torreg(B_0)  \quad (\textrm{by } \eqref{CMreg-kernel})\\
    =&\Torreg(P^{-j})  \quad (\textrm{by Corollary } \ref{ASreg=0})\\
    =&\sdeg(B_0\otimes_B P^{-j}).
\end{aligned}
\end{equation}
Let $c:=\Torreg(B_0)=-\CMreg(B)< \infty$, and $c_1$ be a  positive integer such that $c_1>c$. Since $\sdeg(B_0\otimes_B P^{-j})\geqslant j$ when $P^{-j}\neq 0$, it follows from \eqref{sup-degree-B-tensor-P} that
$$ d+c_1\leqslant \sdeg(B_0\otimes_B P^{-d-c_1})\leqslant \dots \leqslant \sdeg(B_0\otimes_B P^{-d}).$$
Hence, $\sdeg(B_0\otimes_B P^{-d})-d  \geqslant c_1 > c$, which contradicts with  $c=\Torreg(B_0)$.
\end{proof}

 Note that if $A$ is
$\mathbb{N}$-graded AS-regular of dimension $d$ with the average Gorenstein parameter $\ell_{av}^A \in \mathbb{Z}$, then $A$ is graded Morita equivalent to a basic
AS-regular graded algebra $B$ over $B_0$ of dimension $d$ in sense of Minamoto and Mori \cite[Definition 3.1]{MM}.

Keep the assumptions in Theorem \ref{AS regular 11}. If $A$ is $\mathbb{N}$-graded AS-regular of dimension $d$ such that the 
Gorenstein parameters $\ell_1, \ell_2,...,\ell_n$ are not the same, then obviously, $\ASreg(A)=\max\{\ell_i\}-\min\{\ell_i\}>0$.

There are more equivalent conditions of (1) in Theorem \ref{AS regular 11}.
\begin{corollary}\label{coro ASreg=0}
  Suppose that $A$ satisfies Hypothesis \ref{Hy}, $A$ is basic, ring-indecomposable and $A_0$ is semisimple. Then the following are equivalent.
   \begin{itemize}
       \item [(1)] $A$ is $\mathbb{N}$-graded AS-regular of dimension $d$ with the average Gorenstein parameter $\ell_{av}^A \in \mathbb{Z}$.
   \item [(3)] $B=\gEnd_A(\mathop{\bigoplus}\limits_{i=1}^n Ae_{i}(p_i))$ is an $\mathbb{N}$-graded CM-algebra of dimension $d$ for some integers $ p_1, p_2,..., p_n$, such that $B_0$ is semisimple and $\CMreg(X)-\Torreg(X)=c$ is a constant for all $0 \neq X \in \mathbf{D}^b(\gr B)$.
   \item [(4)] $B=\gEnd_A(\mathop{\bigoplus}\limits_{i=1}^n Ae_{i}(p_i))$ is an $\mathbb{N}$-graded CM-algebra of dimension $d$ for some integers $ p_1, p_2,..., p_n$, such that $B_0$ is semisimple and $\CMreg(M)-\Torreg(M)=c$ is a constant for all $0 \neq M \in \gr B$.
   
    \item [(5)] $B=\gEnd_A(\mathop{\bigoplus}\limits_{i=1}^n Ae_{i}(p_i))$ is an $\mathbb{N}$-graded CM-algebra of dimension $d$ for some integers $ p_1, p_2,..., p_n$, such that $B_0$ is semisimple and  $\exreg(X)+\ideg(X)=-c$ is a constant for all $0 \neq X \in \mathbf{D}^b(\gr B)$.
    
    \item [(6)] $B=\gEnd_A(\mathop{\bigoplus}\limits_{i=1}^n Ae_{i}(p_i))$ is an $\mathbb{N}$-graded CM-algebra of dimension $d$ for some integers $ p_1, p_2,..., p_n$, such that $B_0$ is semisimple and $\exreg(M)+\ideg(M)=-c$ for all $0 \neq M \in \gr B$.
   \end{itemize}
\end{corollary}
\begin{proof}
    (1) $\Rightarrow$ (3) By taking $B$ as in Theorem \ref{AS regular 11}, then $\ASreg(B)=0$ and $B_0$ is semisimple. It follows from Corollary \ref{ASreg=0} that $\CMreg(X)-\Torreg(X)=\CMreg(B)=c$ for all $0\neq X \in \mathbf{D}^b(\gr B)$, where $c=\CMreg(B)\in \mathbb{Z}$ is a constant. 

    (3) $\Rightarrow$ (4) Obviously.

    (4) $\Rightarrow$ (1) Obviously, $c=\CMreg(B)$ by taking $M=B$. Then it follows that $\ASreg(B)=0$ by taking $M=B_0$. Hence by Theorem \ref{AS regular 11}, $A$ is $\mathbb{N}$-graded AS-regular of dimension $d$ with average Gorenstein parameters $\ell_{av}^A \in \mathbb{Z}$.

    (1) $\Rightarrow$ (5) By taking $B$ as in Theorem \ref{AS regular 11}, then  $\ASreg(B)=0$ and $B_0$ is semisimple. It follows from Corollary \ref{ASreg=0 2} that $\exreg(X)+\ideg(X)=\Extreg(B_0)=-c$ for all $0\neq X \in \mathbf{D}^b(\gr B)$, where $c=\CMreg(B)$. 

    (5) $\Rightarrow$ (6) Obviously.

    (6) $\Rightarrow$ (1) By taking $M=B_0$, then $-c=\exreg(B_0)=\Extreg(B_0)$. Suppose $R=\omega[d]$ and let $M=\omega$. Then by (6), $\exreg(\omega)+\ideg(\omega)=\Extreg(B_0)$. So, 
$\Extreg(B_0)=\exreg(R)+\ideg(R)=\Extreg(B)-\CMreg(B)=-\CMreg(B),$ that is, $\ASreg(B)=0$.
    Hence by Theorem \ref{AS regular 11}, $A$ is $\mathbb{N}$-graded AS-regular of dimension $d$ with average Gorenstein parameters $\ell_{av}^A \in \mathbb{Z}$.
\end{proof}

In fact, the condition that $B$ is a CM-algebra of dimension $d$ in (2) of Theorem \ref{AS regular 11} is superfluous.
\begin{theorem}\label{AS regular 22}
 Suppose that $A$ satisfies Hypothesis \ref{Hy}, $A$ is basic, ring-indecomposable and $A_0$ is semisimple. Then the following are equivalent.
     \begin{itemize}
         \item [(1)] $A$ is $\mathbb{N}$-graded AS-regular of dimension $d$ with average Gorenstein parameter $\ell_{av}^A \in \mathbb{Z}$.
         \item[(7)] $B=\gEnd_A(\mathop{\bigoplus}\limits_{i=1}^n Ae_{i}(p_i))$ is an $\mathbb{N}$-graded algebra for some integers $ p_1, p_2,..., p_n$, such that $B_0$ is semisimple and $\ASreg(B)=0$.
     \end{itemize}
\end{theorem}
\begin{proof}
    (1) $\Rightarrow$ (7) It follows from (1) $\Rightarrow$ (2) of Theorem \ref{AS regular 11}.

    (7) $\Rightarrow$ (1) By Theorem \ref{AS regular 11}, it suffices to show that $B$ is a CM-algebra of dimension $d$. Let $R$ be the balanced dualizing complex of $B$.
    Suppose first $\pdim {}_BR <\infty$. Let $P^{\bullet}$ be a minimal graded projective resolution of ${}_BR \in \mathbf{D}^b(\gr B)$. Then, 
 $ \pdim {}_BR = -\inf\{ r \in \mathbb{Z} \mid P^r\neq 0\} \geqslant -\inf {}_BR.$
 
  Since $R\Gamma_B(R)\cong \D(B)$, $\depth {}_BR =0$ and $\depth B =\pdim {}_BR \geqslant -\inf {}_BR$. On the other hand,
  $\depth B = \inf \{ r \in \mathbb{Z} \mid R^r\Gamma_B(B)\neq 0\} =- \sup {}_BR.$
  It follows that $\sup {}_BR \leqslant \inf {}_BR$, and so $B$ is a CM-algebra. 

It is left to show that $\pdim {}_BR <\infty$. The proof is the same as the \cite[proof of Theorem 0.8]{KWZ1}.
\end{proof}

\begin{corollary}\label{coro ASreg=0 2}
   Suppose that $A$ satisfies Hypothesis \ref{Hy}, $A$ is basic, ring-indecomposable and $A_0$ is semisimple. Then the following are equivalent.
  \begin{itemize}
      \item [(1)] $A$ is $\mathbb{N}$-graded AS-regular of dimension $d$ with average Gorenstein parameter $\ell_{av}^A \in \mathbb{Z}$.
      
      \item [(8)] $B=\gEnd_A(\mathop{\bigoplus}\limits_{i=1}^n Ae_{i}(p_i))$ is an $\mathbb{N}$-graded algebra for some integers $ p_1, p_2,..., p_n$, such that $B_0$ is semisimple and $\CMreg(X)-\Torreg(X)=c$ is a constant for all $0 \neq X \in \mathbf{D}^b(\gr B)$.
   \item [(9)] $B=\gEnd_A(\mathop{\bigoplus}\limits_{i=1}^n Ae_{i}(p_i))$ is an $\mathbb{N}$-graded algebra for some integers $ p_1, p_2,..., p_n$, such that $B_0$ is semisimple and $\CMreg(M)-\Torreg(M)=c$ is a constant for all $0 \neq M \in \gr B$.
   
    \item [(10)] $B=\gEnd_A(\mathop{\bigoplus}\limits_{i=1}^n Ae_{i}(p_i))$ is an $\mathbb{N}$-graded algebra for some integers $ p_1, p_2,..., p_n$, such that $B_0$ is semisimple and $\exreg(X)+\ideg(X)=-c$ is a constant for all $0 \neq X \in \mathbf{D}^b(\gr B)$.
  \end{itemize}
\end{corollary}
\begin{proof}
   (1) $\Rightarrow$ (8) See (1) $\Rightarrow$ (3) in Corollary \ref{coro ASreg=0}. 

    (8) $\Rightarrow$ (9) Obviously.

    (9) $\Rightarrow$ (1) It follows that $c=\CMreg(B)$ by taking $M=B$. Then, by taking $M=B_0$, $\ASreg(B)=0$ follows. Hence by Theorem \ref{AS regular 22}, $A$ is $\mathbb{N}$-graded AS-regular of dimension $d$ with average Gorenstein parameters $\ell_{av}^A \in \mathbb{Z}$. 

    (1) $\Rightarrow$ (10) See (1) $\Rightarrow$ (5) in Corollary \ref{coro ASreg=0}. 

     (10) $\Rightarrow$ (1) By taking  $M=B_0$, then $-c=\exreg(B_0)=\Extreg(B_0)$.
    It follows from Corollary \ref{ASreg=0 2} that $\ASreg(B)=0$. Hence by Theorem \ref{AS regular 22}, $A$ is $\mathbb{N}$-graded AS-regular of dimension $d$ with average Gorenstein parameters $\ell_{av}^A \in \mathbb{Z}$.
\end{proof}
The following theorem is a direct extension of \cite[Theorem 0.8]{KWZ1} for $\mathbb{N}$-graded algebras. By Lemma \ref{CM+Torreg=0}, $\exreg(A)+\Exreg(A)=0$ means that the Gorenstein parameters of $A$ are the same for any $\mathbb{N}$-graded AS-Gorenstein algebra $A$.
\begin{theorem}\label{AS regular 33} 
   Suppose that $A$ satisfies Hypothesis \ref{Hy} and $A_0$ is semisimple. Then the following are equivalent.
   \begin{itemize}
       \item [(1)] $A$ is $\mathbb{N}$-graded AS-regular of dimension $d$, and $\Exreg(A)=-\exreg(A)$.
       \item [(2)] $A$ is a graded CM-algebra of dimension $d$ such that $\ASreg(A)=0$.
       \item [(3)] $\ASreg(A)=0$.
   \end{itemize}
    If, furthermore, $A$ is basic, then $\ASreg(A)=0$ if and only if $A$ is AS-regular over $A_0$.
\end{theorem}

\begin{proof}
  (1) $\Rightarrow$ (2) Suppose $A$ is $\mathbb{N}$-graded AS-regular of dimension $d$ with Gorenstein parameters $\{\ell_1,\ell_2,\cdots, \ell_n\}$. Then 
  \begin{align*}
    \exreg(A)=&-d -\ideg(\gExt^d_{A}(A_0,A))\\
    =& -d-\ideg (\mathop{\bigoplus}\limits_{i=1}^n (e_{\sigma (i)}A_0(\ell_i))^{r_i}) \quad (\textrm {by } \eqref{def-GAS-Goren})\\
    =& -d+\max\{\ell_i\}, \,\, \textrm{ and }\\
   \Exreg(A)=&d+\sdeg(\gExt^d_{A}(A_0,A))=d-\min\{\ell_i\}.
  \end{align*}
It follows from $\Exreg(A)=-\exreg(A)$ that $\ell := \max\{\ell_i\}=\min\{\ell_i\}$. By Theorem \ref{ex}, $\CMreg(A)=\Exreg(A)=d-\ell$. By Lemma \ref{GAS-Goren}, $ R\cong \mathop{\bigoplus}\limits_{i=1}^n (Ae_{\sigma(i)}(-\ell))^{r_i}[d]$ and $\Torreg(R)=\ell-d$. It follows from Remark \ref{Torreg(R)=Torreg(A_0)} that $\Torreg(A_0)=\ell-d$. Hence $\ASreg(A)=0$.

(2) $\Rightarrow$ (3) Obviously.

(3)$\Rightarrow$ (1) It follows from  the proof of (7) $\Rightarrow$ (1) in Theorem \ref{AS regular 22} that $A$ is graded CM-algebra of dimension $d$. By the proof of (2) $\Rightarrow$ (1) in Theorem \ref{AS regular 11}, $A$ is $\mathbb{N}$-graded AS-regular of dimension $d$ with Gorenstein parameters $\{\ell,\ell,\cdots,\ell\}$. Hence 
$$\Extreg(A)=-d+\ell \,\, \textrm{ and }\, \Exreg(A)=d-\ell.$$
Therefore, $\exreg(A)+\Exreg(A)=0$.

If $A$ is basic, It follows from Definition \ref{def-AS-Goren over A_0}  that $\ASreg(A)=0$ if and only if $A$ is AS-regular over $A_0.$
\end{proof}

\begin{corollary}\label{coro AS-reg=0}
    Suppose that $A$ satisfies Hypothesis \ref{Hy} and $A_0$ is semisimple. Then the following are equivalent.
    \begin{itemize}
        \item [(1)] $A$ is $\mathbb{N}$-graded AS-regular of dimension $d$ such that $\Exreg(A)=-\exreg(A)$.
        
        \item [(4)] $A$ is a graded CM-algebra of dimension $d$ such that $\CMreg(X)-\Torreg(X)$ is a constant $c$ for all $0\neq X \in \mathbf{D}^b(\gr A)$.
        
        \item [(5)] $A$ is a graded CM-algebra of dimension $d$ such that $\CMreg(M)-\Torreg(M)$ is a constant $c$ for all $0\neq M \in \gr A$.
        
        \item [(6)] $A$ is a graded CM-algebra of dimension $d$ such that $\exreg(X)+\ideg(X)$ is a constant $-c$ for all $0\neq X \in \mathbf{D}^b(\gr A)$.
        
        \item [(7)] $A$ is a graded CM-algebra of dimension $d$ such that $\exreg(M)+\ideg(M)$ is a constant $-c$ for all $0\neq M \in \gr A$.
    \end{itemize}
    If one of the above conditions holds, then $c=\CMreg(A)$. 
\end{corollary}

\begin{proof}
  (1) $\Rightarrow$ (4) It follows from Theorem \ref{AS regular 33} that $\ASreg(A)=0$. By Corollary \ref{ASreg=0}, $\CMreg(X)-\Torreg(X)=\CMreg(A)=c$ is a constant for all $0\neq X \in \mathbf{D}^b(\gr A)$.

  (4) $\Rightarrow$ (5) Obviously.

  (5) $\Rightarrow$ (1) It follows that $c=\CMreg(A)$ by taking $M=A$.  By taking $M=A_0$, it follows that $\ASreg(A)=0$. Hence by Theorem \ref{AS regular 33}, $A$ is $\mathbb{N}$-graded AS-regular such that $\exreg(A)+\Exreg(A)=0$.

(1) $\Rightarrow$ (6) By Theorem \ref{AS regular 33}, $\ASreg(A)=0$. It follows from Corollary \ref{ASreg=0 2} that $\exreg(X)+\ideg(X)=\Extreg(A_0)=-c$ is a constant for all $0\neq X \in \mathbf{D}^b(\gr A)$, where $c:=-\Extreg(A_0)=\CMreg(A)$. 

    (6) $\Rightarrow$ (7) Obviously.

    (7) $\Rightarrow$ (1) 
     By taking  $M=A_0$, it follows that $-c=\exreg(A_0)=\Extreg(A_0)$. Suppose $R=\omega[d]$, where $\omega$ is the balanced CM-module. Let $M=\omega$. Then 
    $$\exreg(\omega)=-\ideg(\omega)+\Extreg(A_0),$$
    and 
    $$\exreg(\omega[d])=-\ideg(\omega[d])+\Extreg(A_0).$$
    It follows from $\CMreg(A)=-\ideg(R)$ that 
    $$\ASreg(A)=\CMreg(A)+\Extreg(A_0)=-\ideg(R)+\Extreg(A_0)=\exreg(R)=0.$$
    By Theorem \ref{AS regular 33}, $A$ is $\mathbb{N}$-graded AS-regular with $\exreg(A)+\Exreg(A)=0$.
\end{proof}
\begin{corollary}\label{coro AS-reg=0 2}
Suppose that $A$ satisfies Hypothesis \ref{Hy} and $A_0$ is semisimple. Then the following are equivalent.
    \begin{itemize}
        \item [(1)] $A$ is $\mathbb{N}$-graded AS-regular of dimension $d$ such that $\Exreg(A)=-\exreg(A)$.
         \item [(8)] $\CMreg(X)-\Torreg(X)$ is a constant $c$ for all $0\neq X \in \mathbf{D}^b(\gr A)$.
        \item [(9)] $\CMreg(M)-\Torreg(M)$ is a constant  $c$ for all $0\neq M \in \gr A$.
        \item [(10)] $\exreg(X)+\ideg(X)$ is a constant $-c$ for all $0\neq X \in \mathbf{D}^b(\gr A)$.
        \end{itemize}
    In fact, $c=\CMreg(A)$ if one of the above holds.
\end{corollary}
\begin{proof}
   (1) $\Rightarrow$ (8) See the proof of (1) $\Rightarrow$ (4) in Corollary \ref{coro AS-reg=0}. 

    (8) $\Rightarrow$ (9) Obviously.

    (9) $\Rightarrow$ (1) By taking $M=A$ and $M=A_0$ respectively, it follows that $c=\CMreg(A)$ and $\ASreg(A)=0$. Hence by Theorem \ref{AS regular 33}, $A$ is $\mathbb{N}$-graded AS-regular such that $\exreg(A)+\Exreg(A)=0$. 

    (1) $\Rightarrow$ (10). See the proof of (1) $\Rightarrow$ (6) in Corollary \ref{coro AS-reg=0}. 

     (10) $\Rightarrow$ (1) By taking $M=A_0$, then $-c= \exreg (A_0)=\Extreg(A_0)$. Then, for all $0\neq X \in \mathbf{D}^b(\gr A)$, 
     $$\exreg(X)=-\ideg(X)+\Extreg(A_0).$$
     It follows from Corollary \ref{ASreg=0 2} that $\ASreg(A)=0$.  Hence by Theorem \ref{AS regular 33}, $A$ is $\mathbb{N}$-graded AS-regular such that $\exreg(A)+\Exreg(A)=0$.
\end{proof}

\section*{Acknowledgment} The authors are very grateful to  Xuyang Chen, Mengying Hu, Yimin Huang and Haonan Li for their comments to the first version of the paper, in particular, for pointing out the condition that $B_0$ is semisimple is necessary in Theorem \ref{AS regular 11}, and helping us modify a preparatory result Theorem \ref{AGP}.

\end{document}